\documentclass[10pt]{amsart}
\usepackage[centertags]{amsmath}
\usepackage{amsfonts}
\usepackage{newlfont}
\usepackage[all]{xy}
\usepackage{graphicx}
\usepackage[latin1]{inputenc}
\usepackage[english]{babel}
\usepackage[T1]{fontenc}
\usepackage{fancyhdr}
\usepackage{amscd}
\usepackage{amsthm}
\usepackage{amssymb}
\usepackage{verbatim}
\usepackage[a4paper,top=4.2cm,bottom=3.8cm,left=3.1cm,right=3.1cm]{geometry}

\newlength{\defbaselineskip}
\newcommand{\setlinespacing}[1]%
           {\setlength{\baselineskip}{#1 \defbaselineskip}}

\newtheorem{thm}{Theorem}[subsection]
\newtheorem{prop}[thm]{Proposition}

\newtheorem*{corintro}{Corollary}

\newtheorem{lem}[thm]{Lemma}
\newtheorem{cor}[thm]{Corollary}
\theoremstyle{definition}
\newtheorem{defn}[thm]{Definition}
\newtheorem{ex}[thm]{Example}
\theoremstyle{remark}
\newtheorem{rem}[thm]{Remark}

\date{}

\DeclareMathOperator{\Hom}{Hom} \DeclareMathOperator{\Ext}{Ext}

\newcommand \an{\alpha_{p^n}}

\newcommand \F{\mathbb{F}}

\newcommand \g{\mathcal{G}}
\newcommand \glb{\g^{(\lambda)}}
\newcommand \gmu{\g^{(\mu)}}
\newcommand \gmup{\g^{(\mu^p)}}

\newcommand \glbn{\g^{(\lambda^{p^n})}}
\newcommand\gln{G_{\lambda,n}}
\newcommand \glx[1]{G_{{\lambda},{#1}}}
\newcommand \gmx[1]{G_{{\mu},{#1}}}
\newcommand \gm{\mathbb{G}_m}
\newcommand \Ga{\mathbb{G}_a}

\newcommand \hxg[1]{H^{1}(X,#1)}
\newcommand \hg[2]{H^1(#1,#2)}

\newcommand \lb{\lambda}

\newcommand \sh{\mathcal{L}}

\newcommand \mup{\mu_p}
\newcommand \mun{\mu_{p^{n}}}

\newcommand \N{\mathbb{N}}

\newcommand \oo[1]{\mathcal{O}_{#1}}
\newcommand \ox{\mathcal{O}_{X}}

\newcommand \Z{\mathbb{Z}}

\newcommand \ha{\mbox{$\hookrightarrow$}}
\newcommand \id{\operatorname{id}}

\newcommand \In{\subseteq}

\newcommand \too{\longrightarrow}

\newcommand \mTo {\longmapsto}

\newcommand \nequiv{\not\equiv}

\newcommand \on{\stackrel}
\renewcommand \phi{\mbox{$\varphi$}}

\newcommand \pt   {\otimes}

\renewcommand \rho{\mbox{$\varrho$}}

\def\clA{{\cal A}}  \def\clD{{\cal D}} \def\clE{{\cal E}}
 \def\clG{{\cal G}} \def\clH{{\cal H}} 
 \def\clL{{\cal L}}

\newcommand{\RW[2]}{W_{#1}{\Omega^{#2}_{X}}}

\newcommand \fr{\operatorname{F}}

\def\clA{{\cal A}}  \def\clD{{\cal D}} \def\clE{{\cal E}}
 \def\clG{{\cal G}} \def\clH{{\cal H}} 
 \def\clL{{\cal L}}

\DeclareMathOperator \Pic {Pic} 
\DeclareMathOperator{\Sp}{Spec}

\begin{document}


\title[Effective models and extension of torsors]{Effective models and extension of torsors over a discrete valuation ring of unequal characteristic}

\author{Dajano Tossici}

\address{Dipartimento di Matematica, Universit\`{a} Roma
Tre, Largo S. Leonardo Murialdo 1, 00146, Rome, Italy\\
 E-mail: dajano@mat.uniroma3.it}

\begin{abstract}

Let $R$ be a discrete valuation ring of unequal characteristic
with fraction field $K$ which contains a primitive $p^2$-{th} root
of unity. Let $X$ be a faithfully flat $R$-scheme
 and $G$ be a finite abstract group.  Let  us consider  a $G$-torsor  $Y_K\too X_K$  and
 let $Y$ be  the normalization of $X_K$ in $Y$. If  $G=\Z/p^n\Z$, $n\le 2$, under some hypothesis on
 $X$,
  we attach some invariants to $Y_K\too X_K$. If $p>2$,  we determine, through these
  invariants,
   when $Y\too X$ has a structure of torsor which extends that of
  $Y_K\too X_K$.
  Moreover,  we explicitly calculate  the  effective model (recently defined by Romagny) of the action of $G$ on $Y$.  
\end{abstract}
\keywords{finite group schemes, torsors, unequal characteristic}
 \maketitle

 \tableofcontents
\section*{Introduction}\label{sec:presentation problem}
\textsc{Notation and conventions.} Throughout the paper, except in
\S \ref{sec:effective models}, we denote by $R$ a discrete
valuation ring (d.v.r. in the sequel) of unequal characteristic,
i.e.  a discrete valuation ring  with fraction field $K$ of
characteristic zero and residue field $k$ of characteristic $p>0$.
Moreover, we write $S=\Sp(R)$. If, for $n\in \N$, there exists a
distinguished  primitive $p^n$-th root of unity $\zeta_n$ in  $R$,
 we write $\lambda_{(n)}:=\zeta_n-1$. Moreover, for any $i\le n$,
 we suppose $\zeta_{i-1}=\zeta_i^p$.  We remark
 that $v(\lb_{(n-1)})=p v(\lb_{(n)})$ and $ v(p)=p^{n-1}(p-1)v(\lb_{(n)})$.
 We will denote by
$\pi\in R$ one of its uniformizers. Moreover if $G$ is an affine
$R$-group scheme  we will denote by $R[G]$ the associated Hopf
algebra. All the schemes will be assumed no\oe therian, however,
sometimes we will explicitly stress this fact. If not otherwise
specified, the cohomology is calculated in the fppf topology.

\vspace{1cm}

Let $G$ be a finite group. It is known that there is a nice smooth
proper stack classifying admissible \mbox{$G$-Galois} covers of
stable curves over $K$, with fixed ramification invariants. While
the problem of understanding its reduction $\mod \pi$  seems quite
far from being solved.  The first phenomenon which may occur is
the following. Let us consider a generically smooth stable marked
curve $Y\to S$ with an action of a finite group $G$ with order
divisible by $p$. Denote $X = Y/G$ and let us assume that $G$ acts
freely on $Y$ outside the set of marked points. It could happen
that the reduced action of $G$ on the special fibre is not
faithful. An  attempt to solve this kind of problems is the
introduction of the notions of effective models (by Romagny
\cite{Ro}) and Raynaud's group schemes (by Abramovich
{\cite{abramovich-2003}). They are very similar notions. In this
paper we will utilize the first one.  More precisely we will
investigate in detail on effective models of actions of cyclic
groups of order $p$ and $p^2$.

Another motivation for this work, related in some sense to the
first one, is the construction of a Hurwitz space for
automorphisms of order $p^2$ of the formal disc $\Sp(R[[T]])$. We
recall that for automorphisms of order $p$ it has been done by
Henrio (\cite{H1}). One of the main tools that Henrio used
 is the strong extension (see
later in the introduction for the definition) of $\Z/p\Z$-torsors
over the boundary of the formal disc $\Sp(R[[T]]\{T^{-1}\})$ where
$R[[T]]\{T^{-1}\}:=\{\sum_{i \in \Z}a_{i}T^{i}\text{ such that
}\lim_{i\to -\infty}a_{i}=0\}.$ We will speak again about it at
the end of the introduction.

 We now introduce more
precisely the problem.
%
%
%
  Let $G$ be an abstract finite group,  $X$  a faithfully flat scheme over $R$ and $Y_K\too X_K$ a
$G_K$-torsor. We remark that, since $K$ is of characteristic $0$,
any finite  group scheme is étale; so, up to an extension of $R$,
any group scheme over $K$ is an abstract group. We call $Y$ the
normalization of $X$ in $Y_K$.
A natural question is the following.

\textsc{Coarse question:}\noindent\textsl{ is it  possible to find
a model $\clG$ of $G_K$ over $R$ together with an action on $Y$
such that $Y\too X$ is a $\clG$-torsor and the action of $G$
coincides with that of $G_K$ on the generic fiber?}

If it happens we say that $Y_K\too X_K$ can be \textsl{strongly
extended}. We observe that the $G$-action on $Y_K$ can be extended
to a $G$-action $\mu:G\times_R Y\too Y$.
The action could be not faithful on the special fiber. The
effective model, as we just said, solve this problem.
An effective model 
for $\mu$ is a flat $R$-group scheme, dominated by $G$, with an
action on $Y$ compatible with $\mu$, such that the action is also
faithful. We will recall its precise definition in
\S\ref{sec:effective models}. When it exists, it is unique. This
implies that if there exists $\clG$ as above, then $\clG$ must be
the effective model of $\mu$. So the above question can be
reformulated in the following way.

 \textsc{Question:}
\textsl{ Which is the effective model $\clG$ (if it exists)  for
the $G$-action? When is $Y\too X$  a $\clG$-torsor?}

%

Let us  suppose that  $(|G|,p)=1$. It is classical that, if $X$ is
regular with geometrically integral  fibers then, up to an
extension of $R$, any connected $G$-torsor can be strongly
extended. This follows from the Theorem of Purity of
Zariski-Nagata (\cite[X 3.1]{SGA1}) and from the Lemma of
Abhyankar (\cite[X 3.6]{SGA1}). See \ref{cor:tame extension} for
the same  result when $X$ is not necessarily regular and $G$ is
abelian.

 Let us now
consider the  case when $p| |G|$. The first case is $G=\Z/p\Z$.
For this group
  the effective models of  its actions have been calculated in some
cases.
 For details see the papers of Raynaud (\cite[1.2.1]{Ra2}), when $X=\Sp(R)$ and $R$ complete, of Green-Matignon (\cite[III
1.1]{GM2}), when $X$ is the $p$-adic closed disc, of Henrio
(\cite[1.6]{H1}), for factorial affine $R$-curves complete with
respect to the $\pi$-adic topology, and of Sa\"idi
(\cite[2.4]{Sa}) for formal smooth curves of finite type. We
remark that the above results  are true under the further
assumption that $X_k$ and $Y_k$ are integral. In all the above
cases, the effective model induces a structure of torsor, i.e. the
$\Z/p\Z$-torsor $Y_K\too X_K$ is strongly extendible. In the
affine case, we will extend these results (with weaker hypothesis
on $X$) also in higher dimensions and, moreover, we will treat the
case $G=\Z/p^2\Z$.

The paper is  organized as follows. In the first section we recall
the definition of effective models and  the principal results
about them. These results are  taken from \cite{Ro}. In the second
one, we recall results (taken from \cite{io2}) about the
classification of finite and flat group schemes over $R$ which are
isomorphic to $(\Z/p^n\Z)_K$ ($n\le 2$) on the generic fiber. We
will call these group schemes models of $(\Z/p^{n}\Z)_K$ $(n\le
2)$.   These models are the candidates to be effective models for
actions of $\Z/p\Z$ or $\Z/p^2\Z$. In the third section, we
construct a filtration of $\hxg{\mun}$ which will be useful to
find the effective model of $\Z/p\Z$-actions and
$\Z/p^2\Z$-actions. In \S\ref{sec:estensione debole} we give a
weak answer to the coarse question for commutative group schemes.
For any $m\in \N$ and any scheme $Z$, we define $
_mPic(Z):=\ker(Pic(Z)\on{m}{\too} Pic(Z))$. The precise statement
is the following.
\begin{corintro}$\ref{cor:estensione torsori per gruppi
commutativi}$ Let $G$ be an abelian group of order $m$ and let us
suppose that $R$ contains a primitive $m$-th root of unity. Let
$X$ be a normal faithfully flat scheme over $R$ with integral
fibers and $_mPic(X)={ _m}Pic(X_K)$. Let us consider a connected
$G$-torsor $f_K:Y_K\too X_K$ and let $Y$ be the normalization of
$X$ in $Y_K$. Moreover, we assume that $Y_k$ is reduced. Then
there exists a (commutative) $R$-group-scheme $G'$ and a
$G'$-torsor $Y'\too X$ over $R$ which extends $f_K$.

\end{corintro}
The point is that we do not require  $Y'$ to coincide with $Y$,
i.e. we do not require $Y'$ to be  normal.  In such a case we
speak about \textsl{weak extension.} Clearly, strong extension
implies weak
extension. 

 In \S \ref{sec:deg Z/pZ-torsori} we study the strong extension
of $\Z/p\Z$-torsors.
 Let us suppose that $R$ contains  a primitive ${p}$-{th} root of
unity. We now suppose that  $X=\Sp(A)$ is a normal faithfully flat
$R$-scheme with integral fibers such that 
 $\pi\in \mathcal{R}_A$, where $\mathcal{R}_A$ is the Jacobson
radical of $A$, and $_pPic(X_K)=0$ (e.g. $A$ a local regular
faithfully flat $R$-algebra with integral fibers).
We will prove in \ref{prop:degenerazione Z/pZ torsori} that any
connected \mbox{$\Z/p\Z$-torsor}  $Y_K\too X_K$ is strongly
extendible under the assumption that the special fiber of the
normalization $Y$  of $X$ in $Y_K$ is reduced. 

 In the sixth section,
 we consider the case
$G=\Z/p^2\Z$. We will assume that $R$ contains a primitive
$p^2$-th root of unity. Let $X:=\Sp A$ be a normal essentially
semireflexive scheme
 over $R$ (see \S \ref{sec:effective models}) with integral fibers  such that $\pi \in \mathcal{R}_A$. We
moreover assume \mbox{$_{p^2} Pic(X_K)=0$.}  Let $Y_K\too X_K$ a
connected $\Z/p^2\Z$-torsor, $Y$ be the normalization of $X_K$ in
$Y$, and  assume that  $Y_k$ integral.  We will attach to any such
$\Z/p^2\Z$-torsor an element  of $\N^4$. Our main result is that,
for $p>2$, this element determines  its effective model $\clG$,
which we explicitly describe (see \ref{teoremone}). Moreover, we
will give a criterion to see if $Y$ is a $\clG$-torsor (see
\ref{cor:criterio torsore}).
 Finally, we will give an
example of $\Z/p^2\Z$-torsor over $X_K$ satisfying  the above
hypothesis but non-strongly extendible.

The last section is devoted to study elements of $\N^4$ which
arise from a $\Z/p^2\Z$-torsor over $X_K$ as above. The main
result of this section is the theorem \ref{teo:realizzazione tipi
di degenerazione} about the classification of $4$-uples of natural
integers which are associated to
strongly extendible $\Z/p^2\Z$-torsors over $X_K$. 

 We recall that the case  $R$ of positive characteristic
has been studied by Sa\"idi in \cite{sa2}. He proved a result of
strong extension of $\Z/p\Z$-torsors for formal normal schemes of
finite type of any dimension (\cite[2.2.1]{sa2}) and moreover he
studied the case of $\Z/p^2\Z$-torsors. His approach is slightly
different: he is not interested in effective models but only in
explicit equations of the induced cover on the special fiber.
Moreover, he gave an example  of non strongly extendible
$\Z/p^2\Z$-torsors in \cite{sa3}. Another such example has been
given by Romagny (\cite[2.2.2]{Ro}). We remark that in equal
characteristic
 there is no criterion   to determine if a $\Z/p^2\Z$-torsor is strongly
 extendible.

To study  extensions of  $\Z/p^n\Z$-torsors, with $n\ge 3$, it is
necessary to have  an explicit classification of models of
$(\Z/p^n \Z)_K$ as done in \cite{io2} for models of
$(\Z/p^2\Z)_K$. The strategy used in \cite{io2} could be used, in
principle,  also for models of  $(\Z/p^n \Z)_K$. However, this
could lead to very complicated calculations.

We finally remark about the construction of a Hurwitz space for
automorphisms of order $p^2$ mentioned earlier, that, since the
strong extension does not work for $\Z/p^2\Z$-torsors, it seems
reasonable to restrict to automorphisms of order $p^2$ which
induce, on the boundary of the formal disc, strongly extendible
covers.

\section*{Acknowledgments.} This paper constitutes part of my PhD
thesis. It is a pleasure to thank my advisor Professor Carlo
Gasbarri for  guidance and  constant encouragement.  Part of this
work was done during my visit to
 Professor Michel Matignon in 2005 at the Department of
Mathematics of the University of Bordeaux 1. I am indebted to him
for the useful conversations  and for his great interest in my
work. I am deeply grateful  to Matthieu Romagny  for having
pointed out me his work \cite{Ro} and, above all, for his very
careful reading of this paper and for his several comments,
suggestions and remarks  which helped me to improve an earlier
version of the paper. Finally I thank
Marco Antei for useful discussions. 

\section{Effective models}\label{sec:effective models}
Here we recall some definitions and results about effective models
which will play a key role in our results about extensions of
$\Z/p^2\Z$-torsors. For more details see \cite{Ro}, from which
most of the material of this paragraph has been taken. In this
section $R$, is  a d.v.r. not necessarily of unequal
characteristic.
We firstly recall the definition of model of a group scheme.
\begin{defn}
Let $H_K$ be a group scheme over $K$. Any  flat $R$-group scheme
$G$ such that $G_K\simeq H_K$ is called a \textsl{model} of $H_K$.
If $H_K$ is finite over $K$ we also require that $G$ is finite
over $R$.
\end{defn}
 In
\cite[\S~2]{Ra1}, the relation of \emph{domination} between models
of a group scheme has been introduced.
\begin{defn}
 Let $G_1$ and $G_2$ be finite flat group schemes over $R$
with an isomorphism \mbox{$u_K\colon G_{1,K}\too G_{2,K}$.}  We
say that $G_1$ \emph{dominates} $G_2$ if we are given an
$R$-morphism $u\colon G_1\too G_2$ which restricts to $u_K$ on the
generic fibre. The map $u$ is also called a \textsl{model map}. If
moreover we are given two actions $\mu_i:G_i\times_R Y\too Y$, for
$i=1,2$,  we say that $G_1$ dominates $G_2$ \emph{compatibly}
(with the actions) if $\mu_1=\mu_2\circ (u\times\id)$.
\end{defn}
We now recall the definition of a faithful action.
\begin{defn} Let $G$ be a group scheme   which acts on a scheme $Y$ over a scheme $T$. This action is
\textsl{faithful} if the induced morphism of sheaves of groups, in
the fppf topology of $T$,
$$
G\too \clA ut_T(Y)
$$
is injective.
\end{defn}
%
%
We recall here the definition of effective model given by Romagny.
\begin{defn}
 Let $G$ be a finite flat group scheme over $R$. Let $Y$ be a faithfully flat
scheme over $R$. Let $\mu\colon G\times_R Y\too Y$ be an action,
faithful on the generic fibre. An \emph{effective model} for $\mu$
is a finite flat $R$-group scheme $\clG$ acting on $Y$, dominated
by $G$ compatibly, such that $\clG$ acts  faithfully on $Y$.
\end{defn}
\begin{ex}\label{ex:torsore e modello effettivo}
Let $X$ be a faithfully flat scheme over $R$ and $\clG$ a finite
and flat group scheme over $R$. Let  $Y\too X$ be a $\clG$-torsor
over $R$. Then $\clG$ is already an effective model. Indeed let us
suppose that $ \clG\too \clA ut_R(Y)$ is not injective. Then there
exists
a faithfully flat morphism $U\too \Sp(R)$ 
and $g\in \clG(U)\setminus \{0\}$ such that
$$
Y\times_R U\on{g\times id}{\too} \clG\times_R Y\times_R
U\on{\mu\times\id}{\too} Y\times_R Y\times_R U
$$
is  equal to $\Delta\times \id:Y\times_R U\too Y\times_RY\times_R U$
where $\Delta:Y\too Y\times_RY$ is the diagonal morphism. By the
definition of $\clG$-torsor $\clG\times_R Y\times_R
U\on{\mu\times\id}{\too} Y\times_R Y\times_R U$ is an isomorphism.
Then
$$
Y\times_R U\on{g\times id}{\too} \clG\times_R Y\times_R
U\on{\mu\times\id}{\too} Y\times_R Y\times_R
U\on{(\mu\times\id)^{-1}}\too \clG\times_R Y\times_R U
$$
is the zero section, against assumptions.
\end{ex}
 We report here some results about
effective models.
\begin{prop}\label{prop:unicita modelli effettivi}
An effective model is unique up to unique isomorphism, if it exists.
\end{prop}
\begin{proof}
\cite[1.1.2]{Ro}.
\end{proof}
The following crucial remark is the reason for our interest in
effective models.
\begin{rem}\label{rem:modelli effettivi e torsori}
 Let $G$ be a finite and flat group scheme over $R$ and  $Y$  a faithfully flat
scheme  over $R$.  Let $\mu\colon G\times_R Y\too Y$ be an action.
Moreover we suppose that $Y_K\too Y_K/G_K$ is a $G_K$-torsor. Then
by  \ref{prop:unicita modelli effettivi} we have that the
effective model $\clG$ whose action extends that of $G_K$ is
unique if it exists. By \ref{ex:torsore e modello effettivo} this
means that if there exists a  model
 $G'$  of $G_K$, compatible with the action, such that $Y$ is a $G'$-torsor, then $G'$ is the effective model for $\mu$.
\end{rem}
We recall that an action $\mu:G\times_R Y\too Y$ is admissible if
$Y$ can be
covered by $G$-stable open affine subschemes. 
\begin{prop}\label{prop:proprieta modelli effettivi}
Let $G$ be a finite flat group scheme over $R$. Let $Y$ be a
faithfully flat scheme over $R$ and $\mu:G\times_R Y\too Y$ an
admissible action, faithful on the generic fiber. Assume there
exists an effective model $\clG$. Then
\begin{itemize}\item[(i)] If $H$ is a finite flat subgroup of $G$, the restriction of the action to
$H$ has an effective model $\clH$ which is the schematic image of
$H$ in $\clG$. If $H$ is normal in $G$, then $\clH$ is also normal
in $\clG$. \item[(ii)] The identity of $Y$ induces an isomorphism
$Y/G\simeq Y/\clG$. \item[(iii)] Assume that there exists an open
subset $U\In Y$ which is schematically dense in any fiber of
\mbox{$Y\too \Sp(R)$} such that $\clG$ acts freely on $U$. Then
for any closed normal subgroup $H\vartriangleleft G$ the effective
model of $G/H$ acting on $X/H$ is $\clG/\clH$.
\end{itemize}
\end{prop}

\begin{proof}
\cite[1.1.3]{Ro}.
\end{proof}

\begin{thm} \label{teo:modelli effettivi caso algebrico}
Let $G$ be a finite flat group scheme over $R$. Let $Y$ be a
faithfully flat scheme of finite type over $R$ and let $\mu\colon
G\times Y\too Y$ be an action. We assume that $Y$ is covered by
$G$-stable open affines $U_i$ with function ring separated for the
$\pi$-adic topology, such that $G$ acts faithfully on the generic
fibre $U_{i,K}$. Then, if $Y$ has reduced special fibre, there
exists an effective model for the action of $G$.
\end{thm}
\begin{proof}
\cite[1.2.3]{Ro}.
\end{proof}
We remark that the condition about the separatedness of the
function rings of $U_i$ is assured if, for instance, we assume $Y$
no\oe therian and integral. This follows from the Theorem of Krull
(\cite[1.3.13]{liu}).

If we add some hypothesis on $Y$, then we have a useful criterion
to see if a group scheme which acts on $Y$ is the effective model
for the action.

Recall that a module $M$ over a ring $A$ is called
\textsl{semireflexive} if the canonical map from $M$ to its bidual
is injective. Equivalently, $M$  is a submodule of some product
module $A^I$. Indeed, consider the set $I=\Hom_{A}(M,A)$ and the
morphism $a:M\to A^I$ mapping $x$ to the collection of values
$(f(x))_{f\in I}$ for all linear forms $f$. By definition, if $M$
is semireflexive then for each nonzero $x\in M$ there exists a
linear form such that $f(x)\ne 0$, so $a$ is injective. The
converse is easy.
\begin{rem}\label{rem:semireflexive->flat}A semireflexive module over   $R$
is faithfully flat and separated with respect to the $\pi$-adic
topology. Indeed, since $M\In R^I$ for some set $I$, then $M$ is
torsion free, hence flat over $R$. Moreover let $x\in \cap\pi^m
M$. Then for any linear form $f$ we have $f(x)\in \cap\pi^m R=0$.
Since $M$ is semireflexive over $R$, this implies $\cap\pi^m M=0$.
So $M$ is separated with respect to the $\pi$-adic topology. But,
over a d.v.r,  being flat and separated with respect to the
$\pi$-adic topology implies faithfully flat. Indeed, since
$\cap\pi^m M=0$ implies $M\neq \pi M$, then $M$ is faithfully flat
(see \cite[1.2.17]{liu}). We do not know if the converse is true
too, i.e. if any (faithfully) flat $R$-module separated with
respect to the $\pi$-adic topology is semireflexive over $R$.
\end{rem}
\begin{rem}
If $M$ is semireflexive over $R$ and $M'\In M$ is an inclusion of
$R$-modules  then $M'$ is semireflexive over $R$. It easily
follows by definition.
\end{rem}
\begin{ex}  Any free  $R$-module  is semireflexive over $R$, e.g. any reduced faithfully flat $R$-algebra of finite type
with geometrically irreducible fibers (see
\cite[3.3.4,3.3.5]{GR}). Other examples, not free, are
$R[[T_1,\dots,T_n]]$ or
$$
R[[T_1,\dots,T_n]]\{T_1^{-1},\dots,T_n^{-1}\}:=\{\sum_{(i_1,\dots,i_n)
\in \Z^n}a_{i_1\dots i_n}T^{i_1}\dots T^{i_n}\text{ such that
}\lim_{i_1+\dots+i_n\to -\infty}a_{{i_1\dots i_n}}=0\}.
$$
Finally, if  $R$ is complete, $X\too R$ a smooth surjective
morphism of finite type of dimension $n$ and $x\in X$ a closed
point in the special fiber with residue field $k$ then $\oo{X,x}$
is semireflexive over $R$. This follows from the previous remark
since $O_{X,x}\In{\widehat{O}}_{X,x}\simeq {R}[[T_1,\dots,T_n]]$.
\end{ex}
The following lemma will be useful in \S\ref{sec:deg
Z/p^2Z-torsors}.
\begin{lem}\label{lem:finite-->semireflexive} Let $A$ be an $R$-algebra which is
 semireflexive as  an $R$-module and let $B$ be a  flat
$R$-algebra. 
If there exists a finite $R$-morphism of modules
$$
A\too B,
$$
such that $B_K$ is semireflexive as an $A_K$-module,  then $B$ is
semireflexive as an $R$-module.
\end{lem}
\begin{rem}
In particular any finite and flat $R$-algebra is semireflexive as an
$R$-module.
\end{rem}
\begin{proof}
Let us consider  $B_K$. It is a vector space over $K$, so in
particular it is  semireflexive over $K$. 
Since $A$ and $B$ are flat over $R$,  the natural maps $A\too A_K$
and $B\too B_K$ are injective. We now prove that $B$ is
semireflexive over $A$. Let $b\in B$ and let us take an
$A_K$-linear form $f:B_K\too A_K$ such that $f(b)\neq 0$. It
exists since $B_K$ is semireflexive over $A_K$. Let $b_1,\dots,
b_m$ generators of $B$ as an $A$-module and let $n\in\N$ such that
$\pi^n f(b_i)\in A$ for $i=1,\dots,m$. Then we have  $\pi^n
f(B)\In A$. Moreover, $\pi^n f(b)\neq 0$, since $A$ is flat over
$R$ by \ref{rem:semireflexive->flat}. So $\pi^n f:B\too A$ is a
linear form with $\pi^n f(b)\neq 0$. Then $B$ is semireflexive as
an $A$-module. But $A$ is semireflexive over $R$. Therefore, $B$
is semireflexive over $R$. Indeed for any $b\in B$ let us take an
$A$-linear form $g:B\too A$ with $g(b)\neq 0$. Moreover, let us
consider an $R$-linear form $h:A\too R$ such that $h(g(b))\neq 0$.
Then $h\circ g:B\too R$ is an $R$-linear form with $h\circ
g(b)\neq 0$. Hence, $B$ is semireflexive over $R$.
\end{proof}
\begin{defn}
We will say that a morphism of schemes $f:X\to T$ is {\em
essentially semireflexive} if there exists a cover of $T$ by open
affine subschemes $T_i$, an affine faithfully flat $T_i$-scheme
$T'_i$ for all $i$, and a cover of $X'_i=X\times_T T'_i$ by open
affine subschemes $X'_{ij}$, such that the function ring of
$X'_{ij}$ is semireflexive as a module over the function ring of
$T'_i$.
\end{defn}This is a generalization of the definition of an
essentially free morphism given in \cite{SGA3}. The proofs of the
following two lemmas have been suggested to us by Romagny.
\begin{lem} Let $X$ be essentially semireflexive and separated
over $T$. Let $G$ be a $T$-group scheme acting on $X\to T$. Then the
kernel of the action is representable by a closed subscheme of $G$.
\end{lem}
\begin{proof} Proceeding like in \cite{SGA3} we are reduced to
proving the analogue of \cite[6.4]{SGA3}. Then the proof given in
\cite{SGA3} works in our case, because the only property of free
modules that is used in the proof is that they
are semireflexive. 
\end{proof}
The next lemma is the reason  we are interested in essentially
semireflexive schemes. Indeed, in such a case we have an useful
criterion to check if a finite group scheme is an effective model.
\begin{lem}\label{lem:basta fedelta' su fibra speciale} Let $G$ be a  finite and
flat $R$-group scheme which acts on an essentially semireflexive
$R$-scheme. Then the action of $G$ is faithful if and only if the
action of $G_k$ on the special fibre is faithful.
\end{lem}
\begin{proof} Only the {\em if} part needs a proof. Let $I_G$ be
the augmentation ideal of $G$ and let $J$ be the ideal defining
the kernel $H$ of the action. Since $H$ is a subgroup scheme of
$G$ and $H_k$ is trivial, then
\begin{equation}\label{eq:J in IG}J\In I_G \quad \text{ and }I_G+\pi R[G]=J+\pi
R[G]\end{equation} Moreover, since $R[G]/I_G$ is flat over $R$ then
\begin{equation}\label{eq:piJ=piAcap J}I_G\cap \pi R[G]=\pi I_G.
\end{equation}  We now
claim that
$$
I_G=J+\pi I_G.
$$
Clearly $J+\pi I_G\In I_G$. We now prove the converse. Let $a\in
I_G$, then from \eqref{eq:J in IG} it follows that $a=b+\pi c$ for
some $b\in J$ and $c\in R[G]$. Since, $J\In I_G$ then $\pi c\in
I_G\cap \pi R[G]$. Therefore, by \eqref{eq:piJ=piAcap J} we have
$c\in I_G$. Hence, $I_G\In J+\pi I_G$. We have so proved $I_G=
J+\pi I_G$.
 Then $I_G/J$ is an $R$-module of finite type and
$(I_G/J)\otimes_R k=I_G/(J+\pi I_G)=0$, so $I_G/J=0$ by Nakayama's
lemma. Hence, the kernel is trivial.
\end{proof}

\section{Models of $(\Z/p\Z)_K$ and $(\Z/p^2\Z)_K$}\label{sec:modelli}
\subsection{Some group schemes of order $p^n$ and models of $\Z/p \Z$}
 Let $G$ be a constant
finite $R$-group scheme. By definition, the effective model for a
$G$-action over a scheme $Y$ is, in particular, a $G_K$-model
$\clG$ of $G$ with  a morphism $G\too\clG$.  Since $G$ is étale,
for any model $\clG'$ of $G_K$  the existence of a morphism
$G\too\clG'$ is automatic: $\clG'$ is finite over $R$ so you can
take the normalization map. In this section we recall results
about models of $(\Z/p\Z)_K$ and $(\Z/p^2\Z)_K$.
 The results are taken from \cite{io2}. See there for the
proofs and for more details.

We  introduce some $R$-smooth unidimensional group schemes. 
For any $\lambda\in R$ define the group scheme
$$
\glb=\Sp(R[T,\frac{1}{1+\lambda T}])
$$
The  $R$-group scheme  structure is given by
\begin{align*}
&T\longmapsto 1\pt T+T\pt 1 +\lb T\pt T\qquad&\text{comultiplication,}\\
&T\longmapsto 0 \qquad& \text{counit,}\\
&T\longmapsto -\frac{T}{1+\lb T}\qquad& \text{coinverse,}
\end{align*}

We observe that if $\lb=0$ then $\glb\simeq \Ga$. It is possible
to prove that $\glb\simeq \gmu$ if and only if $v(\lambda)=v(\mu)$
and the isomorphism is given by $T\longmapsto
\frac{\lambda}{\mu}T$. Moreover, it is easy to see that, if
$\lb\in \pi R\setminus\{ 0\}$, then $\glb_k\simeq \mathbb{G}_a$
and $\glb_K\simeq \gm$. It has been proved by Waterhouse and
Weisfeiler, in \cite[2.5]{Wat}, that
 any deformation, as a group scheme, of $\Ga$ to
 $\gm$ is isomorphic to $\glb$ for some $\lb\in \pi R\setminus\{0\}$.
If $\lb\in R\setminus\{0\}$, we can define the  morphism
\begin{align*}
\alpha^{\lb}:\glb\too \gm\\
\end{align*}
given, on the level of Hopf algebras,  by  $T\mTo 1+\lb T$: it is
an isomorphism on the generic fiber. If $v(\lb)=0$ then
$\alpha^{\lb}$ is an isomorphism.
For any flat $R$-scheme $X$, let us consider 
 the exact sequence on the fppf site $X_{fl}$
\begin{equation}\label{eq:succ esatta glb in gm}
0\too \glb\on{\alpha^{\lb}}{\too}{\gm}\too
i_*{{\mathbb{G}}_{m,X_{\lb}}}\too 0,
\end{equation}
 where $i$ denotes the closed immersion $X_\lb=X\pt_R(R/\lb R)\ha X$
(see \cite[1.2]{SS2}).
This gives the following  associated long exact sequence
\begin{equation}\label{eq:succ esatta lunga di g lambda}
\begin{aligned} 0\too H^0(X,\glb)&\too H^0(X,{\gm})\too
H^0({X_\lb},{\gm})\too\\
 \ \ \ \ &\ \ \ \ \ \too\hxg{\glb}\too
\hxg{{\gm}}\too \hg{X_\lb}{{\gm}}\too \dots\\
\end{aligned}
\end{equation}

We now define some finite and flat group schemes of order $p^n$. Let
$\lambda\in R$  satisfy the condition
 $$
(*)\qquad v(p)\ge p^{n-1}(p-1)v(\lambda).
$$ Then the map
\begin{align*}
\psi_{\lb,n}:&\glb\too \g^{(\lb^{p^n})}\\
             &T\longmapsto P_{\lb,n}(T):=\frac{(1+\lambda T)^{p^n}-1}{\lambda^{p^n}}
\end{align*}
is an isogeny of  degree $p^n$. Let
$$
\gln:=\Sp(R[T]/P_{\lb,n}(T))
$$
be  its kernel. It is a commutative finite flat group scheme over
$R$ of rank $p^n$. It is possible to prove that
\begin{align*}
{(G_{\lb,n})}_k\simeq  \mup &\qquad\text{ if $v(\lb)=0$;}\\
{(G_{\lb,n})}_k\simeq \an &\qquad\text{ if $p^{n-1}(p-1)v(\lb)<v(p)$};\\
(G_{\lb,n})_k\simeq \alpha_{p^{n-1}}\times \Z/p\Z &\qquad \text{ if
$p^{n-1}(p-1)v(\lb)=v(p)$}.
\end{align*}

  We observe that $\alpha^{\lambda}$ is compatible with $\psi_{\lb,n}$, i.e
 the following  diagram is commutative
\begin{equation}\label{eq:compatibilità phi_lb e alpha^lb}
\xymatrix@1{\glb\ar[d]_{\psi_{\lb,n}}  \ar[r]^{\alpha^{\lb}}&\ar[d]^{p^{n}}\gm\\
\glbn \ar[r]^(.6){\alpha^{\lb^{p^n}}}&\gm}
\end{equation}
Then it   induces a  map
$$
\alpha^{\lb,n}:\gln\too\mun,
$$
which is an isomorphism on the generic fiber. And if $v(\lb)=0$,
then $\alpha^{\lb,n}$ is an isomorphism.
We remark that
\begin{equation}\label{eq:Hom(glbn,glmn)}
 \Hom_{gr}(G_{\lb,n},G_{\lb',n})\simeq%
\begin{array}{ll}
        \Z/p^{n-r}\Z, \\
\end{array}%
\end{equation}
where $r=\min\{0\le r'\le n \text{\ such that\ } p^{r'}v(\lb)\ge
v(\lb')\}$.  The morphisms are given by
\begin{align*}
G_{\lb,n}&\too G_{\lb',n}\\
 T&\longmapsto \frac{(1+\lb T)^{p^{r}i}-1}{\lb'}
\end{align*}
for $i=0,\dots,p^{n-r}-1$. It  follows easily that
$G_{\lb,n}\simeq G_{\lb',n}$ if and only if $v(\lb)=v(\lb')$.

From the exact sequence of group schemes
$$
0\too G_{\lb,n}\on{i}{\too} \glb\on{\psi_{\lb,n}}{\too}\glbn\too 0
$$
we have the exact sequence of groups
\begin{equation}\label{eq:succ esatta glbn}
\dots\too H^0(X,\glb) \on{(\psi_{\lb,n})_*}{\too}H^0({X},{\glbn})
\too \hxg{\gln}\on{i_*}{\too}\hxg{{\glb}}\too \hg{X}{{\glbn}}\too
\dots
\end{equation}
 In the following when we  speak
about $G_{\lb,n}$ it will be assumed  that $\lb$ satisfies $(*)$. If
$R$ contains a primitive $(p^n)$-th root of unity $\zeta_n$ then,
since \mbox{$v(p)=p^{n-1}(p-1)v(\lb_{(n)})$},  the condition $(*)$
is equivalent to $v(\lb)\le v(\lb_{(n)})$.
The following result classifies models of $\Z/p\Z$.
\begin{thm}\label{teo:modelli di Z/pZ}
Let us suppose that $R$ contains a primitive $p$-th root of unity.
If $G$ is a finite and flat $R$-group scheme such that $G_K\simeq
\Z/p\Z$, then $G\simeq G_{\lb,1}$ for some $\lb\in R\setminus\{0\}$.
\end{thm}
\begin{proof}
See \cite[III 3.2.2]{RoPhD}. 
\end{proof}

\subsection{Models of $\Z/p^2\Z$}
We firstly recall the definition of  some smooth $R$-group schemes
of dimension $2$, which are extensions of $\gmu$ by $\glb$, for some
$\mu,\lb\in R\setminus \{0\}$. We remark that it is possible to
prove that any action of $\gmu$ on $\glb$ is trivial (\cite[I
1.6]{SS6}). 
  Let $Y=\Sp(R[T_1,\dots,T_m]/(F_1,\dots, F_n))$ be an affine
$R$-scheme of finite type. We recall that, for any $R$-scheme $X$,
we have that $\Hom_{Sch}(X,Y)$ is in bijective correspondence with
the set
$$\{(a_1,\dots,a_m)\in H^0(X,\oo{X})^m| F_1(a_1,\dots,
a_m)=0,\dots,F_n(a_1,\dots, a_m)=0 \}.$$  With an abuse of notation
we will identify these two sets. If $X$ and $Y$ are $R$-group
schemes
  we will also identify $\Hom_{gr}(X,Y)$ with a subset of $$\{(a_1,\dots,a_m)\in H^0(X,\oo{X})^m| F_1(a_1,\dots,
a_m)=0,\dots,F_n(a_1,\dots, a_m)=0 \}.$$

For any $\lb \in R \setminus\{0\}$ let us define $S_\lb:=\Sp(R/\lb
R)$. We now fix presentations for the group schemes $\gm$ and
$\glb$ with $\lb \in \pi R$. Indeed we write
\mbox{$\gm=\Sp(R[S,1/S])$} and $\glb=\Sp(R[S,1/1+\lb S])$. We
observe that by definition we have that
$$
\Hom_{gr}({\gmu}_{|S_{\lb}},{\gm}_{|S_\lb})=\{F(S)\in (R/\lb
R[S,\frac{1}{1+\mu S}])^*|F(S)F(T)=F(S+T+\mu ST) \}
$$
 We remark that throughout the
paper will be a conflict of notation since  $S$ will denote both
$\Sp(R)$ and an indeterminate. But it should not cause any
problem. If we apply the functor $\Hom_{gr}(\gmu,\cdot)$ to the
sequence \eqref{eq:succ esatta glb in gm}, we obtain, in
particular,  a
 map 
\begin{equation*}
\begin{aligned}
\Hom_{gr}(\gmu_{|S_\lb},{\gm}_{|S_\lb})
\on{\alpha}{\too}\Ext^1(\gmu&,\glb).
\end{aligned}
\end{equation*}
This map is given by
$$
F\mTo \clE^{(\mu,\lb;F)},
$$
where
$$
\clE^{(\mu,\lb;F)}
$$
is a smooth affine  commutative group defined as follows: let
$\tilde{F}(S)\in R[S]$ be a lifting of $F(S)$, then
$$
\clE^{(\mu,\lb;F)}= \Sp(R[S_1,S_2,\frac{1}{1+\mu
S_1},\frac{1}{\tilde{F}(S_1)+\lb S_2}])
$$
\begin{enumerate}
    \item law of multiplication
    \begin{align*}
    S_1\longmapsto &S_1\pt 1+1\pt S_1+\mu S_1\pt S_1\\
    S_2\longmapsto &S_2\pt \tilde{F}(S_1)+\tilde{F}(S_1)\pt S_2 +\lb S_2\pt S_2+\\
                   & \qquad   \quad    \frac{\tilde{F}(S_1)\pt \tilde{F}(S_1)-\tilde{F}(S_1\pt 1+1\pt S_1+\mu S_1\pt S_1)}{\lb}
    \end{align*}
    \item unit
    \begin{align*}
    &S_1\longmapsto 0\\
    &S_2\longmapsto \frac{1-\tilde{F}(0)}{\lb}
    \end{align*}
    \item inverse
    \begin{align*}
    &S_1\longmapsto -\frac{S_1}{1+\mu S_1}\\
& S_2 \longmapsto\frac{\frac{1}{\tilde{F}(S_1)+\lb
S_2}-\tilde{F}(-\frac{S_1}{1+\mu S_1})}{\lb}
\end{align*}
\end{enumerate}
Moreover, we  define the following  homomorphisms of group schemes
$$
\glb=\Sp(R[S,(1+\lb S)^{-1}])\too \clE^{(\mu,\lb;F)}
$$
by
\begin{align*}
&S_1\longmapsto 0\\
&S_2\longmapsto S +\frac{1-\tilde{F}(0)}{\lb}
\end{align*}
and
$$
\clE^{(\mu,\lb;F)}\too \gmu=\Sp(R[S,\frac{1}{1+\mu S}])
$$
 by
\begin{align*}
&S\too S_1.
\end{align*}
It is easy to see that
\begin{equation}\label{eq:estensioni lisce}
0\too \glb\too \clE^{(\mu,\lb;F)}\too \gmu\too 0
\end{equation}
is exact. A different choice of the lifting $\tilde{F}(S)$ gives
an isomorphic extension. We recall the following theorem.
\begin{thm}\label{teo:ss ext1}For any $\lb,\mu\in  R\setminus\{0\}$,  $$\alpha:\Hom_{gr}(\gmu_{|S_\lb},{\gm}_{|S_\lb}){\too}
\Ext^1(\gmu,\glb)$$  is a surjective morphism of groups. And
$\ker(\alpha)$ is generated by the class of $1+\mu S$. In particular
any extension of $\gmu$ by $\glb$ is commutative.
\end{thm}

\begin{proof}
See \cite[\S 3]{SS9} or \cite[II 1.2]{SS6} for the case $\mu,\lb\in
\pi R\setminus\{0\}$ and \cite[I 2.7, II 1.4]{SS6} for the others.
\end{proof}
\begin{ex}If $v(\mu)=0$ or $v(\lb)=0$ we have  that
$\Ext^1(\gmu,\glb)=0$.
\end{ex}
We are interested in the study of the models of  $(\Z/p^2\Z)_K$.
Let us suppose that $R$ contains a primitive $p^2$-th root of
unity. So $(\Z/p^2\Z)_K\simeq (\mu_{p^2})_K$. By \ref{teo:modelli
di Z/pZ}, it is easy to show that any such model is an extension
of $G_{\mu,1}$ by $G_{\lb,1}$ for some $\mu,\lb\in
R\setminus\{0\}$. So the first step is to study the group
$\Ext^1(G_{\mu,1},G_{\lb,1})$. We remark that it is possible to
show that any action of $\gmx{1}$ on $\glx{1}$ is trivial. We now
define some extensions of $G_{\mu,1}$ by $G_{\lb,1}$.
%

\begin{defn}\label{def:clE(mu,lb,F,i)}
 Let $F\in \Hom_{gr}({\gmx{1}}_{|S_\lb},{\gm}_{|S_\lb})$, $j\in
\Z/p\Z$ such that  $$F(S)^p(1+\mu S)^{-j}=1\in
\Hom_{gr}({\gmx{1}}_{|S_{\lb^p}},{\gm}_{|S_{{\lb^p}}}).$$ Let $\tilde{F}(S) \in R[S]$  be a  lifting of $F$. 
We denote by $\clE^{(\mu,\lb;F,j)}$  the subgroup scheme of
$\clE^{(\mu,\lb;F)}$ given
 on the level of schemes  by
$$
\clE^{(\mu,\lb;F,j)}=\Sp\bigg(R[S_1,S_2]/\big(\frac{(1+\mu
S_1)^p-1}{\mu^p},\frac{(\tilde{F}(S_1)+\lb S_2)^p(1+\mu
S_1)^{-j}-1}{\lb^p}\big)\bigg).
$$
\end{defn}

 We moreover  define the following
homomorphisms of group schemes
$$
\glx{1}\too \clE^{(\mu,\lb;F,j)}
$$
by
\begin{align*}
&S_1\longmapsto 0\\
&S_2\longmapsto S +\frac{1-\tilde{F}(0)}{\lb}
\end{align*}
and
$$
\clE^{(\mu,\lb;F,j)}\too G_{\mu,1}
$$
 by
\begin{align*}
&S\mTo S_1.
\end{align*}
It is easy to see that
$$
0\too \glx{1}\too \clE^{(\mu,\lb;F,j)}\too G_{\mu,1}\too 0
$$
is exact. A different choice of the lifting $\tilde{F}(S)$ gives an
isomorphic extension. It is easy to see that, as a group schemes,
$(\clE^{(\mu,\lb;F,j)})_K\simeq (\Z/p^2\Z)_K $ if  $j\neq 0$ and
$(\clE^{(\mu,\lb;F,0)})_K\simeq (\Z/p\Z\times \Z/p\Z)_K $.
\begin{rem}In the above definition the integer $j$ is uniquely
determined by $F\in \Hom_{gr}({\gmx{1}}_{|S_\lb},{\gm}_{|S_\lb})$
if and only if ${\lb^p} \nmid \mu$.
\end{rem}
We recall the following result.

\begin{prop}\label{lem:suriettività mappa tra hom}
Let    $\lb,\mu \in R\setminus\{0\}$. 
The restriction map
$$
i^*:\Hom_{gr}(\gmu_{|S_\lb},{\gm}_{|S_\lb})\too
\Hom_{gr}({G_{\mu,1}}_{|S_\lb},{\gm}_{|S_\lb})
$$
induced by
$$
i:G_{\mu,1}\hookrightarrow \gmu
$$
is surjective. If $p>2$,
 $v(p)\ge(p-1)v(\mu)>0$ and $v(\mu)\ge
v(\lb)>0$, then we have an isomorphism of groups
$$
(\xi^0_{R/\lb R})_p:(R/\lb R)^{\fr}
\too\Hom_{gr}({\gmx{1}}_{|S_{\lb}},{\gm}_{|S_\lb})
$$
given by
$$
a\longmapsto E_p(a S):=\sum_{i=0}^{p-1}\frac{a^i}{i!}S^i.
$$
\end{prop}
\begin{proof}
See \cite[4.15, 4.22]{io2}. We remark that there is a similar
statement for $v(\mu)<v(\lb)$ which we omit since,  as we will
see, it is not of  interest in this paper. \end{proof} We now
remark about   definition \ref{def:clE(mu,lb,F,i)}. Applying the
functor $\Hom_{gr}(\cdot,{\gm}_{|S_{\lb}})$ to   the exact
sequence over $S_\lb$
$$ 0\too G_{\mu,1}\on{i}{\too}
\gmu\on{\psi_{\mu,1}}{\too}\gmup\too 0,
$$
we have that
\begin{equation}\label{eq:ker i_*}
\ker \bigg(i^*:\Hom_{gr}(\gmu_{|S_\lb},{\gm}_{|S_\lb})\too
\Hom_{gr}({G_{\mu,1}}_{|S_\lb},{\gm}_{|S_\lb})\bigg)={\psi_{\mu,1}}_*\Hom_{gr}(\gmup_{|S_\lb},{\gm}_{|S_\lb}).
\end{equation}
So let $F(S)\in
\Hom_{gr}({\gmx{1}}_{|S_{\lb}},{\gm}_{|S_{{\lb}}})$. By
\ref{lem:suriettività mappa tra hom} we can choose a
representative of $F(S)$ in $
\Hom_{gr}({\gmu}_{|S_{\lb}},{\gm}_{|S_{{\lb}}})$ which we denote
again $F(S)$ for simplicity. We remark that $F(S)^p\in
\Hom_{gr}({\gmu}_{|S_{\lb^p}},{\gm}_{|S_{\lb^p}})$. Therefore, by
\eqref{eq:ker i_*}, we have that
$$F(S)^p(1+\mu S)^{-j}=1\in
\Hom_{gr}({\gmx{1}}_{|S_{\lb^p}},{\gm}_{|S_{{\lb^p}}})$$ is
equivalent to saying that there
exists $G\in\Hom_{gr}({\gmup}_{|S_{\lb^p}},{\gm}_{|S_{\lb^p}})$
with the property that $$F(S)^p(1+\mu S)^{-j}=G(\frac{(1+\mu
S)^p-1}{\mu^p})\in
\Hom_{gr}({\gmu}_{|S_{\lb^p}},{\gm}_{|S_{{\lb^p}}}).$$ This
implies that $\clE^{(\mu,\lb;F,j)}$  can be seen as the kernel of
the isogeny
\begin{align*}
\psi^j_{\mu,\lb,F,G}:\clE^{(\mu,\lb;F)}&\too \clE^{(\mu^p,\lb^p;G)}\\
          S_1&\longmapsto \frac{(1+\mu S_1)^p-1}{\mu^p}\\
          S_2&\longmapsto \frac{(\tilde{F}(S_1)+\lb
S_2)^p(1+\mu S_1)^{-j}-\tilde{G}(\frac{(1+\mu
S_1)^p-1}{\mu^p})}{\lb^p}
\end{align*}
where $\tilde{F},\tilde{G}\in R[S]$  are liftings of $F$ and $G$.

\begin{ex}\label{ex:equazione per eta_pi} Let us define
$$
\eta=\sum_{k=1}^{p-1}\frac{(-1)^{k-1}}{k}\lb_{(2)}^k.
$$
We remark that $v(\eta)=v(\lb_{(2)})$. We consider
$$
F(S)=E_p(\eta S)=\sum_{k=1}^{p-1}\frac{(\eta S)^k}{k!}.
$$
It was showed in \cite[\S 5]{SS4} that, using our notation,
$$
\Z/p^2\Z\simeq \clE^{(\lb_{(1)},\lb_{(1)};E_p(\eta T),1)}
$$
as group schemes. A similar description of $\Z/p^2\Z$ was
independently found by Green and Matignon (\cite{GM1}). 
\end{ex}
\begin{ex}\label{ex:Gmun} It is easy to see that  $\glx{2}$ 
is isomorphic, as a group scheme, to $ \clE^{(\lb^p,\lb;1,1)}.$
\end{ex}

 Let us define, for any $\mu,\lb\in R$ with $v(\mu), v(\lb)\le v(\lb_{(1)})$, the group
\begin{align*}
rad_{p,\lb}(<1+\mu S>):=\bigg\{(F(S),j)\in
\Hom_{gr}&({\gmx{1}}_{|S_\lb},{\gm}_{|S_\lb})\times \Z/p\Z \text{
such that }\\
 &F(S)^p(1+\mu S)^{-j}=1\in
\Hom_{gr}({\gmx{1}}_{|S_{\lb^p}},{\gm}_{|S_{{\lb^p}}})\bigg\}.
\end{align*}
We can define the map $$ \beta: rad_{p,\lb}(<1+\mu S>) \on{}{\too}
\Ext^{1}(G_{\mu,1},\glx{1})$$ by
$$ (F(
S),j)\longmapsto \clE^{(\mu,\lb;F( S),j)}
$$
It has been proved in \cite[4.38, 4.39]{io2} that it is an injective
morphism of groups. The image of $\beta$ is the set
$\{\clE^{(\mu,\lb;F( S),j)}\}$, which is therefore a group
isomorphic to $rad_{p,\lb}(<1+\mu S>)$.

\begin{thm}
Let us suppose $p>2$. Let $\mu, \lb \in R\setminus \{0\}$ be with
$v(\lb_{(1)})\ge v(\mu)\ge v(\lb)$. Then, \mbox{$
rad_{p,\lb}(<1+\mu S>) $} is isomorphic to the group
\begin{align*}
\Phi_{\mu,\lb}:=\bigg\{(a,j)\in (R/\lb R)^{\fr}\times\Z/p\Z&
\text{ such that } pa-j\mu=\frac{p}{\mu^{p-1}}a^p\in R/\lb^p
R\bigg\},
\end{align*}
through the map
$$
(a,j)\longmapsto (E_p(a S),j) 
$$
\end{thm}
\begin{proof}
\cite[4.47]{io2}
\end{proof}
\begin{rem}\label{rem:0 in Phi}
It is clear that if $(0,j)\in \Phi_{\mu,\lb}$, with $j\neq 0$,
then $\mu\equiv 0 \mod \lb^p$.
\end{rem}
\begin{rem}\label{rem:G=E(ap S)}From the proof of this theorem, it also follows that
$\clE^{(\mu,\lb,F,j)}$, with $F(S)=E_p(a S)\in
\Hom_{gr}({\gmx{1}}_{|S_{\lb}},{\gm}_{|S_{{\lb}}})$, is the kernel
of the isogeny $\psi^j_{\mu,\lb,F,G}$, with $G(S)=E_p(a^p S)\in
\Hom_{gr}({\gmx{1}}_{|S_{\lb^p}},{\gm}_{|S_{{\lb^p}}})$.
\end{rem}
We recall that we are interested in group schemes which are
generically isomorphic to $\Z/p^2\Z$ and not in extensions.   We
call $$\Phi_{\mu,\lb}^{1}:=\{a\in R/\lb R| (a,1)\in
\Phi_{\mu,\lb}\};$$ it is only a set. The following result says
that every extension represented by an element of
$\Phi_{\mu,\lb}^{1}$ is a model of $(\Z/p^2\Z)_K$ as a group
scheme.
\begin{thm}\label{cor:clE se lb divide mu} Let us suppose that $R$
contains a primitive $p^2$-th root of unity.
Let $G$ be a model of 
$(\Z/p^2\Z)_K$. Then there exist unique $v(\lb_{(1)})\ge m\ge n\ge
0$ and $a\in \Phi_{\pi^m,\pi^n}^1$ such that $
G\simeq\clE^{(\pi^m,\pi^n;E_p(a S),1)}$  and moreover
$\Phi_{\pi^m,\pi^n}^1$ is one of the following.
\begin{itemize}
\item[a)] If $m<pn$, $\Phi_{\pi^m,\pi^n}^1 $ is nonempty if and
only if $pm-n\ge v(p)$, and it  is the set
\begin{align*}
\bigg\{\eta \frac{\pi^m}{\lb_{(1)}}+ \alpha\in R/\pi^n R \text{
s.t., for any
  lifting }&\text{$\tilde{\alpha}\in R$ of $\alpha\in R/\pi^n R$},\\
pv(\tilde{\alpha})\ge\max&\{pn+(p-1)m-v(p),n\}\bigg\}.
\end{align*}
For the definition of $\eta$ see \ref{ex:equazione per eta_pi}. 
\item[b)] If $m\ge pn$ then  $ \Phi_{\pi^m,\pi^n}^1$ is the set
$$\bigg\{\alpha\in R/\lb R\text{ s. t., for any
  lifting } \tilde{\alpha}\in R, \
pv(\tilde{\alpha})\ge\max\{pn+(p-1)m-v(p),n\}\bigg\}
$$
\end{itemize}
Conversely if $a\in \Phi^1_{\pi^m,\pi^n}$ then
$\clE^{(\pi^m,\pi^n;E_p(a S),1)}$ is a model of $(\Z/p^2\Z)_K$.
\end{thm}
\begin{proof}
See \cite[4.51 and 4.58]{io2}. The notation $\Phi_{\pi^m,\pi^n}^1$
is not used there.
\end{proof}
\begin{rem}\label{rem:valutazione di elementi di Phi1}Let us suppose $m<pn$. Let $b\in \Phi^1_{\pi^m,\pi^n}$. By \ref{rem:0 in Phi}, then $b\neq 0$. Let $\tilde{b}\in
R$ be any of its lifting. Then
$v(\tilde{b})=v(\eta\frac{\pi^m}{\lb_{(1)}})=m-\frac{v(p)}{p}$.
Indeed, by the theorem, we have
$\tilde{b}=\eta\frac{\pi^m}{\lb_{(1)}}+\alpha$ for some $\alpha
\in R/\pi^n R$ with
$v(\tilde{\alpha})>v(\eta\frac{\pi^m}{\lb_{(1)}})=m-\frac{v(p)}{p}$,
where $\tilde{\alpha}\in R$ is any lifting of $\alpha$.
\end{rem}

The following proposition tells us, in particular, when two
extensions as above are isomorphic as group schemes.
\begin{prop}\label{lem:abbasso valutazione con blow-up} 
Let $a_i \in \Phi_{\mu,\lb}^1$, for $i=1,2$. We set
$F_i(S)=E_p(a_i S)$. Let
$$
f:\clE^{(\mu_1,\lb_{1};F_1,1)}\too \clE^{(\mu_2,\lb_{2};F_2,1)},
$$
 be a model map between models of $(\Z/p^2\Z)_K$. Then it is a morphism of extensions and it is given by
\begin{equation}\label{eq:def f}
\begin{aligned}
S_1&\mTo \frac{(1+\mu_1 S_1)^{\frac{rj_1}{j_2}}-1}{\mu_2} \\
S_2&\mTo \frac{(F_1(S_1)+{\lb_1}S_2)^{r}(1+\mu_1
S_1)^{s}-F_2(\frac{(1+\mu_1
S_1)^{\frac{rj_1}{j_2}}-1}{\mu_2})}{{\lb_2}}
\end{aligned}
\end{equation}
for some $r\in (\Z/p\Z)^*$ and $s\in \Z/p\Z$. Moreover, it exists
if and only if  $v(\mu_1)\ge v(\mu_2)$, $v(\lb_1)\ge v({\lb_2})$
and
\begin{equation}\label{eq:a_1=k_1/k_2...
} a_1\equiv \frac{\mu_1}{\mu_2}a_2\mod{\lb_2}.
\end{equation}
 Moreover, any  such  $f$ is an isomorphism if  and only if $v(\mu_1)=v(\mu_2)$ and
$v(\lb_1)=v({\lb_2})$.
\end{prop}
\begin{proof}
See \cite[4.54, 4.55]{io2}. There the statements  are slightly
more general.
\end{proof}



\subsection{Torsors under $G_{\lb,n}$ and  $\clE^{(\mu,\lb;F(S),j)}$}\label{sec:clE
torsors} We begin describing explicitly $G_{\lb,n}$-torsors.

\begin{prop}\label{prop:h1 glb uguale in zar e flat} Let $X$ be a faithfully flat
$R$-scheme and  let $f:X_{fl}\too X_{Zar}$ be the natural
continuous morphism of sites.   Then $R^1f_*(\glb)=0$. In
particular $H^1(X,\glb)=H^1(X_{Zar},\glb)$.
\end{prop}

\begin{proof}
It is sufficient to prove that $H^1(\Sp(A),\glb)=0$ for any local
ring $A$ flat over $R$. If  $\lb =0$, then $\glb\simeq \Ga$ and the
statement is classical  (\cite[III 2.14,3.7]{Mi}). Let us now
suppose $\lb\neq 0$.
 We recall that we have the following exact
sequence (see \eqref{eq:succ esatta lunga di g lambda})
\begin{multline*}
0\too H^0(\Sp(A),\glb)\too H^0(\Sp(A),{\gm})\too
H^0({\Sp(A/\lb A)},{\gm})\too\\
 \ \ \ \ \ \ \ \ \ \too H^1(\Sp(A),{\glb})\too
H^1(\Sp(A),{\gm})\too \hg{\Sp(A/\lb A)}{{\gm}}.\\
\end{multline*}
But, since $A$ is local, then $ H^0(\Sp(A),{\gm})\too
H^0({\Sp(A/\lb A)},{\gm})$ is surjective. Moreover by Hilbert 90
(see \cite[III.4.9]{Mi}) we have
$\hg{\Sp(A)}{{\gm}}=\Pic(\Sp(A))$. But, since $A$ is local then
$\Pic(\Sp(A))=0$. So $H^1(\Sp(A),\glb)=0$, as necessary. Since
$R^1f_*(\glb)=0$, it follows by the Leray spectral sequence that
$$H^1(X,\glb)=H^1(X_{Zar},\glb).$$
\end{proof}

If $Y\too X$ is a $G$-torsor, we will denote by $[Y]$ the class of
the $G$-torsor  in $\hxg{G}$. Now let $Y\too X$ be a
$\gln$-torsor. Let us consider the exact sequence \eqref{eq:succ
esatta glbn}. By \ref{prop:h1 glb uguale in zar e flat} we can
take a covering $\{U_i=\Sp(A_i)\}$ of $X$ by affine subschemes
such that the class $(i_*[Y])_{|U_i}$ is trivial, where
$i:\gln\too \glb$. This means that
\begin{equation}\label{eq:glbn torsors locally}Y_{|U_i}=\Sp\bigg(A_i[T_i]/(\frac{(1+\lb
T_i)^{p^n}-1}{\lb^{p^n}}-f_i)\bigg) \end{equation} for some
$f_i\in A_i$ such that $1+\lb^{p^n}f_i\in A_i^*$. Moreover, $1+\lb
T_i=f_{ij}(1+\lb T_j)$ with $\{f_{ij}=1+\lb g_{ij}\}=i_*[Y]\in
\hxg{\glb}$.

\begin{rem}
Andreatta and Gasbarri (see \cite{AG2}) have given a description
of $G_{\lb,n}$-torsors from which they deduced that a
$G_{\lb,n}$-torsor is locally
 \eqref{eq:glbn torsors
locally}. From this fact they deduced that
$H^1(X,\glb)=H^1(X_{Zar},\glb).$
\end{rem}

We now give an explicit description of
$\clE^{(\mu,\lb;F(S),j)}$-torsors. 
Let $\clE^{(\mu,\lb;F(S),j)}$ be a group scheme as in the previous
subsection. We have seen that  
there is the following exact sequence
$$
0\too \clE^{(\mu,\lb;F(S),j)}\on{\iota}{\too}
\clE^{(\mu,\lb;F(S))}\on{\psi^j_{\mu,\lb,F,G}}{\too}
\clE^{(\mu^p,\lb^p;G(S))}\too 0
$$
for  $G(S)\in
\Hom_{gr}({G^{(\mu^p)}}_{|S_{\lb^p}},{\gm}_{|S_{\lb^p}})$ such
that $F(S)^p(1+\mu S)^{-j}=G(\frac{(1+\mu S)^p-1}{\mu^p})\in
\Hom_{gr}({\gmu}_{|S_{\lb^p}},{\gm}_{|S_{{\lb^p}}})$. The
associated long exact sequence is
\begin{equation}\label{eq:succ esatta lunga torsori}
\begin{aligned}
\dots \too
H^0(X,\clE^{(\mu,\lb;F(S))})\on{{(\psi^j_{\mu,\lb,F,G})}_*}\too
H^0(X,\clE^{(\mu^p,\lb^p;G(S))})\on{\delta}{\too}
H^1(X,\clE^{(\mu,\lb;F(S),j)})&\on{\iota_*}{\too}\\
\too H^1(X,\clE^{(\mu,\lb;F(S))})\too
H^1(X,&\clE^{(\mu^p,\lb^p;G(S))})\too \dots
\end{aligned}
\end{equation}
\begin{cor}\label{cor:R1f=0 per gruppi dim 2}
Let $X$ be a faithfully flat $R$-scheme and  let $f:X_{fl}\too
X_{Zar}$ be the natural continuous morphism of sites. For any
$R$-group scheme $\clE^{(\mu,\lb;F)}$, we have
$R^1f_*(\clE^{(\mu,\lb;F)})=0$. In particular,
$H^1(X,\clE^{(\mu,\lb;F)})=H^1(X_{Zar},\clE^{(\mu,\lb;F)})$.
\end{cor}
\begin{proof}
Let us consider the exact sequence \eqref{eq:estensioni lisce}, in
the fppf topology,
$$
0\too \glb\too \clE^{(\mu,\lb;F)}\too \gmu\too 0.
$$
If we apply the functor $f_*$, we obtain
$$
\dots \too R^1f_*\glb\too R^1f_*(\clE^{(\mu,\lb;F)})\too
R^1f_*\gmu\too \dots
$$
By \ref{prop:h1 glb uguale in zar e flat} it follows that
$R^1f_*(\glb)=R^1f_*(\gmu)=0$. Hence,
$R^1f_*(\clE^{(\mu,\lb;F)})=0$. As in \ref{prop:h1 glb uguale in
zar e flat}, we conclude that
$H^1(X,\clE^{(\mu,\lb;F)})=H^1(X_{Zar},\clE^{(\mu,\lb;F)})$.
\end{proof}

 Let
$\tilde{F}(S),\tilde{G}(S)$  be liftings of $F(S)$ and $G(S)$ in
$R[S]$. We remark that
\begin{align*}
H^0(X,\clE^{(\mu,\lb;F(S))})=\{(f_1,f_2)&\in H^0(X,\ox)\times
H^0(X,\ox)| \\
1+\mu f_1\in &H^0(X,\ox)^* \text{ and } \tilde{F}(f_1)+\lb f_2\in
H^0(X,\ox)^*\},
\end{align*}
 and
\begin{align*}
H^0(X,\clE^{(\mu,\lb;G(S))})=\{(f_1,f_2)&\in H^0(X,\ox)\times
H^0(X,\ox)|\\ 1+\mu^p f_1&\in H^0(X,\ox)^* \text{ and }
\tilde{G}(f_1)+\lb^pf_2\in H^0(X,\ox)^*\}.
\end{align*}
 The map
$(\psi^j_{\mu,\lb,F,G})_*:H^0(X,\clE^{(\mu,\lb;F(S))})\too
H^0(X,\clE^{(\mu,\lb;G(S))})$ is given by
$$
(f_1,f_2)\longmapsto (\frac{(1+\mu
f_1)^p-1}{\mu^p},\frac{(\tilde{F}(f_1)+\lb f_2)^p(1+\mu
f_1)^{-j}-\tilde{G}(\frac{(1+\mu f_1)^p-1}{\mu^p})}{\lb^p}).
$$
Now let us suppose that $X=\Sp(A)$. We describe explicitly the map
$\delta$. We have
$$
H^0(X,\clE^{(\mu,\lb;G(S))})=\{(f_1,f_2)\in A\times A| 1+\mu^p
f_1\in A^* \text{ and } \tilde{G}(f_1)+\lb^pf_2\in A^*\}.
$$and  $\delta((f_1,f_2))$ is, as a scheme,
\begin{equation}\label{eq:clE-torsori localmente}
Y=\Sp A[T_1,T_2]/\bigg(\frac{(1+\mu
T_1)^p-1}{\mu^p}-f_1,\frac{(\tilde{F}(T_1)+\lb T_2)^p(1+\mu
T_1)^{-j}-\tilde{G}(f_1)}{\lb^p}-f_2\bigg),
\end{equation}
and  the $\clE^{(\mu,\lb;F(S),j)}$-action over $Y$ is given by
\begin{align*}
    T_1\longmapsto &S_1+ T_1+\mu S_1T_1\\
    T_2\longmapsto &S_2 \tilde{F}(T_1)+\tilde{F}(S_1) T_2 +\lb S_2T_2+\\
                   & \qquad   \quad    \frac{\tilde{F}(S_1)\tilde{F}(T_1)-\tilde{F}(S_1+T_1+\mu S_1T_1)}{\lb}
    \end{align*}

Now let $X$ be  a faithfully $R$-scheme. If $Y\too X$ is a
$\clE^{(\mu,\lb;F(S),j)}$-torsor then, by \ref{cor:R1f=0 per
gruppi dim 2} there exists a Zariski covering $\{U_i=\Sp(A_i)\}$
such that $(\iota_*)_{U_i}([Y])=0$, for any $i$. This means that
$$Y_{|U_i}=\Sp A_i[T_1,T_2]/\bigg(\frac{(1+\mu
T_1)^p-1}{\mu^p}-f_{1,i},\frac{(\tilde{F}(T_1)+\lb T_2)^p(1+\mu
T_1)^{-j}-\tilde{G}(f_1)}{\lb^p}-{f_{2,i}}\bigg),$$ for some
$f_{1,i},f_{2,i}$ as above. By a standard argument we can see that
the cocycle $\iota_*([Y])\in H^1(X_{Zar},\clE^{(\mu,\lb;F)})$
permits to patch together the torsors $Y_{|U_i}$ to obtain $Y$. In
\S\ref{sec:deg Z/p^2Z-torsors}, we will consider  the case $X$
affine and we are interested only in
$\clE^{(\mu,\lb;F,j)}$-torsors of the form \eqref{eq:clE-torsori
localmente}.

\section{A filtration of $\hxg{\mun}$}\label{sec:filtration}
\stepcounter{subsection}\setcounter{subsection}{0}

We firstly recall the following result.
\begin{prop}\label{prop:iniezione tra h1}
Let $X$ be a normal integral scheme.  Let $f:Y\too X$ be a
morphism with a rational section and let $g:G\too G'$ be a map of
finite and flat commutative  group schemes over $X$, which is an
isomorphism over $\Sp(K(X))$. Then
$$
f^*g_*:\hxg{G}\too H^1(Y,G'_{Y})
$$
is injective.
\end{prop}
\begin{proof}
See \cite[4.29]{io2}.
\end{proof}

\begin{rem}The previous result will be applied to the case
$f=\id_X$ or to the case $f:U\too X$ an open immersion and
$g=\id_G$.
\end{rem}

We refer to \S \ref{sec:modelli} for the notation. For any
$\lb,\mu\in R$ with $v(\lb)\ge v(\mu)$ we have a morphism
\mbox{$\glb \too \g^{(\mu)}$} defined by $T\mTo \frac{\lb}{\mu}
T$. If $v(p)\ge p^{n-1}(p-1)v(\lb)$, it is compatible with
$\psi_{\lb,n}$ and $\psi_{\mu,n}$. So it induces a morphism $\gln
\too G_{\mu,n}$ such that
$$
\xymatrix{ \gln \ar[rd]^{\alpha^{\lb,n}} \ar[rr]  &
\ar[r]&\ar[ld]_{\alpha^{\mu,n}
} G_{\mu,n}\\
                                     & \mun}
$$
commutes. Let $X$ be a normal integral faithfully flat $R$-scheme.
We obtain the following commutative diagram
$$
\xymatrix{ \hxg{\gln}\ar[rd]^{(\alpha^{\lb,n})_*}\ar[rr] & &\ar[ld]_{(\alpha^{\mu,n})_*}\hxg{G_{\mu,n}}\\
                                                          & \hxg{\mun}}
$$
and, applying  \ref{prop:iniezione tra h1}, we have that
$(\alpha^{\lb,n})_*$ and $(\alpha^{\mu,n})_*$ are injective.
Hence, $\hxg{\gln}{\too} \hxg{G_{\mu,n}}$ is injective. So we have
 proved the following proposition.
\begin{prop}\label{cor:filtration}
Let $X$ be a normal integral faithfully flat $R$-scheme and  let
$i_0=\max\{i|v(p)\ge p^{n-1}{(p-1)}v(\pi^i)\}$. Then, for any $n$,
we have the following filtration
$$
0\In\hxg{G_{\pi^{i_0},n}}\In\hxg{G_{\pi^{i_0-1},n}}\In\dots
\In\hxg{G_{\pi,n}}\In \hxg{\mun}.
$$

\end{prop}

When comparable the above filtration coincides with that of
\cite[5.2]{Sa}. Moreover it coincides with that of
\cite[p.262]{Rob} if $n=1$ and $X=\Sp(A)$ with $A$ the integer
ring of a local number field.
%
%
If $R$ contains  a primitive ${p^n}$-{th} root of unity, then
 $i_0=v(\lb_{(n)})$.

\section{Weak extension of torsors under commutative group
schemes}\label{sec:estensione debole} The aim of this section is to
prove a result of weak extension for torsors under commutative group
schemes over normal schemes with some hypothesis.

\subsection{Preliminary results}

 We here state some results which
will be useful in what follows.

\begin{prop}\label{lem:morfismo di
gruppi deriva da quello di azione} Let $i=1,2$.  Let $Z_i$ be a
faithfully flat $S$-scheme and  let $G_i$ be an affine flat
$S$-group scheme, together with an admissible action, over a
faithfully flat $Z_i$-schemes $Y_i$. Moreover, we suppose that
$Y_2\too Z_2$ is a $G_2$-torsor and that there exists a morphism
$$
\phi_K:(G_1)_K\too (G_2)_K.
$$
Let us suppose we have a commutative diagram
$$
\xymatrix{&Y_1\ar[d]\ar[rr]^f&\ar[r]&\ar[d] Y_2\\
            &Z_1\ar[rr]&\ar[r]    &Z_2&}
$$
of $S$-schemes such that $f_K$ is an  isomorphism compatible with
the actions. Then there exists a unique morphism
$$
\phi:G_1\too G_2
$$
which extends $\phi_K$ and such that  $f$ is compatible with the
actions.
\end{prop}
\begin{proof}
For $i=1,2$ we call $\sigma_i:G_i\times_R Y_i\too Y_i$ the
actions. 
Since $Y_2\too Z_2$ is a $G_2$-torsor, then $\sigma_2\times \id$
is
an isomorphism. 
So by
$$
G_1\times_R Y_1\on{\sigma_1\times \id}{\too} Y_1 \times_{Z_1}
Y_1\on{f\times f}{\too} Y_2 \times_{Z_2} Y_2 \on{(\sigma_2\times
\id)^{-1}}\too G_2 \times_R Y_2
$$
we obtain a morphism
$$
G_1\times_R Y_1 \too G_2\times_R Y_2.
$$
 If we compose it with
the projection $p_1:G_2\times_R Y_2\too G_2$, we obtain a morphism
$$
G_1\times_R Y_1\too G_2.
$$
Moreover, we consider the projection
$$
p_2:G_1\times_R Y_1\too Y_1.
$$
Therefore, we have  a map
\begin{equation}\label{eq:morfismi da G_Y}
\phi_{Y_1}:G_1\times_R Y_1\too G_2\times_R Y_1.
\end{equation}
We now prove that it  is compatible with $\phi_K$, i.e.
$\phi_{Y_1}$ and $\phi_K$ induce the same morphism
\mbox{$G_1\times_R (Y_1)_K\too G_2\times_R (Y_1)_K$.} We observe
that $\phi_{Y_1}$ and $\phi_K$ induce two morphisms,
$(\phi_{Y_1})_K$ and $(\phi_K)_{Y_1}$ respectively,
which are compatible with $f_K$. 
For any $\psi:G_1\times_R (Y_1)_K\too G_2\times_R (Y_1)_K$, to be
compatible with $f_K$ means that the following diagram
$$
\xymatrix{G_1\times_R (Y_1)_K\ar[d]_{(\id\times f_K)\circ\psi}\ar[r]^(.6){\sigma_1}&\ar[d]^{f_K} (Y_1)_K \\
             G_2\times_R(Y_2)_K\ar[r]^(.6){\sigma_2}& (Y_2)_K                                  }
$$
commutes, i.e.  $\sigma_2\circ(\id\times f_K)\circ\psi=f_K\circ
\sigma_1$. So we have
$$
\sigma_2\circ(\id\times
f_K)\circ(\phi_{Y_1})_K=\sigma_2\circ(\id\times
f_K)\circ(\phi_K)_{Y_1}=f_K\circ \sigma_1.
$$
Moreover, since $p_2\circ (\id\times f_K)=f_K\circ p_2$ and
$p_2\circ (\phi_{Y_1})_K=p_2\circ {(\phi_K)}_{Y_1}= p_2$ it
follows that
$$
p_2\circ(\id\times f_K)\circ (\phi_{Y_1})_K=p_2\circ (\id\times
f_K)\circ (\phi_K)_{Y_1}=f_K\circ p_2.
$$
Since $(Y_2)_K\too (Z_2)_K$ is a $(G_2)_K$-torsor, then
$$
(\id\times f_K)\circ(\phi_{Y_1})_K=(\id\times
f_K)\circ(\phi_K)_{Y_1}.
$$
 For $i=1,2$, let $p_i$ be the projections from $G_2\times_R(Y_1)_K$
and let  $p'_i$ be the projections from $G_2\times_R(Y_2)_K$. Then
$$
p_1\circ(\phi_{Y_1})_K=p'_1\circ(\id\times
f_K)\circ(\phi_{Y_1})_K=p'_1\circ(\id\times
f_K)\circ(\phi_K)_{Y_1}=p_1 \circ(\phi_K)_{Y_1}.
$$
and
$$
f_K\circ p_2\circ (\phi_{Y_1})_K=p'_2\circ(\id\times
f_K)\circ(\phi_{Y_1})_K=p'_2\circ(\id\times
f_K)\circ(\phi_K)_{Y_1}=f_K\circ p_2\circ (\phi_K)_{Y_1}.
$$
Since $f_K$ is an isomorphism, then $p_2\circ
(\phi_{Y_1})_K=p_2\circ (\phi_{K})_{Y_1}$.  Hence,
$(\phi_{Y_1})_K=(\phi_{K})_{Y_1}$, i.e. $\phi_K$ is compatible
with $\phi_{Y_1}$.
By the next descent lemma  we have a unique morphism of schemes
$\phi:G_1\too G_2$ which extends $\phi_K$ and $\phi_{Y_1}$. Since
$G_1$ is flat over $R$, $\phi_K$ is a morphism of group schemes
and $G_2$ a separated scheme over $R$,
then $\phi$ is a morphism of group schemes. By construction it is
clear that, through $\phi$, the morphism $f$ preserves the
actions.
\end{proof}

We now prove  the  descent lemma  used in the previous proof.
\begin{lem}\label{prop:morfismo di gruppi estendibile}
Let $S'\to S$ be a faithfully flat morphism.  Let $X_1,X_2$ be
affine $S$-schemes with $X_2$ flat over $S$. Given two morphisms
$\varphi_K:(X_1)_K\to (X_2)_K$ and $\varphi_{S'}:X_1\times_S S'\to
X_2\times_S S'$ that coincide on $S'_K$, there is a unique
morphism $\varphi:X_1\to X_2$ that extends them.
\end{lem}
\begin{proof}
Up to restricting ourselves to  an affine  subscheme of $S'$, we can
suppose \mbox{$S'=\Sp(A)$.} For $i=1,2$, let us consider
$X_i=\Sp(B_i)$. In terms of function rings we have two morphisms
\begin{equation}\label{eq:phi A}
\phi_{A}^{\sharp}:=\phi^{\sharp}_{S'}:B_2\pt_R A \too B_1\pt_R A
\end{equation}
and
\begin{equation}\label{eq:phi K}
\phi^{\sharp}_K: B_2\pt_R K \too B_1\pt_R K.
\end{equation} Moreover, by compatibility, it follows that the above morphisms
induce the same map
$$
 \phi^\sharp_{A_K}:=\phi^\sharp_A\pt \id_K=(\phi^\sharp_K)\pt{\id_{A_K}}:B_2\pt_R A\pt_R K
 \too
B_1\pt_R A \pt_R K.
$$
First of all we prove the uniqueness of $\phi$. Since  $X_2$ is
flat over $S$, then the inclusion $(X_2)_K\too X_2$ induces an
injection $B_2\hookrightarrow B_2\pt_R K$. Therefore, if (any)
$\phi$ exists then it is given, on the level of function rings, by
the restriction of $\phi^\sharp_{K}$ to $B_2$. Therefore it is
unique.

We now prove the existence of $\phi$.  
We have the following commutative diagram with the obvious maps
 $$\xymatrix{&             B_i\pt_R A\ar@^{->}[rd]^{}\ar@_{>}[rd]& \\
                       B_i\ar@^{->}[ru]\ar@_{>}[ru]\ar@^{->}[rd]\ar@_{>}[rd]  && B_i\pt_R A\pt_R K\\
                                   & B_i\pt_R K \ar@^{->}[ru] \ar@_{>}[ru]    }$$
Since $S'\too S$ is flat, in particular, the induced map $R\too A$
is injective. Moreover, $B_2$ is a flat $R$-algebra, then for
$i=2$ all the maps of the above diagram are injective.
 We remark that \eqref{eq:phi A}
and \eqref{eq:phi K} imply $\phi^\sharp_{A_K}(B_2)\In (B_1\pt_R
K)\cap (B_1\pt_R A)$. We claim that $\phi^\sharp_{A_K}(B_2)\In B_1$.
Let us suppose that there exists $b\in B_2$ such that
$\phi^\sharp_{A_K}(b)\not\in B_1$. Then there exists $n\ge 1$ such
that
 $\pi^n\phi^\sharp_{A_K}(b)\in B_1$ and $\pi^{n-1}\phi^\sharp_{A_K}(b)\not\in B_1$. 
Hence, 
$$
 \pi^n \phi^\sharp_{A_K}(b)\in B_1\cap \pi ( B_1\pt_R A ).
$$
We remark that since $S'\too S$ is surjective, then
$S'_k=\Sp(A/\pi A)$ is nonempty. Now, since any scheme over a
field is flat,
$$
B_1/\pi B_1\too ( B_1\pt_R A) /\pi (B_1\pt_R A )\simeq B_1/\pi
B_1\pt_k A/\pi A
$$
is injective. Therefore,
$$
 B_1\cap \pi( B_1\pt_R A )=\pi B_1,
$$
which implies $\pi^{n-1}\phi^\sharp_{A_K}(b)\in B_1$. This is a
contradiction. So $\phi^\sharp_{A_K}$ induces a morphism
$$
\phi^\sharp:B_2\too B_1.
$$
 We have so proved that $\phi_K:(Y_1)_K\too (Y_2)_K$
is extendible to a morphism
\mbox{$\phi:Y_1\too Y_2$.} 
%
\end{proof}
The following easy consequence of the previous lemma  will not be
used in the rest of the paper.
\begin{prop} Let  $G$ be an affine flat and commutative $S$-group scheme. Then
$$H^1(S,G)\too H^1({K},G_K)$$ is injective.
\end{prop}
\begin{proof} Let $f:Y\too S$ be a $G$-torsor. This means that
there exists a faithfully flat $S$-scheme $T$ such that
$Y_T:=Y\times_X T\longmapsto T$ is trivial (for instance we can
chose $T=Y$).  Then it has a section $\phi_T:T\too Y_T$. Moreover
let us suppose that $Y\too S$ is trivial as $G_K$-torsor on $X_K$.
Then there is a section $\phi_K:\Sp(K)\too Y_K$ of $Y_K\too
\Sp(K)$. Since $G$ is affine then $f:Y\too S$ is an affine
morphism. So $Y$
is affine. From the previous lemma the thesis follows. 
\end{proof}
\begin{lem}\label{lem:B/A flat}
Let $X,Y$ be integral faithfully flat schemes over $S$. Moreover,
let us suppose that $X$ is normal. If $f:Y\too X$ is an integral
dominant $S$-morphism then, $f_k$ is schematically dominant, i.e.
$f_k^{\sharp}:\oo{X_k}\too f_*\oo{Y_k}$ is injective. In
particular, if $Y_k$ is integral (resp. reduced) then $X_k$ is
also integral (resp. reduced).
\end{lem}
\begin{proof}Since any integral morphism is affine by the definition,  it is enough to prove the lemma in the affine
case. So we  can suppose  $X=\Sp(A)$, $Y=\Sp(B)$ with an integral
injection
$A\hookrightarrow B$.  We will prove that
$A_k\hookrightarrow B_k$. 
This is equivalent to proving $\pi B\cap A=\pi A$. One inclusion is
obvious. Now let $a\in \pi B\cap A$, then $a=\pi b$ with $ b\in B$.
We remark that $b=\frac{a}{\pi}\in A_K\cap B$ is integral over $A$.
But $A$ is integrally closed by hypothesis. Therefore $b\in A$.

It follows immediately the last statement.
\end{proof}
\begin{lem}\label{lem:normality's criterion}
 Let $X$ be a flat $S$-scheme. If $X_K$ is normal and $X_k$
reduced, then $X$ is normal.
\end{lem}
\begin{proof}
For a  proof  see \cite[4.1.18]{liu}.
 \end{proof}

\subsection{Weak extension}
 We now consider   an integral  normal  faithfully
flat  $S$-scheme $X$. 
Let $Y_K\too X_K$ be a connected $G_K$-torsor, for some finite
group-scheme $G_K$ over $K$, and $Y$ the normalization of $X$ in
$Y_K$. We remark that $Y_K$ is also normal (\cite[I 9.10]{SGA1}).
In particular $Y_K$ is integral, hence $Y$ is integral. We denote
by $g$ the morphism $Y\too X$.

\begin{lem}
 $Y$ is normal.
\end{lem}
\begin{proof} By definition of normalization we have that $g$ is an affine morphism.
Since  being normal is a local property we can suppose that
$X=\Sp(A)$ and $Y=\Sp(B)$. We have to prove that  $B$ is
integrally closed.
 Since $Y_K$ is normal, then $B_K$ is integrally closed.
So, if $b\in Frac(B)$ is integral over $B$, then $b\in B_K$.
However, $B$ is integral over $A$, then $b$ is integral over $A$.
But $B$ is the integral closure of $A$ in $B_K$, therefore $b\in
B$.
\end{proof}
From the proof it follows that, since $X$ is normal,  in fact $Y$
is the normalization of $X$ in $K(Y)$, the function field of $Y$.
Therefore  $g:Y\too X$ is finite (see \cite[4.1.25]{liu}). We
remark that $Y$ has also the following property.
\begin{lem}\label{lem:universal property of normal closure}
\noindent Let $Z$ be a faithfully flat $S$-scheme, let $f:Z\too X$
be a dominant integral $S$- morphism  and let $h_K:Y_K\too Z_K$ be
a dominant $X_K$-morphism. Then there exists an integral
$X$-morphism $h:Y\too Z$ which extends $h_K$. Moreover if $Z$ is
normal and $h_K$ is an isomorphism then $h$ is an isomorphism.
\end{lem}
\begin{proof}
We have that  $g:Y\too X$ and $f:Z\too X$ are affine.  So, first,
we consider the case $X=\Sp(A)$, $Y=\Sp(B)$ and $Z=\Sp(C)$. By
hypothesis, we can suppose
$$
A\In C\In C_K\In B_K
$$
with $C$ integral over $A$. But, since $B$ is the integral closure
of $A$ in $B_K$, then $C\In B$. So we have
$$
A\In C\In B.
$$
These inclusions are functorial, so
 we have an
injective morphism of $\ox$-algebras $f_*(\oo{Z})\In g_*(\oo{Y})$.
We remark that $Y=\mathcal{S}pec(g_*\oo{Y})$ and
$Z=\mathcal{S}pec(f_*\oo{Z})$. This implies that there exists an
integral morphism  $h:Y\too Z$ such that $g=f \circ h$.

Now if $Z$ is normal then,  by definition of the integral closure
of $X$ in $Y_K$, clearly $f_*\oo{Z}=g_*\oo{Y}$ and so we have that
$h$ is an isomorphism.


\end{proof}
 We now prove a result of weak extension of $\Z/m\Z$-torsors. For
 any $m\in \N$ and any scheme $Z$ we define $_m Pic(Z):=\ker(Pic(Z)\on{m}\too
 Pic(Z))$.
\begin{prop}\label{lem:estensione di mu_p-torsori}Let $m\ge 1$ be an
integer.
Let $X$ be a normal faithfully flat scheme over $R$ with integral fibers and $_mPic(X_K)= {_m}Pic(X)$. 
  Let $f_K:Y_K\too X_K$ be a connected  $\Z/m\Z$-torsor.  Let $Y$ be
the normalization of $X$ in $Y_K$ and suppose that  $Y_k$ is
reduced. If $R$ contains a primitive $m$-{th} root of unity, there
exists a unique $\mu_{m}$-torsor $Y'$ over $X$ which
extends $f_K$. 
\end{prop}
\begin{rem} Let us consider the restriction map ${_m}Pic(X)\too { _m}Pic(X_K)$. If $X$ is a normal
faithfully flat scheme over $R$ with integral fibers, then the
above map is injective for any $m$ (see \cite[II 6.5]{Har1} and
\cite[7.2.14]{liu}; the hypothesis of separatedness cited in
\cite[II 6.5]{Har1} is not necessary). The above map is an
isomorphism, for any $m$, if, for instance, $X$ is also separated
and
 locally factorial (e.g. regular). See \cite[II 6.5, II 6.11]{Har1}.
\end{rem}

\begin{rem}
It has been proved by Epp(\cite{ep}) that if, for instance, $X$ is
of finite type over $R$ (or $X=\Sp(A)$ with $A$ the localization
or the $\pi$-adic completion of an $R$-algebra of finite type)
then, up to an extension of $R$, it is possible to suppose $Y_k$
reduced. Moreover the hypothesis $_mPic(X_K)= {_m}Pic(X)$ is
necessary. Indeed, to any
 $\sh\in \ _mPic(X_K)$, by the Kummer theory, we can associate (at least)  a $\mu_m$-torsor  $Y_K\too X_K$. It is clear that if
 $\sh$ is not extendible the same is true for the $\mu_m$-torsor $Y_K\too
 X_K$. 

\end{rem}
\begin{proof}
First, we consider the affine case $X=\Sp(A)$ and
$$Y_K=\Sp(A_K[Z]/(Z^{m}-f))$$ with $f\in A_K^*$. Since $Y\too X$ is a finite morphism, then
$Y=\Sp(B)$ for some normal and finite $A$-algebra
$B$. 
Multiplying $f$ by an $m$-{th} power of $\pi$  if necessary, which
does not change the $\mu_{m}$-torsor $Y_K\too X_K$, we can suppose
$f\in A$ and $f=\pi^n f_0$, with $0\le n<m$ and $f_0\not \in \pi
A$. This is possible since $A$ is integral and noetherian, so from
Krull's Theorem (\cite[1.3.13]{liu}) we have $\cap_{m\ge 0} \pi^m
A=0$.


We call $Y'=\Sp(A[Z]/(Z^{m}-f))$. We prove that $Y'$ is a
$\mu_{m}$-torsor over $X$, i.e. $f\in A^*$. Since $Y'$ is flat
over $R$ and $Y_K$ is connected and normal then $Y'$ is integral.
Moreover $Y'\too X$ is an integral morphism, so, by the previous
lemma, $Y'$ is dominated by $Y$. So $Z\in B$. Let us now suppose
$Z\in \pi B$.
 From $Z^{m}=f$ in $B$, it follows that
$ f\in \pi^m B.$ Applying  $m$-times \ref{lem:B/A flat} we have
$f\in \pi^m A$, which is a contradiction. So $Z\not\in \pi B$.
Since $Y_k$ is reduced, then $Z^m=f\not\in \pi B$. In particular,
$f\not\in \pi A$.
And since 
$f\in A_K^*$  there exists $g\in A\setminus \pi A$ and $l\in \N$
such that $f \frac{g}{\pi^l}=1$. So $fg=\pi^l$. But $X_k$ is
integral. Then $l=0$, which implies $f\in A^*$.

 Now by the Kummer theory, we associate
to any $\mu_m$-torsor over $X_K$ a line bundle $\clL$ over $X_K$
such that $\clL^{m}\simeq \oo{X_K}.$ Since $_mPic(X)=
{_m}Pic(X_K)$ we can assume that $\clL\in {_m}Pic(X)$. Then let
$\{U_i=\Sp(A_i)\}$ an affine covering of $X$ such that
$\clL_{|U_i}$ is trivial. If $\{g_{ij}\}\in H^1(X,\oo{X}^*)$
represents $\clL$,
 we have that there exists $f_i\in
H^0(U_{i,K},\oo{U_{i,K}}^*)$  such that
$(Y_K)_{U_i}=\Sp(A_{i,K}[T_i]/(T_i^{m}-f_i))$ and
$g_{ij}^{m}=\frac{f_i}{f_j}$. As seen before for any
$U_i=\Sp(A_i)$
 we can suppose $f_i\in A_i^*$. 
So
$\{Y_i'=\Sp(A_{i}[T_i]/(T_i^{m}-f_i))\}$ is a $\mu_m$-torsor which
extends the $\Z/m\Z$-torsor $Y_K\too X_K$. The uniqueness comes
from \ref{prop:iniezione tra h1}.
\end{proof}
\begin{rem}
 We remark that $Y$  does not usually coincide with $Y'$. This means that $Y'$  is
 possibly not normal.
\end{rem}
\begin{rem}\label{rem:weak estensione tipo Tp-f} From the first part of the proof also follows that if $X=\Sp(A)$
is affine, then the proposition remains true if we remove the
hypothesis ${_m}Pic(X)={_m}Pic(X_K)$ and we consider only
$\mu_m$-torsors of the type $A_K[T]/(T^m-f)$ with $f\in A_K^*$.
\end{rem}
\begin{cor}\label{cor:estensione torsori per gruppi commutativi}
Let $G$ be an abelian group of order $m$ and let us suppose that
$R$ contains a primitive $m$-th root of unity. Let $X$ be a normal
faithfully flat scheme over $R$ with integral fibers and
$_mPic(X)={ _m}Pic(X_K)$. Let us consider a connected $G$-torsor
$f_K:Y_K\too X_K$ and let $Y$ be the normalization of $X$ in
$Y_K$. Moreover, we assume that $Y_k$ is reduced. Then there
exists a (commutative) $R$-group-scheme $G'$ and a $G'$-torsor
$Y'\too X$ over $R$ which extends $f_K$.
\end{cor}
\begin{proof}
By the classification of abelian groups,  we have that
$G=\Z/{m_1}\Z\times \dots \times \Z/{m_r}\Z$ for some
$m_1,\dots,m_r\in \N$.   We remark that $R$ contains a primitive
$m_i$-{th} root of unity for $i=1,\dots,r$. So
\mbox{$(\Z/m_i\Z)_K\simeq (\mu_{m_i})_K$} for $i=1,\dots,r$.
Moreover, from hypothesis it follows, for $i=1,\dots,r$, that \mbox{$_{m_i}Pic(X)={ _{m_i}}Pic(X_K)$.} 
We firstly state the following lemma.
 \begin{lem}\label{lem:incollo torsori}Let $G_1,\dots, G_r$  be   flat group schemes over a scheme $X$.
 Let $Y_i\too X$ be a
$G_i$-torsor for any $i$. Then $\tilde{Y}=Y_1\times_X\dots \times_X
Y_r$ is a $G_1\times \dots \times G_r$-torsor, with the action
induced by those of $G_i$.
\end{lem}
\begin{proof}
We skip the proof which easily
 follows  by definition of torsor.

\end{proof}
We now come back to the proof of the corollary.  We  call
$G_i=\Z/{m_i}\Z$ for $i=1,\dots, r$.
 Moreover,
we call $\tilde{G}_i=G_1\times \dots\times \hat{G_i}\times \dots
\times G_{r}$. Let us define $(Y_i)_K=Y_K/(\tilde{G}_i)$, then
$(Y_i)_K\too X_K$ is a ${G}_i$-torsor. Moreover, $(Y_i)_K$ is
integral and normal. For any $i$, we call $\sigma_i$ the action of
$G_i$ induced by that of $G$ on $(Y_i)_K$. Hence,
$$
\sigma_i\times \id :G_i\times_{X_K} (Y_i)_K\too
(Y_i)_K\times_{X_K} (Y_i)_K
$$
is an isomorphism. By the above lemma, we have that
$(Y_1)_K\times_{X_K}\dots \times_{X_K} (Y_{r})_K$ is a $G$-torsor.
 Moreover,  the natural map
$$
q:Y_K\too (Y_1)_K\times_{X_K}\dots \times_{X_K} (Y_{r})_K
$$
preserves the $G$-actions; therefore it
is 
 a morphism of $G$-torsors. But, as it is well known, any morphism
of $G$-torsors is an isomorphism of schemes; hence, $q$ is an
isomorphism.


For $i=1,\dots, r$, we denote by $Y_i$ the normalization of $X$ in
$(Y_i)_K$. Then $Y_i$ is integral and normal for any $i$.
The projection $Y_K\simeq(Y_1)_K\times_K\dots \times_K(Y_r)_K\too
(Y_i)_K$ induces, by \ref{lem:universal property of normal
closure}, an integral morphism $Y\too Y_i$. Hence, we have, by
\ref{lem:B/A flat}, that $(Y_i)_k$ is reduced. So, by
\ref{lem:estensione di mu_p-torsori}, for any $i=1,\dots,r$, there
exists a $\mu_{{m_i}}$-torsor \mbox{$Y'_i\too X$} which extends
$(Y_i)_K\too X_K$.  Now let us consider
$Y'=Y'_1\times_X\dots\times_{X} Y'_{r}$. Using the above lemma
again, it follows that $Y'$ is a $\mu_{m_{1}}\times\dots \times
\mu_{m_{r}}$-torsor.
\end{proof}
\begin{rem}
Let $X$ be an integral scheme  faithfully flat over $R$ and $x$ an
$R$-point  of $X$. 
Let us consider the fundamental group schemes of Gasbarri
$\pi(X,x)$ over $R$ (see \cite{ga}) and $\pi(X_K,x_K)$. Then
Antei, in \cite{Mar},  proved that the natural morphism
$$
\phi:\pi(X_K,x_K)\too \pi(X,x)_K
$$
is a quotient morphism and that $\ker(\phi)=0$ if and only if  any
 $G_K$-torsor $Y_K\too X_K$ (with a $K$-section), such that the induced morphism
$\pi(X_K,x_K)\too G_K$ is a dominant morphism,  is weakly
extendible to a $G$-torsor $Y\too X$ over $R$ (with an
$R$-section), for some model $G$ of $G_K$. 
We stress that in that context you have also to extend  the
section of $Y_K$. But  if, for instance, $X$ is proper over $R$
hence $Y$ is proper over $R$ and the $K$-section of $Y_K$ can be
always extended.
\end{rem}
\begin{rem}\label{ex:non weak extension senza estendere}
 If $Y_k$ is not reduced we have an example of a torsor
not weakly extendible. For instance, take $X=\Sp(R[Z,1/Z])$ and
$Y_K=\Sp(K[Z,1/Z][T]/(T^p-\pi Z))$ as
$\Z/p\Z$-torsor over $X_K$. 
It is not too hard to see that $Y=\Sp(R[Z,1/Z][T]/(T^p-\pi Z))$ is
normal (see for example \cite[8.2.26]{liu}), so it is the
normalization of $X$ in $Y_K$. Moreover the action of
$\mu_p=\Sp(R[S]/(S^p-1))$ over $Y$ given by $T\mTo ST$ is clearly
faithful. So $\mu_p$ is the
effective model.  
Using \ref{lem:morfismo di gruppi deriva da quello di azione} it
follows that, if $Y_K\too X_K$ is weakly extendible by a
$G'$-torsor, then there is a model map $\mu_p\too G'$. Hence
$G'\simeq \mu_p$, because $\mu_p$ does not dominate any group
scheme except itself. We now claim that there is no $\mu_p$-torsor
$Y'\too
X$ which extends $Y_K\too X_K$. 
Since $\Pic(R[Z,1/Z])=0$,  we would have, by the Kummer theory,
$Y'=\Sp(R[Z,1/Z][T]/(T^p-f))$, for some $f\in R[Z,1/Z]^*$ such
that there exists $g\in K[Z,1/Z]^*$ with $fg^p=\pi Z$. But, since
$R[Z,1/Z]$ is factorial and $X_k$ is integral it is easy to see
that this not possible. 
In particular $Y_K\too X_K$ is
also not strongly extendible.
\end{rem}
We give here a proof of strong extension of  torsors under finite
abelian groups $G$ with $(|G|,p)=1$, for some schemes not
necessarily regular. For regular schemes refer to the
introduction.
\begin{cor}\label{cor:tame extension}Let $G$ be an  abelian group of order $m$ with  $(m,p)=1$ and let us suppose that $R$ contains a primitive $m$-th root of unity.
Let $X$ be a normal faithfully flat scheme over $R$ with  integral
fibers and \mbox{$_m Pic(X_K)={ _m}Pic(X)$.}
 Let $Y_K\too X_K$ be a connected $G$-torsor and
let $Y$ be the normalization of $X$ in $Y_K$. Moreover we assume
 that $Y_k$ is
reduced. Then $Y_K\too X_K$ is strongly extendible. 
\end{cor}
\begin{proof}
Indeed, by \ref{lem:estensione di mu_p-torsori}, there exists a
commutative finite flat $R$-group scheme $G'$ and a $G'$-torsor
$Y'\too X$ which extends $Y_K\too X_K$. Since the order of $G'$ is
coprime with $p$, then $G'$ is étale. So $Y'\too X$ is étale.
Since $X$ is normal, it follows from \cite[I 9.10]{SGA1} that $Y'$
is normal. Hence, from \ref{lem:universal property of normal
closure}, $Y'=Y$ and the proof is complete.
\end{proof}

We now conclude the section proving a lemma which will be
essential in the  sections that follow.


Let us consider  a normal faithfully flat  $R$-scheme $X=\Sp(A)$
with integral fibers. 
And we suppose that $\pi \in \mathcal{R}_A$, where $\mathcal{R}_A$
is the Jacobson radical. This condition means that the closed
points of $X$ are in the special fiber. Moreover, it is equivalent
to say that, for any $\lb\in R\setminus R^*$, any lifting of $a\in
(A/\lb A)^*$ is invertible in $A$. In particular, $A^*\too (A/\lb
A)^*$ is surjective.
From \eqref{eq:succ esatta lunga di g lambda} it follows that
$\hxg{\glb}\In \hxg{\gm}$. Then applying the snake lemma to the
following diagram (see \eqref{eq:succ esatta glbn})
\begin{equation}
\label{eq:diagramma succ esatta glbn} \xymatrix {\glb(A)
\ar[d]\ar[r]^{(\psi_{\lb,n})_*}&{\clG^{(\lb^{p^n})}(A)}
\ar[r]\ar[d]&
\hxg{G_{\lb,n}}\ar[r]^{i_*}\ar[d]&\hxg{{\glb}}\ar[r]\ar[d] &\hg{X}{\clG^{(\lb^{p^n})}}\ar[d]\\
A^* \ar[r]^{p^n}&A^* \ar[r]&
\hxg{\mu_{p^n}}\ar[r]^{i_*}&\hxg{{\gm}}\ar[r]& \hg{X}{\gm} }
\end{equation}
we have that the following commutative diagram 
$$
\xymatrix{\glbn(A)/{(\psi_{\lb,n})}_*(\glb(A))\ar[d]^{\alpha^{\lb^{p^n}}(A)}\ar[r]&H^1(X,\gln)\ar[d]\\
A^*/(A^*)^{p^n}\ar[r]&\hxg{\mu_{p^n}}.}
$$
is cartesian with $\alpha^{\lb^{p^n}}(A)$  injective.


\begin{rem}\label{lem:iniettività di una mappa}
 Explicitly the injection  $\alpha^{\lb^{p^n}}(A)$ means that, for any  $x\in A^*$,
$x^{p^n}\equiv 1 \mod \lb^{p^n}$ if and only if $x\equiv 1 \mod
\lb$.
\end{rem}
 \begin{lem}\label{lem:tecnico}
 Let $X$ be as above and let us suppose that $R$ contains a primitive $p^n$-th root of unity. Let $Y_K\too X_K$ be a $\mu_{p^n}$-torsor with $Y_K=\Sp(A_K[T]/(T^{p^n}-f))$, $f\in A^*$.
 Let us denote by $Y=\Sp(B)$ the normalization of $X$ in $Y_K$. We moreover suppose that $Y_k$ is reduced.
 Then, from \ref{lem:estensione di mu_p-torsori} and \ref{rem:weak estensione tipo Tp-f}, there exists a $\mu_{p^{n}}$-torsor $Y'\too X$ which weakly extends
 $Y_K\too X_K$.
 \begin{itemize}\item[i)]
 Using the filtration of \ref{cor:filtration},
    $[Y']\in \hxg{G_{\pi^{j},n}}$ if
 and only if there exists $g\in A^*$ such that
 \mbox{$fg^{p^n}=1+\pi^{jp^n}f_0$} for some $f_0\in A$.

\item[ii)] Let us suppose moreover $j<v(\lb_{(n)})$. If  $[Y']\in
\hxg{G_{\pi^{j},n}} \setminus \hxg{G_{\pi^{j+1},n}}$ and
$fg^{p^n}=1+\pi^{jp^n}f_0$, for some $f_0\in A$ and $g\in A^*$,
then $f_0$ is not a $p^n$-th power $\mod\pi$.
\end{itemize}
\end{lem}
\begin{proof}

We first prove (i). We observe that $[Y']\in H^1(X,\mu_{p^n})$ is
represented by $f\in A^*/{A^*}^{p^n}$. Since  the  above
commutative diagram is cartesian, then $[Y']\in
H^1(X,G_{\pi^j,n})$ if and only if $f\in \alpha^{\lb^{p^n}}(A)$.
And by definition, $f\in \alpha^{\lb^{p^n}}(A)$ if and only if
there exists $g\in A^*$
such that $fg^{p^n }=1+\pi^{jp^n}f_0$ for some $f_0\in A. $
We now prove the second statement.
  Let us suppose that  $[Y']\in
\hxg{G_{\pi^{j},n}}\setminus \hxg{G_{{\lambda_{(n)}},n}}$. We take
any  $g\in A^*$ such that $fg^{p^n}=1+\pi^{jp^n}f_0$,
 for some $f_0\in A$. If $[Y']\not\in  \hxg{G_{\pi^{j+1},n}}$, then,
 by (i),
$f_0\not \equiv 0\mod \pi^{p^n}$. In fact we will prove
$f_0\not\equiv 0\mod
 \pi$ in $A$. 
Since the torsor
$Y_1=\Sp(A[T]/(\frac{(1+\pi^jT)^{p^n}-1}{\pi^{jp^n}}-f_0))$,
associated with $[Y']\in  \hxg{G_{\pi^{j},n}}$, is integral over
$X$ and its generic fiber is isomorphic to $Y_K$, then, by
\ref{lem:universal property of normal closure}, the morphism
$Y\too X$ factors through  $Y_1$. Moreover, $Y\too Y_1$ is a
dominant morphism between integral affine schemes; hence $T\in
B\setminus \{0\}$. The fact that $f_0\nequiv 0\mod \pi^{p^n}A$
implies $T\nequiv 0\mod \pi B$. Otherwise, if $T=\pi T_0$ for some
$T_0\in B$, then $T^{p^n}\equiv 0\mod \pi^{p^n}B$. And, since
$j+1\le v(\lb_{(n)})$, we have
$$
f_0=\frac{(1+\pi^{j+1}T_0)^{p^n}-1}{\pi^{jp^n}}\equiv 0\mod
\pi^{p^n}B.
$$
So, by \ref{lem:B/A flat}, $f_0\equiv 0\mod \pi^{p^n}A$ against
the assumptions. Therefore, $T\nequiv 0\mod \pi B$.
 Now if $f_0 \equiv 0\mod \pi A$, then, since $j<v(\lb_{(n)})$, $T^{p^n}\equiv
0\mod \pi B$. But, as we just proved, $T\nequiv 0\mod \pi B$,
which contradicts the fact that $Y_k$ is reduced. So $f_0\nequiv
0\mod \pi A$.  We finally prove that $f_0$ is not a $p^n$-th power
$\mod \pi$. Indeed, if $f_0\equiv g_0^{p^n}\mod \pi$ for some
$g_0\in A\setminus \pi A$ then $f\equiv (1+\pi^j g_0)^{p^n}\mod
\pi^{jp^n+1}$. But $1+\pi^jg_0$ is invertible $\mod \pi$ then,
since $\pi \in \mathcal{R}_A$, $1+\pi^jg_0$ is invertible.
Multiplying $f$ by $(1+\pi^j g_0)^{-p^n}$, we can suppose $f\equiv
1 \mod \pi^{jp^n+1}$, which implies, by what we have just proved,
that $f\equiv 1 \mod \pi^{(j+1)p^n}$. Hence, by $i)$,  $ [Y']\in
\hxg{G_{\pi^{j+1},n}}$; this contradicts the hypothesis of
maximality of $j$.
\end{proof}

\section{Strong extension of $\Z/p\Z$-torsors}\label{sec:deg Z/pZ-torsori}
\stepcounter{subsection}\setcounter{subsection}{0}

 Let us suppose that $R$ contains  a primitive ${p}$-{th} root of
unity. We now suppose that  $X=\Sp(A)$ is a normal faithfully flat
$R$-scheme with integral fibers such that 
 $\pi\in \mathcal{R}_A$, where $\mathcal{R}_A$ is the Jacobson
radical of $A$, and $_pPic(X_K)=0$ (e.g. $A$ a local regular
faithfully flat $R$-algebra with integral fibers or $A$ a
factorial faithfully flat $R$-algebra complete with respect to the
$\pi$-adic topology and with integral fibers).
  We recall we suppose from the beginning that $X$ no\oe therian. Let us consider 
a nontrivial $\Z/p\Z$-torsor 
$$
Y_K\too X_K.
$$

  Let $Y$ be the normalization of $X$ in $Y_K$. We suppose that
$Y_k$ is reduced. There exists, by \ref{lem:estensione di
mu_p-torsori}, a
 unique
 $\mu_{p}$-torsor $Y'\too X$ such that $Y'_K\too X_K$ is isomorphic to
 the $\Z/p\Z$-torsor $Y_K\too X_K$.   So $Y_K\too X_K$
defines uniquely an
  element $[Y']\in H^1(X,\mu_p)$.

\begin{thm}\label{prop:degenerazione Z/pZ torsori}
 Notation as above.  Let us consider the filtration of
\ref{cor:filtration}. If $[Y_K] \in
\hxg{G_{{\pi^\gamma},1}}\setminus \hxg{G_{{{\pi^{\gamma+1}},1}}}$,
for some $0\le\gamma\le v(\lb_{(1)})$,
 then $Y$ is a $G_{{\pi^\gamma},1}$-torsor.
 In \mbox{particular} $Y_k$ is integral if $\gamma< v(\lb_{(1)})$. Moreover
 the valuation of
 the different of
 the extension $\oo{Y,\eta}/\oo{X,(\pi)}$ is
 $(p-1)(v(\lb_{(1)})-\gamma)$, for the generic point $\eta$ of any
 irreducible component of
 $Y_k$.
\end{thm}
\begin{rem}The trivial
$\Z/p\Z$-torsor over $X_K$ is  strongly extendible by the trivial
$\Z/p\Z$-torsor over $X$.
\end{rem}
\begin{proof}

Since $_pPic(A_K)=0$,  from the bottom line of \eqref{eq:diagramma
succ esatta glbn} with $v(\lb)=0$,
 it follows that $Y_K=\Sp(A_K[T]/(T^p-f))$ with $f\in
A_K^*$. 
From \ref{rem:weak estensione tipo Tp-f}, we can suppose $f\in
A^*$. Moreover, $Y'=\Sp(A[T]/{(T^p-f)})$ is
 the $\mup$-torsor  which extends $Y_K$. 

  If
$\gamma=0$, then $f$ is not a $p$-power $\mod \pi$ (otherwise, by
\ref{lem:tecnico}, up to a multiplication by a $p$-{th} power, we
can suppose $f\equiv 1\mod \pi$ and so $\gamma>0$). So $Y'$ is
normal by \ref{lem:normality's criterion}. Since $Y'$ is integral
over $X$ and its generic fiber is isomorphic to $Y_K$ then, from
\ref{lem:universal property of normal closure},  we have
 $Y\simeq Y'$.

If $\gamma= v(\lb_{(1)})$,
 then it is an étale torsor and the proof is complete.

If $v(\lb_{(1)})>\gamma>0$ we can suppose, by \ref{lem:tecnico},
$f=1+\pi^{p\gamma} f_0$
with $f_0\nequiv 0 \mod \pi$ in $A$. 
Let us consider the $G_{\pi^{\gamma},1}$-torsor
$$Y_1=\Sp(A[T]/(\frac{(1+\pi^\gamma T)^{p}-1}{\pi^{p\gamma}}-f_0)).
$$
By \ref{lem:tecnico} , $f_0\not\equiv g_0^p\mod \pi$ for any
$g_0\in A\setminus \pi A$. 
So $(Y_1)_k$ is reduced and, by \ref{lem:normality's criterion},
$Y_1$ is normal. Hence again from  \ref{lem:universal property of
normal closure}, we can conclude that $Y_1\simeq Y$. If
$\gamma<v(\lb_{(1)})$, then $\clG_k$ is radicial and $Y_k\too X_k$
is an inseparable morphism. Therefore, since $X_k$ is irreducible,
$Y_k$ is also irreducible.

The statement about the valuation of the different is clear.

\end{proof}

\begin{rem}\label{rem:estensione Z/pZ torsori senza fattoriale} The
theorem remains true if we remove the hypothesis $_pPic(A_K)=0$
 and we suppose that $Y_K=\Sp(A_K[T]/(T^p-f))$ with $f\in
A_K^*$. Indeed  $_pPic(A_K)=0$  only needs to ensure that any
$\mu_p$-torsor of $X_K$ is of the form $Y_K=\Sp(A_K[T]/(T^p-f))$
with $f\in A_K^*$.
\end{rem}

The following corollary will not used for the rest of the paper.

\begin{cor}
Let $G$ be an abelian group of order $m$ and let us suppose that
$R$ contains a primitive $m$-th root of unity. Let $X=\Sp(A)$ be a
normal faithfully flat
$R$-scheme with integral fibers such that 
 $\pi\in \mathcal{R}_A$ and $_m Pic(X_K)=0$. Let $h_K:Y_K\too X_K$ be a
connected $G$-torsor with $G$ an abelian group. Let $Y$ be the
normalization of $X$ in $Y_K$ and let us assume that $Y_k$ is
integral then $h:Y\too X$ is faithfully flat.

\end{cor}
\begin{proof}
 Since $G$ is abelian then $G\simeq \Z/{{m_1}}\Z\times \dots
\times \Z/{m_r}\Z\times \Z/{p^{m_{r+1}}\Z}\times \dots \times
\Z/{p^{m_s}\Z}$, with \mbox{$(m_i,p)=1$} for $i=1,\dots,r$. We
remark that any cyclic group is isomorphic, over $K$, to $
\mu_{l}$ for some $l$. So, as above, since $_m Pic(X_K)=0$, we
conclude that
$$
Y_K=\Sp(A_K[T_1,\dots,T_s]/(T_1^{{{m_1}}}-f_1,\dots,T_r^{{{m_r}}}-f_r,T_{r+1}^{p^{m_{r+1}}}-f_{r+1},\dots,
T_s^{p^{m_s}}-f_s))
$$
with $f_i\in A_K^*$. Let $n=r+m_{r+1}+\dots+m_{s}$. We have a
factorization
$$
h_K:(Y_n)_K:=Y_K\on{(h_{n})_K}\too (Y_{n-1})_K\on{(h_{n-1})_K}\too
\dots {\too} (Y_{1})_K\on{(h_{1})_K}{\too} X_K=:(Y_0)_K
$$
with $(h_i)_K$ a $\Z/p\Z$-torsor, if  \mbox{$i=1,\dots, n-r$} and
$(h_i)_K$ is a $\Z/m_{n-i+1}\Z$-torsor if \mbox{$i=n-r+1,\dots,
n$.} Any $(Y_i)_K$ is normal integral and affine. So if
$(Y_{i-1})_K=\Sp((B_{i-1})_K)$, it is easy to show that
$(Y_{i})_K=\Sp((B_{i-1})_K[T]/(T^{p}-f))$ with $f\in
(B_{i-1})_K^*$, if $i=1,\dots, n-r$ or
$(Y_{i})_K=\Sp((B_{i-1})_K[T]/(T^{m_{n-i+1}}-f))$ with $f\in
(B_{i-1})_K^*$, if \mbox{$i=n-r+1,\dots,n$.} We observe that not
necessarily $Pic((Y_i)_K)=0$.
Let $Y_i=\Sp(B_i)$ be the normalization of $X$ in $(Y_i)_K$. We
also obtain that $Y_i$ is the normalization of $Y_{i-i}$ in
$(Y_i)_K$. We have a factorization
$$
h:Y\on{h_{n}}\too Y_{n-1}\on{h_{n-1}}\too \dots {\too}
Y_{1}\on{h_{1}}{\too} X.
$$
Since $Y_k$ is integral, then $(Y_i)_k$ is integral for any $i$
(see \ref{lem:B/A flat}). 
 Moreover $h_i$ is finite, in particular it is closed. So  if
$\pi \in \mathcal{R}_{A}$ then $\pi \in \mathcal{R}_{B_i}$ also
for any $i$. From the above theorem, \ref{rem:estensione Z/pZ
torsori senza fattoriale}, \ref{cor:tame extension} and
\ref{rem:weak estensione tipo Tp-f} we have, for any $i$, that
$Y_i\too Y_{i-1}$ is a torsor under some finite flat group scheme.
In particular $h_i$ is faithfully flat for any $i$. Therefore, $h$
is faithfully flat.

\end{proof}

\
\section{Strong extension of $\Z/p^2\Z$-torsors}\label{sec:deg Z/p^2Z-torsors}
\subsection{Setup and degeneration types}Now  we
suppose that $R$ contains a primitive $p^2$-th root of unity.
Therefore we have \mbox{$(\Z/p^2\Z)_K\simeq (\mu_{p^2})_K$.} We
moreover
suppose $p>2$. 
Let $X:=\Sp A$ be  a normal  essentially semireflexive scheme
 over $R$ (see \S \ref{sec:effective models}) with integral fibers  such that $\pi \in \mathcal{R}_A$. We
also assume \mbox{$_{p^2} Pic(X_K)=0$.}
 Let $h_K:Y_K\too X_K$ be  a connected $\Z/{p^2\Z}$-torsor. Then we
consider the factorization
$$
Y_K\stackrel{(h_2)_K}{\too} (Y_1)_K\stackrel{(h_1)_K}{\too} X_K
$$
with both $(h_1)_K,(h_2)_K$ nontrivial $\Z/p\Z$-torsors.  Let
$Y_1=\Sp(B_1)$ be the normalization of $X$ in $(Y_1)_K$ and
$Y=\Sp(B)$ the normalization of $X$ in $Y_K$.  We  suppose that
$Y_k$ is integral.  Since $X$ is normal, by \ref{lem:universal
property of normal closure}, it follows that $Y$ is normal and
that $Y$ is the integral closure of $Y_1$ in $Y_K$. So we have the
factorization
$$
h:Y\on{h_2}\too Y_1\on{h_1}\too X
$$
with $h_1$ and $h_2$ degree $p$ morphisms.
Again by \ref{lem:universal property of normal closure}, it follows that  $Y_1$ is normal. 
By \ref{lem:B/A flat}, 
we have that, since $Y_k$ is integral, then $(Y_1)_k$ is integral
too. 

Since $_{p^2}Pic(A_K)=0$, we can suppose
$Y_K=\Sp(A_K[T]/(T^{p^2}-f))$ for some $f\in A_K^*\setminus
(A_K^*)^{p^2}$. By \ref{lem:estensione di mu_p-torsori},  we can
suppose $f\in A^*$. We can also write
$$
Y_K=\Sp(A[T_1,T_2]/(T_1^{p}-f,\frac{T_2^p}{T_1}-1)).
$$
 Therefore, we have
$$
(Y_1)_K=\Sp(A_K[T_1]/(T_1^{p}-f))
$$ and
\begin{equation}\label{eq:B=Tp-f} Y_K=\Sp((B_1)_K[T_2]/(\frac{T_2^{p}}{T_1}-1)).
\end{equation}
We remark that $B_K$ is finite and free as an $A_K$-module. In
particular it is semireflexive over $A_K$. From
\ref{lem:finite-->semireflexive} it follows that $Y$ is an
essentially semireflexive scheme over $\Sp(R)$. Therefore we can
apply \ref{lem:basta fedelta' su fibra speciale} to check if a group
scheme is an effective model for the $\Z/p^2\Z$-action on $Y$.
We now want to attach to the cover $Y_K\to X_K$ four invariants.
 We
have seen in the previous section that there exists an invariant,
which we called $\gamma$, that  is sufficient to solve the problem
of strong extension of $(\Z/p\Z)_K$-torsors. So the first two
invariants are simply the invariants $\gamma$ which arise from the
two $(\Z/p\Z)_K$-torsors $Y_K\too (Y_1)_K$ and $(Y_1)_K\too X_K$.
In precise by the above discussion it follows that $h_1$ satisfies
hypothesis of \S\ref{sec:deg Z/pZ-torsori}, hence we can apply
\ref{prop:degenerazione Z/pZ torsori}. Then, if we define
$\gamma_1\le v(\lb_{(1)})$ such that
$$[(Y_1)_K]\in\hxg{G_{{\pi^{\gamma_1}},1}}\setminus
\hxg{G_{{{\pi^{\gamma_1+1}},1}}},$$ it follows that $Y_1\too X$ is
a $G_{{\pi^{\gamma_1}},1}$-torsor.  From \ref{prop:degenerazione
Z/pZ torsori}, we have that $\gamma_1$ is also determined by the
valuation of the different $\mathcal{D}(h_1)$ of
$h_1:\Sp(\oo{Y_1,(\pi)})\too \Sp(\oo{X,(\pi)})$. We indeed have
$$
v(\mathcal{D}(h_1))=v(p)-(p-1)\gamma_1.
$$
From \ref{prop:degenerazione Z/pZ torsori}, we have that $(Y_1)_k$
is integral. Moreover, since $h_1$ is a closed morphism, then $\pi
\in {\mathcal{R}}_{B_1}$.
But not necessarily $_pPic((Y_1)_K)=0$. However, from
\eqref{eq:B=Tp-f} and the remark \ref{rem:estensione Z/pZ torsori
senza fattoriale} we can also apply \ref{prop:degenerazione Z/pZ
torsori}  to the morphism $h_2:Y\too Y_1$. Then, if we define
$\gamma_2\le v(\lb_{(1)})$ such that
$$[(Y)_K]\in H^1(Y_1,{G_{{\pi^{\gamma_2}},1}})\setminus
H^1(Y_1,{G_{{{\pi^{\gamma_2+1}},1}}}),$$ it follows that $Y\too
Y_1$ is a $G_{{{\pi^{\gamma_2}},1}}$-torsor. The invariant
$\gamma_2$ is determined by the different of
\mbox{$h_2:\Sp(\oo{Y,(\pi)})\too \Sp(\oo{Y_1,(\pi)})$,} too.
Indeed
$$
v(\mathcal{D}(h_2))=v(p)-(p-1)\gamma_2.
$$
 The third
invariant is linked to the filtration of \ref{cor:filtration}. It
is the integer $j\le v(\lb_{(2)})$ such that
$[Y_K]\in\hxg{G_{{\pi^j},2}}\setminus
\hxg{G_{{{\pi^{j+1}},2}}}$. 
We observe  that there exists
 a $G_{\pi^j,2}$-torsor $Y''$ which extends $Y_K\too X_K$. By \ref{lem:universal property of normal closure}, we have  morphisms
 $Y\too Y''$ and $Y_1\too Y''/G_{\pi^{j},1}$ such that the following
 diagram commutes

\begin{equation}\label{eq:Diagramma che commuta}
\xymatrix{Y\ar[d]\ar[rr]&\ar[r]&\ar[d] Y''\\
           Y_1 \ar[rr]\ar[rd]&\ar[r]    &\ar[ld]Y''/G_{\pi^{j},1}&\\
           &X&}
\end{equation}
\begin{lem}\label{lem:restriction degeneration type}
We have the following relations.
\begin{enumerate}
\item[i)]$ pj\le \gamma_1\le v(\lb_{(1)})$, 
\item[ii)] 
 $j\le \gamma_2\le v(\lb_{(1)})$,
\end{enumerate}

\end{lem}

\begin{proof}
By definition of $\gamma_i$, for $i=1,2$, we have $\gamma_i\le
 v(\lb_{(1)})$. We now prove the remaining statements. Let us consider  the diagram \eqref{eq:Diagramma che commuta}.
\begin{enumerate}
 \item [i)] We
 recall that $Y_1\too X$ is a $G_{\pi^{\gamma_1},1}$-torsor and $Y''/G_{\pi^{j},1}\too
 X$ is a $G_{\pi^{pj},1}$-torsor. So, by \ref{lem:morfismo di gruppi deriva da quello di azione}, we have a morphism $G_{\pi^{\gamma_1},1}\too G_{\pi^{pj},1}$. Therefore $\gamma_1\ge
 pj$. 
 \item[ii)] We recall that $Y\too Y_1$ is a $G_{\pi^{\gamma_2},1}$-torsor and $Y''\too Y''/G_{\pi^{j},1}$ is a $G_{\pi^{j},1}$-torsor.
 Again by \ref{lem:morfismo di gruppi deriva da quello di azione}, we have a morphism
 $G_{\pi^{\gamma_2},1}\too G_{\pi^{j},1}$. Therefore, $\gamma_2\ge j$.  
\end{enumerate}
\end{proof}

 By definition of $j$, up to a multiplication of $f$ by an element of ${(A^{*})}^{p^2}$, which
does not change the $\mu_{p^2}$-torsor on the generic fiber,
we can suppose $f=1+\pi^{p^2j}f_0$ with $f_0\in A$. And, if $j
<v(\lb_{(2)})$, by \ref{lem:tecnico} $f_0$ is not a $p^2$-th power
$\mod \pi$.

Before  introducing the last invariant, we describe explicitly the
scheme $Y$. By definition of $\gamma_1$ and by the proof of
\ref{prop:degenerazione Z/pZ torsori}, there exists $g\in A^*$
(not unique) such that $fg^{-p}=1+\pi^{p\gamma_1}f_1$ with $f_1\in
A$.  
Let us consider $\Sp(B_1)$ with
$$B_1=A[T_1]/(\frac{(1+\pi^{\gamma_1}T_1)^p-1}{\pi^{p\gamma_1}}-f_1).
$$
If $\gamma_{1}=v(\lb_{(1)})$ then $\Sp(B_1) \too X$ is étale, then
$\Sp(B_1)$ is normal and hence, by \ref{lem:universal property of
normal closure}, $Y_1=\Sp(B_1)$. While if $\gamma_1<v(\lb_{(1)})$,
by \ref{lem:tecnico},  $f_1$ is not a $p$-{th} power $\mod \pi$.
Then, by \ref{lem:normality's criterion}, $\Sp(B_1)$ is normal and
hence again $Y=\Sp(B_1)$.

If ${Y_1'}\too Y_1$ is the $\mu_p$-torsor which extends $Y_K\too
(Y_1)_K$, then  $Y$ is the normalization of $Y_1$ in
$$
({Y}_1')_K=\Sp((B_1)_K[T_2]/(\frac{T_2^p}{1+\pi^{\gamma_1 }T_1}-g)).
$$
Then,  reasoning as above, there exists $H(T_1)\in B_1^*$, such that
$$
g(1+\pi^{\gamma_1} T_1)(H(T_1))^{-p}\equiv 1 \mod
\pi^{p\gamma_2}B_1,
$$
and $Y=\Sp(B)$ with
$$
B=B_1[T_2]/\bigg(\frac{(1+\pi^{\gamma_2}T_2)^p-1}{\pi^{p\gamma_2}}-\frac{gH(T_1)^{-p}(1+\pi^{\gamma_1}T_1)-1}{\pi^{p\gamma_2}}\bigg).
$$


We remark that the definition of $g$ and $H(T_1)$ depends on the
choice of the representative $f$ of  $[Y_K]\in \hxg{\mu_{p^2}}$
and they are not uniquely determined. We now see how they vary as
$f$ varies. We stress that we require $f\equiv 1 \mod \pi^{p^2j}$.
Let us substitute $a^{p^2}f$ to $f$, with $a\in A^*$ and
$a^{p^2}f\equiv 1\mod \pi^{p^2j}$. It follows from
\ref{lem:iniettività di una mappa} that $a^{p^2}f\equiv 1\mod
\pi^{p^2j}$ is equivalent to $a\equiv 1\mod \pi^{j}$. Now it is
immediate to see that we have to substitute $a^p g$ to $g$  and
$aH(T_1)$ to $H(T_1)$. We now prove that, for a fixed $f$, the
elements $g$ and $H(T_1)$ are uniquely determined in a certain
sense.
 \begin{lem}\label{lem:r=j} 
 Notation as  above. Let us fix a  representative $f=1+\pi^{p^2j}f_0$ of $[Y_K]\in
\hxg{\mu_{p^2}}$.  We have the following results. 
\begin{enumerate}
\item  The element $g$ is uniquely determined $\mod
\pi^{\gamma_1}$. 
Moreover, we can suppose
 \begin{gather*}
g=1+\pi^{jp}g_0 \text{ with } g_0\not \in \pi A.
\end{gather*}
\item The element $H(T_1)$ is
uniquely determined $\mod\pi^{\gamma_2}$. Any $H(T_1)$ as above is
of the form
$$
H(T_1)=1+\pi^{j} H_1(T_1) \text{ with } H_1(T_1) \not\in  \pi B_1.
$$
if $0\le j<\gamma_2$, and we can suppose $H(T_1)=1$ if
$j=\gamma_2$.
Finally, if $j>0$,  up to a change of the representative $f$, we
can suppose $H(0)=1$.
\end{enumerate}
\end{lem}
\begin{proof}
\begin{enumerate}
\item
We firstly prove the uniqueness of $g\mod \pi^{\gamma_1}$. By
definition of $g$ we have $g^p\equiv f\mod \pi^{p\gamma_1}$.  Let
us take  $g'\in A^*$. Then  $g'^p\equiv g^p\equiv f\mod
\pi^{p\gamma_1}$ if and only if
$$
(\frac{g}{g'})^p\equiv 1\mod\pi^{p\gamma_1}
$$
if and only if, from \ref{lem:iniettività di una mappa},
$$
g\equiv g'\mod \pi^{\gamma_1}.
$$
This proves the uniqueness of $g\mod \pi^{\gamma_1}$. In fact we
proved something more. Indeed, since $g$ is determined by the
property $g^p\equiv f\mod \pi^{\pi^{p\gamma_1}}$, we have also
proved that we can replace $g$ with any $g'\in A^*$ such that
$g'\equiv g\mod\pi^{\gamma_1}$. (Clearly,  $H(T_1)$ will also be
different).

We now prove that, up to a change of $g$ by $g(1+\pi^{\gamma_1}
h)$ with $h\in A$, we can suppose
 \begin{gather*}
g=1+\pi^{jp}g_0 \text{ with } g_0\not \equiv 0\mod \pi A.
\end{gather*}

Since $f=1+\pi^{p^2j}f_0$ and $pj\le \gamma_1$ (see
\ref{lem:restriction degeneration type}), we have
$$
g^{p}\equiv 1\mod \pi^{p^2j}.
$$
Again by \ref{lem:iniettività di una mappa}, we obtain
$$
g=1+\pi^{pj}g_0
$$
for some $g_0\in A$.  If $pj<\gamma_1$ then $g_0\not\in \pi A$,
otherwise, since $g^p\equiv 1+\pi^{p^2j}f_0\mod \pi^{p\gamma_1}$,
we would have, again by \ref{lem:iniettività di una mappa},
$f_0\equiv 0\mod \pi$, against hypothesis on $f_0$.  While, if
$pj=\gamma_1$ and $g_0\equiv 0\mod \pi$, by what proved before, we
can replace $g$ with $g'=g(1+\pi^{pj}h)$, with $h\not \in \pi A$.
So $g'=1+\pi^{pj}(g_0+h+\pi^{pj}g_0 h)$ with
$g_0+h+\pi^{pj}g_0h\equiv h\not \equiv 0\mod \pi$.
\item As above, we can prove that $H(T_1)$ is unique $\mod
\pi^{\gamma_2}$ and that we can replace it with any
$\tilde{H}(T_1)$ such that $\tilde{H}(T_1)\equiv H(T_1)\mod
\pi^{\gamma_2}$. And since $g\equiv 1 \mod \pi^{pj}B_1$ and
$\gamma_1\ge pj$ (see \ref{lem:restriction degeneration type})
then, $H(T_1)^p\equiv 1\mod \pi^{jp}B_1$. Hence, as above, we can
conclude
 $H(T_1)=1+\pi^jH_1(T_1)$. We now prove that
 $H_1(T_1)\not \equiv 0\mod \pi B_1$ if  $j<\gamma_2$ and that we can assume $H(T_1)=1$ if $j=\gamma_2$. If $j<\gamma_2$, then from
 $$
(1+\pi^{j} H_1(T_1))^p\equiv (1+\pi^{pj} g_0)(1+\pi^{\gamma_1}T_1)
\mod \pi^{p\gamma_2}
 $$
it follows that, if $H_1(T_1)\equiv 0\mod \pi B_1$, then, using
\ref{lem:iniettività di una mappa}, we have
$g_0+\pi^{\gamma_1-pj}T_1\equiv 0\mod\pi B_1$. Since $B_1$ is
finite and free over $R$ it follows that $g_0\equiv 0\mod \pi$,
against what just proved in $1$. Now if $\gamma_2=j$ then we can
take $H(T_1)=1$ which clearly satisfies the condition
$$
H(T_1)^p\equiv (1+\pi^{pj}g_0)(1+\pi^{\gamma_1}T_1)\mod
\pi^{p\gamma_2}.
$$
We now prove the last statement. Let $j>0$. From what we just
proved we know that $H(0)\equiv 1 \mod \pi^{j} A$. Since $\pi\in
\mathcal{R}_A$, then  $H(0)$ is invertible. 
If we  change $f$ into $f{H}(0)^{-p^2}$, from the discussion
before the
lemma, we have  
to replace 
${H}(T_1)$ with $\frac{{H}(T_1)}{H(0)}$. So the proof is complete.

\end{enumerate}
\end{proof}

Now, given $H(T_1)=\sum_{k=0}^{p-1}a_k{T_1}^k\in B_1^*$,  let us
consider $H'(T_1)$ as its formal derivative. Using the above
lemma, we suppose $a_0=1$ if $j>0$. For any $m\ge \gamma_1$, we
will say that $a\in \pi R$ satisfies $(\triangle)_m$ if
$$
aH(T_1)\equiv  \pi^{m-\gamma_1}H'(T_1)\mod \pi^{\gamma_2}.
$$
We finally give the definition of the fourth invariant.
\begin{defn} We will call \textit{effective threshold} the number
$$
\kappa=\min\{m\ge \gamma_1| \exists a\in \pi R
 \text{
which satisfies } (\triangle)_{m} \}.
$$
\end{defn}If we take $m\ge \gamma_1+\gamma_2$ and $a=0$, we see that
such a minimum exists.
\begin{lem}\label{rem:unicità di a}For any $m\ge \kappa$ there exists a unique solution, $\mod \pi^{\gamma_2}$, of
$ (\triangle)_m$. We will call  $\alpha_m\in \pi R$ any of its
lifting. If $H(0)=a_0\equiv 0 \mod \pi A$ then $\alpha_m\equiv 0\mod
\pi^{\gamma_2}$.
\end{lem}
\begin{rem}
By \ref{lem:r=j} it follows that the case $H(0)\equiv 0\mod\pi$
can  happen only if $j=0$.
\end{rem}
\begin{proof}
Let us firstly suppose $a_0\not\equiv 0\mod \pi A$.
If $b_i$, for $i=1,2$, is  solution of  $(\triangle)_m$, 
it follows that for any $m\ge \gamma_1$ we have in particular $b_i
a_0\equiv  \pi^{m-\gamma_1}a_1\mod\pi^{\gamma_2}$. Therefore
$$
a_0(b_1-b_2)\equiv 0\mod\pi^{\gamma_2}.
$$
But $a_0\not\in \pi A$, $X$ is flat over $R$ and $X_k$ is
integral; therefore $b_1\equiv b_2\mod\pi^{\gamma_2}$.

 We  now consider the case $a_0\in \pi A$. Since
$H(T)\in B_1^*$ and $a_0\in \pi A$, then  there exists $0<i\le
p-1$ such that $a_i\not \in \pi A$. Let $\bar{i}$ be the least
integer with this property.
Let $a$ be a solution solution of $(\triangle)_m$ and suppose that
$a\not\equiv 0 \mod\pi^{\gamma_2}$. In particular
$$
a a_{\bar{i}}\equiv (\bar{i}+1)a_{\bar{i}+1}\pi^{m-\gamma_1}\mod
\pi^{\gamma_2}
$$
and
$$
a a_{\bar{i}-1}\equiv \bar{i}a_{\bar{i}}\pi^{m-\gamma_1}\mod
\pi^{\gamma_2}.
$$
Therefore,  by the minimality of $\bar{i}$ and by the fact that
$a\not\equiv 0\mod \pi^{\gamma_2}$,
\begin{gather*}\label{eq:v(a)>v(mu/mu')}
v(a)\ge {m-\gamma_1}>v(a)
\end{gather*}
which is a contradiction.  Therefore, if $a_0\in \pi A$, then
$a\equiv
 0\mod\pi^{\gamma_2}$.\end{proof}


\begin{defn}Using the previous notation, we say that the
degeneration type of $Y_K\too X_K$ is
$(j,\gamma_1,\gamma_2,\kappa)$.
 \end{defn}
\begin{lem}\label{lem:restriction degeneration type II}
We have
$ \gamma_1\le \kappa\le \gamma_1+\gamma_2-j$. In particular
$\gamma_2=j$ implies $\kappa=\gamma_1$.
\end{lem}

\begin{proof}
By \ref{lem:r=j}, $H'(T)\equiv 0\mod \pi^j$. Therefore, if we take
$m=\gamma_1+\gamma_2-j$, then
$$
\pi^{m-\gamma_1}H'(T)\equiv 0\mod \pi^{\gamma_2}.
$$
Therefore, $a=0$ satisfies $(\triangle)_m$. This implies
$\kappa\le \gamma_1+\gamma_2-j$. Now, if $\gamma_2=j$, then
$\kappa \le \gamma_1$. But, by definition of $\kappa$, we have
$\kappa \ge \gamma_1$. Hence $\kappa=\gamma_1$.
\end{proof}
\subsection{The main theorem}
We here prove the main theorem of the paper.
\begin{thm}\label{teoremone}
Let us suppose that $R$ contains a primitive $p^2$-th root of
unity and that
 $p>2$. 
Let $X:=\Sp A$ be  a normal  essentially semireflexive scheme
 over $R$  with integral fibers  such that $\pi \in \mathcal{R}_A$. We
moreover assume \mbox{$_{p^2} Pic(X_K)=0$.} Let $Y_K\too
X_K$ be a connected $\Z/p^2 \Z$-torsor and 
 $Y$ be the normalization of $X$ in $Y_K$.
 Let us suppose that $Y_k$ is integral. 
If $Y_K$ has $(j,\gamma_1,\gamma_2,\kappa)$ as degeneration type,
then its effective model is
$$\clE^{(\pi^{\kappa},\pi^{\gamma_2};E_p(\alpha_\kappa S),1)}.$$
Moreover, if $\alpha_\kappa\not\equiv 0\mod\pi^{\gamma_2}$, then
$v(\alpha_\kappa)=\kappa-\gamma_1+j$. Otherwise
$\kappa-\gamma_1+j= \gamma_2$.
\end{thm}

\begin{proof}
As we proved in the previous subsection   $Y=\Sp(B)$ with
$$
B=A[T_1,T_2]/\bigg(\frac{(1+\pi^{\gamma_1}T_1)^p-1}{\pi^{p\gamma_1}}-f_1,\frac{(1+\pi^{\gamma_2}T_2)^p-1}{\pi^{p\gamma_2}}-\frac{g
H(T_1)^{-p}(1+\pi^{\gamma_1}T_1)-1}{\pi^{p\gamma_2}}\bigg).
$$
By the definition of integral closure of $X$ in $Y_K$,
 the $\Z/p^2{\Z}$-action on $Y_K$ can be
extended to an action on $Y$. We now explicitly describe this
action. If we set
$$
\eta_{\pi}=\frac{\pi^{v(\lb_{(1)})}}{\lb_{(1)}}\eta=\frac{\pi^{v(\lb_{(1)})}}{\lb_{(1)}}\sum_{k=1}^{p-1}\frac{(-1)^{k-1}}{k}\lb_{(2)}^k
$$
then we can write, by \ref{ex:equazione per eta_pi} and
\ref{lem:abbasso valutazione con blow-up},
\begin{equation}\label{eq:Z/p^2 Z con etapi}
\Z/p^2\Z=\Sp(A[S_1,S_2]/\bigg(\frac{(1+\pi^{v(\lb_{(1)})}
S_1)^p-1}{\pi^{pv(\lb_{(1)})}},\frac{\frac{(E_p(\eta_{\pi}S_1)+\pi^{v(\lb_{(1)})}
S_2)^p}{1+\pi^{v(\lb_{(1)})} S_1}-1}{\pi^{pv(\lb_{(1)})}}\bigg)).
\end{equation}
Since $Y_K$ is a $\mu_{p^2}$- torsor,  on the generic fiber, the
action is given by
\begin{align*}
&1+\pi^{\gamma_1} T_1\longmapsto (1+\pi^{v(\lb_{(1)})}S_1)(1+\pi^{\gamma_2}T_1)\\
&(1+\pi^{\gamma_2}T_2)H(T_1)\longmapsto (E_p(\eta_\pi
S_1)+\pi^{v(\lb_{(1)})}S_2)(1+\pi^{\gamma_2}T_2){H(T_1)},
\end{align*}
so  it is globally given by
\begin{align*}
&T_1\longmapsto \pi^{v(\lb_{(1)})-\gamma_1}S_1+T_1+\pi^{v(\lb_{(1)})}S_1T_1\\
&T_2\longmapsto \frac{(E_p(\eta_\pi
S_1)+\pi^{v(\lb_{(1)})}S_2)\bigg(\frac{(1+\pi^{\gamma_2}T_2){H(T_1)}}{H({\pi^{v(\lb_{(1)})-\gamma_1}}S_1+T_1+\pi^{v(\lb_{(1)})}S_1T_1)}\bigg)-1}{\pi^{\gamma_2}}
\end{align*}
The proof of the theorem is obtained as a consequence of several
lemmas.

\begin{lem}\label{lem:modelli effettivi di Z/p^2Z sono estensioni}
If an effective model $\clG$ for the action of  $\Z/p^2\Z$ exists,
then it is   of the form $\clE^{(\pi^m,\pi^{\gamma_2};F,1)}$  
with $v(\lb_{(1)})\ge m \ge \max\{\gamma_2,\gamma_1\}$.
\end{lem}
 \begin{proof}
Since the effective model is in particular a model of
$(\Z/p^2\Z)_K$, it follows, by \ref{cor:clE se lb divide mu}, that
the effective model $\clG$ is of the form
$\clE^{(\pi^m,\pi^{\gamma_2};F,1)}$ with $v(\lb_{(1)})\ge m\ge
\gamma_2$. Moreover, $\clG/G_{\pi^{\gamma_2},1}\simeq G_{\pi^m,1}$
has an $X$-action over $Y_1$. But $Y_1\too X$  is a
$G_{\pi^{\gamma_1},1}$-torsor. So, by \ref{lem:morfismo di gruppi
deriva da quello di azione}, we have a model map
$G_{\pi^{m},1}\too
G_{\pi^{\gamma_1},1}$. 
Then $m\ge \gamma_1$. 
\end{proof}
Let us now consider a group scheme of type
$\clE^{(\pi^m,\pi^{\gamma_2},F,1)}$. We consider the normalization
map $\phi:\Z/p^2\Z\too \clE^{(\pi^m,\pi^{\gamma_2},F,1)}$.
We  give necessary and sufficient conditions to have an action of
$\clE^{(\pi^m,\pi^{\gamma_2};F,1)}$ on $Y$ compatible with $\phi$.
By \ref{cor:clE se lb divide mu}, we have that
$$F(S)=E_p(a S)\in
((R/\pi^{\gamma_2}R)[S]/(\frac{(1+\pi^{\gamma_1}S)^p-1}{\pi^{p\gamma_1}}))^*$$
for some $a\in R/\pi^{\gamma_2} R$. In the following, we take a
lifting $\tilde{a}\in R$ of $a\in R/\pi^{\gamma_2}R$ and consider
\mbox{$\tilde{F}(S)=\sum_{i=0}^{p-1}\frac{\tilde{a}^i}{i!}S^i\in
R[S]$} as a lifting of  $F(S)$.

\begin{lem} \label{lem:cond suff e nec per avere azione.}
There exists an action of $\clE^{(\pi^m,\pi^{\gamma_2};F,1)}$ on $Y$
compatible with $\phi$ if and only if
\begin{gather*}
\tilde{F}(S)H(T)-H(\pi^{m-\gamma_1}S+T+\pi^m ST)\equiv 0 \mod
\pi^{\gamma_2}
\end{gather*}
\end{lem}
\begin{proof}Let us suppose that such an action exists.
Reasoning as above, it is possible to show that the action is given
by \begin{align*}
&T_1\longmapsto {\pi^{m-\gamma_1}}S_1+T_1+\pi^m S_1T_1\\
 &T_2\longmapsto
\frac{(\tilde{F}(S_1)+\pi^{\gamma_2}
S_2)\bigg(\frac{(1+\pi^{\gamma_2}T_2)H(T_1)}{H(\pi^{m-\gamma_1}S_1+T_1+\pi^m
S_1T_1)}\bigg)-1}{\pi^{\gamma_2}}
\end{align*}

 Then, in particular,  
  $\frac{(\tilde{F}(S_1)+\pi^{\gamma_2}
S_2)\bigg(\frac{(1+\pi^{\gamma_2}T_2)H(T_1)}{H(\pi^{m-\gamma_1}S_1+T_1+\pi^m
S_1T_1)}\bigg)-1}{\pi^{\gamma_2}} $ belongs to
$$
B\pt A[S_1,S_2]/\bigg(\frac{(1+\pi^m
S_1)^p-1}{\pi^{mp}},\frac{\frac{(\tilde{F}(S_1)+\pi^{\gamma_2}
S_2)^p}{1+\pi^m S_1}-1}{\pi^{p\gamma_2}}\bigg)
$$
So 
we have
\begin{equation}\label{eq:azione su Y}
\begin{aligned}
&\frac{(\tilde{F}(S_1)+\pi^{\gamma_2}
S_2)\bigg(\frac{(1+\pi^{\gamma_2}T_2)H(T_1)}{H(\pi^{m-\gamma_1}S_1+T_1+\pi^m
S_1T_1)}\bigg)-1}{\pi^{\gamma_2}} =\\
& \frac{\tilde{F}(S_1)H(T_1)-H(\pi^{m-\gamma_1}S_1+T_1+\pi^m
S_1T_1)}{\pi^{\gamma_2} H(\pi^{m-\gamma_1}S_1+T_1+\pi^m
S_1T_1)}+T_2\frac{\tilde{F}(S_1)H(T_1)}{H(\pi^{m-\gamma_1}S_1+T_1+\pi^m
S_1T_1)}+\\
&+
S_2\frac{(1+\pi^{\gamma_2}T_2)H(T_1)}{H(\pi^{m-\gamma_1}S_1+T_1+\pi^m
S_1T_1)}.
\end{aligned}
\end{equation}
This implies
\begin{gather*}
\tilde{F}(S_1)H(T_1)-H(\pi^{m-\gamma_1}S_1+T_1+\pi^m S_1T_1)\equiv 0 \mod \pi^{\gamma_2}.\\
\end{gather*}
But it is clear that this  condition is also sufficient to define
the required action.
\end{proof}
The next  lemma, 
together with \ref{lem:cond suff e nec per avere azione.}, links the
definition of the effective threshold with the existence of an
action of a model of $\Z/p^2\Z$ on $Y$.

\begin{lem}\label{lem:H(ct)=F(T)}
Let $\tilde{b}\in \pi R$. 
Let us consider
$\tilde{G}(S)=\sum_{i=0}^{p-1}\frac{\tilde{b}^i}{i!}S^i\in R[S]$. 
The following statements are equivalent.
\begin{itemize}
 \item[(i)] $\tilde{G}(S)H(T)\equiv H(\pi^{m-\gamma_1}S+T+\pi^{m} ST) \mod\pi^{\gamma_2}$;
 \item[(ii)]$\tilde{b} H(T)\equiv \pi^{m-\gamma_1} H'(T)\mod\pi^{\gamma_2}$,   
where $H'$ is the formal derivative of $H$.
\end{itemize}
Moreover, they imply the following assertions.
\begin{enumerate}
\item Let us consider $R[G_{\pi^m,1}]=R[S]/(\frac{(1+\pi^m
S)^p-1}{\pi^{mp}})$. Then
$$\tilde{G}(S)\in
\Hom_{gr}({G_{\pi^m,1}}_{|S_{\pi^{\gamma_2}}},{\gm}_{|S_{\pi^{\gamma_2}}})$$
and
$$
\tilde{G}(S)^p\equiv 1+\pi^m S\mod \pi^{p\gamma_2}
R[G_{\pi^{m},1}].
$$
\item If $m>\gamma_1$ then
$$\frac{\tilde{G}(S)H(T)-H(\pi^{m-\gamma_1}S+T+\pi^m
ST)}{\pi^{\gamma_2}}\equiv \frac{\tilde{b}
H(T)-\pi^{m-\gamma_1}H'(T)}{\pi^{\gamma_2}}\mod\pi
$$
\end{enumerate}
\end{lem}
\begin{rem}\label{rem:soluzione equazione diff} We assume the theorem in this remark.
Let us suppose   $\tilde{b} H(T)\equiv H'(T)\mod\pi^{\gamma_2}$.
In particular, by definition, $\kappa=\gamma_1$. And from
\ref{teoremone}, we also have $\kappa\ge \gamma_2$. We remark that
if $j=0$, then, again from \ref{teoremone}, it follows that
$\gamma_2=j=0$. So in any case, from \ref{lem:r=j}($\emph{2}$)
 we can suppose $H(0)=1$. So if we consider only the
constant term, as polynomials in $T$, of the equality $(i)$, we
obtain
$$
H(S)
\equiv\tilde{G}(S)\equiv\sum_{i=0}^{p-1}\frac{\tilde{b}^i}{i!}S^i\mod
\pi^{\gamma_2}.
$$
Moreover, from \ref{lem:H(ct)=F(T)}$(i)$, it follows that we can
think $H(S)$ as an element of
$\Hom_{gr}({\gmx{1}}_{|S_\lb},{\gm}_{|S_\lb})$. Then from
\ref{lem:suriettività mappa tra hom}, it follows that
$\tilde{b}^p\equiv 0\mod\pi^{\gamma_2}$.
\end{rem}
\begin{proof}
$(i)\Rightarrow (ii)$. Let us suppose
$$
\tilde{G}(S)H(T)\equiv H(\pi^{m-\gamma_1}S+T+\pi^m ST)\mod
\pi^{\gamma_2}.
$$
We consider both members as polynomials in $S$ with coefficients
in $ R [T]$. Then, if we compare the coefficients of $S$, we
obtain $(ii)$.

$(i)\Leftarrow (ii)$. Let $H^{(k)}(T)$ denote the $k^{th}$ formal
derivative of $H(T)$. We remark that $(i)$ is equivalent to
\begin{equation*}
\tilde{b}^k H(T)\equiv (\pi^{m-\gamma_1})^{k}H^{(k)}(T)\mod
\pi^{\gamma_2}
\end{equation*}
for $1\le k\le p-1$.  We prove a little more. We prove  that
\begin{equation}\label{eq:per j}
\tilde{b}^k H(T)\equiv (\pi^{m-\gamma_1})^{k}H^{(k)}(T)\mod
\pi^{{\gamma_2}+\min\{(k-1)v(\tilde{b}),(k-1)(m-\gamma_1)\}}
\end{equation}
For $k=1$ it is obvious.
 Let us now suppose it is true for $k$, we
prove it for $k+1$. If we multiply \eqref{eq:per j} by
$\tilde{b}$, we obtain
\begin{equation}\label{eq:moltiplicato per b}
\tilde{b}^{k+1} H(T)\equiv
\tilde{b}(\pi^{m-\gamma_1})^{k}H^{(k)}(T)\mod
\pi^{{\gamma_2}+\min\{(k-1)v(\tilde{b}),(k-1)(m-\gamma_1)\}+v(\tilde{b})}.
\end{equation}
Moreover, if we differentiate the equation (ii) $k$ times,  we
obtain
\begin{equation}\label{eq:derivato j volte}
\tilde{b}H^{(k)}(T)\equiv \pi^{m-\gamma_1}H^{(k+1)}(T)\mod
\pi^{\gamma_2}.
\end{equation}
Multiplying \eqref{eq:derivato j volte} by $\pi^{k(m-\gamma_1)}$,
we obtain
\begin{equation}\label{eq:derivato j volte per pi}
\tilde{b}\pi^{k(m-\gamma_1)}H^{(k)}(T)\equiv
(\pi^{m-\gamma_1})^{(k+1)}H^{(k+1)}(T)\mod
\pi^{\gamma_2+k(m-\gamma_1)}.
\end{equation}
 Then \eqref{eq:moltiplicato per b} and
\eqref{eq:derivato j volte per pi} give
$$
\tilde{b}^{k+1} H(T)\equiv (\pi^{m-\gamma_1})^{k+1}H^{(k+1)}(T)\mod
\pi^{\gamma_2+\min\{kv(\tilde{b}),k(m-\gamma_1)\}}.
$$
as we required. So $(i)$ and  $\emph{2}$ are proved.
%
 Let us now suppose (i)  true.
\begin{enumerate}

 \item[\emph{1}.] We recall that $H(T)=\sum_{i=0}^{p-1}a_iT^i\in (A[T]/(\frac{(1+\pi^{\gamma_1}
T)^p-1}{{\pi^{\gamma_1}}^p}-f_1))^*=B_1^*$. If $a_0\in \pi A$ then
$\tilde{b}\equiv 0\mod\pi^{\gamma_2}$ by \ref{rem:unicità di a}.
 Let us now suppose that $a_0\not\in \pi A$. We  think $H(S)\in A[S]/(\frac{(1+\pi^{\gamma_1}
S)^p-1}{{\pi^{\gamma_1}}^p})=R[G_{\pi^{\gamma_1},1}]\pt_R A$. We
consider the morphism $\psi_{\pi^m,\pi^{\gamma_1}}:G_{\pi^m,1}\too
G_{\pi^{\gamma_1},1}$, given by $S\mapsto \pi^{m-\gamma_1}S$. Then
$$\psi_{\pi^m,\pi^{\gamma_1}}^*(H(S))= H(\pi^{m-\gamma_1}S)$$
However if we compare the coefficients of $T$ in (i) we obtain
\begin{equation}\label{eq:exp} H(\pi^{m-\gamma_1}S)\equiv
a_0\tilde{G}(S)\mod\pi^{\gamma_2}.
\end{equation}
Therefore,
$$
\psi_{\pi^m,\pi^{\gamma_1}}^*(H(S))\equiv
a_0\tilde{G}(S)\mod\pi^{\gamma_2}.
$$
 Let us
now consider $\id \times
\psi_{\pi^m,\pi^{\gamma_1}}:G_{\pi^m,1}\times G_{\pi^{m},1}\too
G_{\pi^m,1}\times G_{\pi^{\gamma_1},1}$. Hence, if we apply $\id
\times \psi_{\pi^m,\pi^{\gamma_1}}^*$ to $(i)$ we obtain, using
\eqref{eq:exp},
$$
a_0 {\tilde{G}}(S)\tilde{G}(T)\equiv H(\pi^{m-\gamma_1}(S+T+\pi^m
ST))\equiv a_0\tilde{G}(S+T+\pi^m ST)\mod\pi^{\gamma_2}
$$
which implies, since $a_0\not \in \pi A$ and $G_{\pi^m,1}\times
G_{\pi^m,1}$ is flat over $A$,
$$
 {\tilde{G}}(S)\tilde{G}(T)\equiv  \tilde{G}(S+T+\pi^m ST)\mod\pi^{\gamma_2}.
$$
 This means ${\tilde{G}}(S)\in
\Hom_{gr}({G_{\pi^m,1}}_{|S_{\pi^{\gamma_2}}},{\gm}_{|S_{\pi^{\gamma_2}}})$.
Moreover, we know that
\begin{equation}\label{eq:Hp=g (1+pil T)}H(T)^p\equiv g(1+\pi^{\gamma_1} T)\mod
\pi^{p\gamma_2} B_1.
\end{equation} Hence,
$$
(H(T)\tilde{G}(S))^p\equiv g(1+\pi^{\gamma_1} T) \tilde{G}(S)^p
\mod \pi^{p\gamma_2} (R[G_{\pi^{m},1}]\pt_R B_1).
$$
Moreover, it is easy to see that, since $m\ge \gamma_1$,
$$
 R[G_{\pi^{m},1}]\pt_R B_1= R[S,T]/\bigg(\frac{((1+\pi^{\gamma_1}T)(1+\pi^m
S))^{p}-1}{\pi^{p{\gamma_1} }}-f_1,\frac{(1+\pi^m
S)^p-1}{\pi^{mp}}\bigg).
$$
Then we can substitute  $\frac{(1+\pi^{\gamma_1}T)(1+\pi^m
S)-1}{\pi^{\gamma_1}}$ to $T$ in \eqref{eq:Hp=g (1+pil T)} and we
obtain
$$
(H(\pi^{m-\gamma_1}S+T+\pi^m ST))^p\equiv g (1+\pi^m
S)(1+\pi^{\gamma_1} T)\mod \pi^{p\gamma_2} (R[G_{\pi^{m},1}]\pt_R
B_1).
$$
By hypothesis we have that
$$
\tilde{G}(S)H(T)\equiv H(\pi^{m-\gamma_1}S+T+\pi^m ST)\mod
\pi^{\gamma_2}(R[G_{\pi^{m},1}]\pt_R B_1)
$$
and therefore, using \eqref{eq:Hp=g (1+pil T)},
$$
g(1+\pi^{\gamma_1} T) \tilde{G}(S)^p\equiv g (1+\pi^m
S)(1+\pi^{\gamma_1} T)\mod \pi^{p\gamma_2}(R[G_{\pi^{m,1}}]\pt_R
B_1).
$$
This implies
$$
\tilde{G}(S)^p\equiv (1+\pi^m S)\mod \pi^{p\gamma_2}
R[G_{\pi^{m,1}}].
$$


\end{enumerate}
\end{proof}



%
%
%
We are now able to find a candidate to be the effective model.
\begin{lem}\label{lem:restrictions on gamma2}
If an effective model for the $\Z/p^2\Z$-action exists, it must be
the group scheme
$\clE^{(\pi^{\kappa},\pi^{\gamma_2};E_p(\alpha_{\kappa} S),1)}$.
In particular, we must have $\alpha_\kappa \in
\Phi_{\pi^{\kappa},\pi^{\gamma_2}}^1$. Moreover $\gamma_2\le
\kappa\le v(\lb_{(1)})$.
\end{lem}
\begin{proof}Since, as we have seen, $\Z/p^2\Z$ acts on $Y$ then, by
\ref{lem:cond suff e nec per avere azione.} and the previous
lemma,
 it follows that $\eta_\pi$
satisfies $(\triangle)_{v(\lb_{(1)})}$. Therefore, $\kappa\le
v(\lb_{(1)})$.

By \ref{lem:modelli effettivi di Z/p^2Z sono estensioni} it follows
that the effective model is of the form
$\clE^{(\pi^m,\pi^{\gamma_2};F,1)}$ for some $m\le v(\lb_{(1)})$ and
$F\in
\Hom_{gr}({G_{\pi^{m},1}}_{|S_{\pi^{\gamma_2}}},{\gm}_{|S_{\pi^{\gamma_2}}})$.
By \ref{lem:cond suff e nec per avere azione.} and
\ref{lem:H(ct)=F(T)}, 
we have that if  a group scheme
$\clE^{(\pi^m,\pi^{\gamma_2};F,1)}$ acts on $Y$ then
$F=E_p(\alpha_m S)$ with $\alpha_m\in \pi R$ which satisfies
$(\triangle)_m$. Conversely, if $m\le v(\lb_{(1)})$ and
$\alpha_m\in \pi R$ satisfies $(\triangle)_m$, then by
\ref{lem:cond suff e nec per avere azione.},
\ref{lem:H(ct)=F(T)}($\emph{1}$) and \ref{def:clE(mu,lb,F,i)}  we
can construct the group scheme
$\clE^{(\pi^{m},\pi^{\gamma_2};E_p(\alpha_m S),1)}$ and it
acts on $Y$ compatible with the action of $\Z/p^2\Z$. 
We remark that by \ref{rem:unicità di a} the equation
$(\triangle)_m$ has (unique)
solution if and only if $m\ge \kappa$. 
Moreover, for any $v(\lb_{(1)})\ge m'\ge m$ there exists a model
map $\clE^{(\pi^{m'},\pi^{\gamma_2};E_p(\alpha_{m'} S),1)}\too
\clE^{(\pi^{m},\pi^{\gamma_2};E_p(\alpha_m S),1)}$. Indeed, by
definition of $m$, we have that there exists $\alpha_m\in \pi R$
such that
$$
\alpha_m H(T)\equiv\pi^{m-\gamma_1}H'(T)\mod \pi^{\gamma_2}.
$$
Therefore
$$
\pi^{m'-m}\alpha_m H(T)\equiv\pi^{m'-\gamma_1}H'(T)\mod
\pi^{\gamma_2}.
$$
But we know that
$$
\alpha_{m'} H(T)\equiv\pi^{m'-\gamma_1}H'(T)\mod \pi^{\gamma_2}.
$$
And, as seen in \ref{rem:unicità di a}, the solution of the above
equation is  unique $\mod \pi^{\gamma_2}$. Therefore,
$\pi^{m'-m}\alpha_m\equiv \alpha_{m'}\mod\pi^{\gamma_2}$. So,
 by \ref{lem:abbasso valutazione con blow-up}, there exists a model
map
$$
\clE^{(\pi^{m'},\pi^{\gamma_2};E_p(\alpha_{m'} S),1)}\too
\clE^{(\pi^{\kappa},\pi^{\gamma_2};E_p(\alpha_{m}S),1)}.
$$
 We recall
that for any $m\ge \kappa$ the action of
$\clE^{(\pi^{m'},\pi^{\gamma_2};E_p(\alpha_{m'} S),1)}$ is given by
\begin{align*}
T_1\longmapsto & {\pi^{m-\gamma_1}}S_1+T_1+\pi^m S_1T_1\\
 T_2\longmapsto &
 \frac{F(S_1)H(T_1)-H(\pi^{m-\gamma_1}S_1+T_1+\pi^m
S_1T_1)}{\pi^{\gamma_2} H(\pi^{m-\gamma_1}S_1+T_1 +\pi^m
S_1T_1)}+\\&+T_2\frac{F(S_1)H(T_1)}{H(\pi^{m-\gamma_1}S_1+T_1+\pi^m
S_1T_1)}+
S_2\frac{(1+\pi^{\gamma_2}T_2)H(T_1)}{H(\pi^{m-\gamma_1}S_1+T_1+\pi^m
S_1T_1)}
\end{align*}
The above model map is compatible with the actions on $Y$. In
particular, we have,
 for any $v(\lb_{(1)})\ge m> \kappa$, a model
map
$$
\clE^{(\pi^{m},\pi^{\gamma_2};E_p(\alpha_{m}S),1)}\too
\clE^{(\pi^{\kappa},\pi^{\gamma_2};E_p(\alpha_{\kappa}S),1)}.
$$
compatible with the actions. Since the above model map is not an
isomorphism, there is a nontrivial kernel $\tilde{H}$ of the
morphism restricted to the special fiber. Since the map is
compatible with the actions, then  $\tilde{H}\In
(\clE^{(\pi^{m},\pi^{\gamma_2};E_p(\alpha_{m}S),1)})_k$ acts
trivially on $Y_k$. So
$$
\clE^{(\pi^{m},\pi^{\gamma_2};E_p(\alpha_{m}S),1)}
$$
is not the effective model of the $\Z/p^2\Z$-action if $m>\kappa$.
Hence, if an effective model exists, it must be
$\clE^{(\pi^{\kappa},\pi^{\gamma_2};E_p(\alpha_{\kappa}S),1)}.$
Since the group
$\clE^{(\pi^{\kappa},\pi^{\gamma_2};E_p(\alpha_{\kappa}S),1)}$
exists, it follows, by \ref{cor:clE se lb divide mu}, that
$\kappa\ge \gamma_2$ and $\alpha_{\kappa}\in
\Phi_{\pi^\kappa,\pi^{\gamma_2}}^1.$
\end{proof}

We remark that if $X$ was of finite type over $R$ then $Y$ would
be of finite type over $R$. So applying the theorem of existence
of effective models \ref{teo:modelli effettivi caso algebrico} we
would have finished.
 We now prove that
 $\clE^{(\pi^{\kappa},\pi^{\gamma_2};E_p(\alpha_{\kappa}S),1)}$ is
 the effective model for the action of $\Z/p^2\Z$ in the general
 case.
By construction the action is faithful
on the generic fiber. 
Since $X$ is  essentially semireflexive over $R$ it is sufficient,
by \ref{lem:basta fedelta' su fibra speciale}, to check the
faithfulness on the special fiber. Let us suppose that the map
$$
\clG_k=(\clE^{(\pi^\kappa,\pi^{\gamma_2};E_p(\alpha_\kappa
S),1)})_k\too Aut_k(Y_k)
$$
has nontrivial kernel $\tilde{K}$. Since the action of
$(G_{\pi^{\gamma_2},1})_k$ on $Y_k$ is faithful by definition of
$\gamma_2$, then $\tilde{K}\times_{\clG_k}
(G_{\pi^{\gamma_2},1})_k$ is the trivial group scheme. Therefore,
$\tilde{K}$ is a group scheme of order $p$ and
$$
(\clE^{(\pi^{\kappa},\pi^{\gamma_2};E_p(\alpha_\kappa
S),1)})_k\simeq (G_{\pi^{\gamma_2},1})_k\times_k \tilde{K}
$$
 We  distinguish two cases.

\begin{tabular}{|c|}
  \hline
  $\kappa=\gamma_1$. \\
  \hline
\end{tabular}
Since $Y\too Y_1$ is a $G_{\pi^{\gamma_2},1}$-torsor and
$\tilde{K}$ acts trivially on $Y_k$, we have
\begin{equation}\label{eq:X=Y/cle}(Y_1)_k\simeq Y_k/(G_{\pi^{\gamma_2},1})_k\simeq
Y_k/((G_{\pi^{\gamma_2},1})_k\times_k \tilde{K})\simeq
Y_k/(\clE^{(\pi^{\kappa},\pi^{\gamma_2};E_p(\alpha_\kappa
S),1)})_k.
\end{equation} 
But by definition of $\gamma_1$, $Y_1\too X$ is a
$G_{\pi^{\gamma_1},1}$-torsor.
 So, using the fact that
$\kappa=\gamma_1$,
$$X_k\simeq (Y_1)_k/(G_{\pi^{\gamma_1},1})_k\simeq \bigg(Y_k/(G_{\pi^{\gamma_2},1})_k\bigg)/\bigg((\clE^{(\pi^{\kappa},\pi^{\gamma_2};E_p(\alpha_\kappa
S),1)})_k/(G_{\pi^{\gamma_2},1})_k\bigg) \simeq
Y_k/(\clE^{(\pi^{\kappa},\pi^{\gamma_2};E_p(\alpha_\kappa
S),1)})_k,$$ which contradicts \eqref{eq:X=Y/cle}, since $X_k\neq
(Y_1)_k$.


\begin{tabular}{|c|}
  \hline
  $\kappa>\gamma_1$. \\
  \hline
\end{tabular} 
We remark that necessarily $\gamma_2>0$. Indeed, if $\gamma_2=0$,
then by \ref{lem:restriction degeneration type}(ii) and
\ref{lem:restriction degeneration type II} necessarily
$\kappa=\gamma_1$. It is also clear that $\kappa>0$.
Now,  
from \ref{lem:H(ct)=F(T)}(2) and \eqref{eq:azione su Y}, it
follows that the action on the special fiber is given by the
reduction $\mod\pi$ of
\begin{align*}
&T_1\longmapsto T_1\\
 &T_2\longmapsto
\frac{\alpha_\kappa
H(T_1)-\pi^{\kappa-\gamma_1}H'(T_1)}{\pi^{\gamma_2}H(T_1)}S_1+T_2 +
S_2
\end{align*}
We now prove that
\begin{equation}\label{eq:no sottogruppi stabilzzatori}
\frac{\alpha_\kappa H(T_1)-
\pi^{\kappa-\gamma_1}H'(T_1)}{\pi^{\gamma_2}H(T_1)}S_1\not \equiv
bS_1 \mod \pi \end{equation} for any $b\in R$. Let us suppose
$\frac{\alpha_\kappa H(T_1)-
\pi^{\kappa-\gamma_1}H'(T_1)}{\pi^{\gamma_2}H(T_1)}S_1\equiv bS_1
\mod \pi$ with
$b\in R$. 
Then
$$\frac{\alpha_{\kappa} H(T_1)-
\pi^{\kappa-\gamma_1}H'(T_1)}{\pi^{\gamma_2}}S_1\equiv
bH(T_1)S_1\mod\pi$$ with $b\in R$. Therefore,
$$
\frac{(\alpha_{\kappa}-b\pi^{\gamma_2}) H(T_1)-
\pi^{\kappa-\gamma_1}H'(T_1)}{\pi^{\gamma_2}}\equiv 0\mod\pi.
$$
It clearly   follows that
$$
\frac{\alpha_\kappa-b\pi^{\gamma_2}}{\pi}H(T_1)\equiv
\pi^{\kappa-1-\gamma_1}H'(T_1)\mod \pi^{\gamma_2}
$$
Then $\alpha_{\kappa-1}=\frac{\alpha_\kappa-b\pi^{\gamma_2}}{\pi}$
satisfies $(\triangle)_{\kappa-1}$; it  is easy to see that this
implies $\alpha_{\kappa-1}\in \pi R$. The  minimality of $\kappa$ is
contradicted. So we have proved \eqref{eq:no sottogruppi
stabilzzatori}.

We now consider three different cases.
 If $\gamma_2,\kappa<v(\lb_{(1)})$, then $$(\clE^{(\pi^{\kappa},\pi^{\gamma_2};E_p(\alpha_\kappa
S),1)})_k\simeq \alpha_p\times_k
\alpha_p=\Sp(k[S_1,S_2]/(S_1^p,S_2^p)).$$ Its subgroups  of order
$p$ different from $(G_{\pi^{\gamma_2},1})_k$ are  the subgroups
$S_2+b S_1=0 $ with $b\in k$. If
\mbox{$\gamma_2<\kappa=v(\lb_{(1)})$,} then
$$(\clE^{(\pi^{\kappa},\pi^{\gamma_2};E_p(\alpha_\kappa
S),1)})_k\simeq \alpha_p\times_k
\Z/p\Z=\Sp(k[S_1,S_2]/(S_1^p-S_1,S_2^p))$$ and the only subgroup
isomorphic to $\tilde{K}\simeq \Z/p\Z$ is $S_2=0$. Finally, if
$\gamma_2=\kappa=v(\lb_{(1)})$, then
$$(\clE^{(\pi^{\kappa},\pi^{\gamma_2};E_p(\alpha_\kappa
S),1)})_k\simeq \Z/p\Z\times_k
\Z/p\Z=\Sp(k[S_1,S_2]/(S_1^p-S_1,S_2^p-S_2))$$ and the only
subgroups isomorphic to $\tilde{K}\simeq \Z/p\Z$ different from
$(G_{\pi^{\gamma_2},1})_k$ are the subgroups $S_2+b S_1=0$ with
$b\in \F_p$. In any case, by \eqref{eq:no sottogruppi
stabilzzatori},
 the action
restricted to any subgroup of
$(\clE^{(\pi^{\kappa},\pi^{\gamma_2};E_p(\alpha_\kappa S),1)})_k$ is
not trivial.

%
 We now
prove the last sentence of the theorem. We have, by definition,
\begin{equation}\label{eq:equazione di alpha kappa}
\alpha_\kappa H(T_1)\equiv \pi^{\kappa-\gamma_1}H'(T_1)\mod
\pi^{\gamma_2}
\end{equation}
 Moreover, $H(T_1)\in B_1^*$ and, if we consider $H'(T_1)\in
{ B_1}_{(\pi)}$, we have
\begin{equation}\label{eq:v(H')=j}v(H'(T_1))=j,\end{equation} by \ref{lem:r=j}.
If $\alpha_\kappa\equiv 0\mod\pi^{\gamma_2}$, then by
\eqref{eq:equazione di alpha kappa} and \eqref{eq:v(H')=j}, it
follows that $\pi^{\kappa-\gamma_1+j}\equiv 0\mod\pi^{\gamma_2}$.
Therefore, $\kappa-\gamma_1+j\ge \gamma_2$. So, by
\ref{lem:restriction degeneration type II}, we have
$\kappa-\gamma_1+j= \gamma_2$. While, if $\alpha_\kappa\not\equiv
0\mod\pi^{\gamma_2}$, it follows again from \eqref{eq:equazione di
alpha kappa} and \eqref{eq:v(H')=j} that
$v(\alpha_k)=\kappa-\gamma_1+j$. The theorem is proved.
\end{proof}
We here give a criterion to determine when $Y$ has a structure of
torsor.
\begin{cor}\label{cor:criterio torsore}
Let us suppose we are in the hypothesis of the theorem. Then
$Y\too X$ is a $G$-torsor under some finite and flat group scheme
$G$ if and only if $\kappa=\gamma_1$. Moreover, $\kappa=\gamma_1$
if and only if  $\gamma_1\ge \gamma_2$ and $H(T)\equiv E_p(a
T)\mod \pi^{\gamma_2}$, for some $a\in \pi R$ such that $a^p\equiv
0\mod\pi^{\gamma_2}$. In such a case,
$G=\clE^{(\pi^{\gamma_1},\pi^{\gamma_2};H,1)}$. 
\end{cor}
\begin{rem}\label{rem:degeneration type of G torsors}The degeneration type of any  $\clE^{(\pi^{\gamma_1},{\pi^{\gamma_2}};E_p(\alpha_\kappa
S),1)}$-torsor is
$$(v(\alpha_\kappa),\gamma_1,\gamma_2,\gamma_1),$$
if $\alpha_\kappa\not\equiv 0\mod \pi^{\gamma_2}$ and
$$
(\gamma_2,\gamma_1,\gamma_2,\gamma_1)
$$
otherwise.  This follows from \ref{teoremone} and
\ref{cor:criterio torsore}.
\end{rem}
\begin{proof}
  We remarked  in \ref{rem:modelli effettivi e torsori} that if $Y\too X$ is a $G$-torsor for some
finite and flat group scheme  then $G$ must coincide with the
effective model $\clG$ of $\Z/p^2\Z$ acting on $Y$. In other
words, $Y\too X$ is a $G$-torsor  if and only if it is a
$\clG$-torsor. By the theorem we have that the effective model for
the $\Z/p^2\Z$-action is
$\clE^{(\pi^{\kappa},\pi^{\gamma_2};E_p(\alpha_\kappa S),1)}$.
Moreover, there is the following exact sequence
$$
0\too G_{\pi^{\gamma_2},1}\on{i}{\too} \cal{G}\on{p}{\too}
G_{\pi^{\kappa},1} \too  0
$$
By \ref{prop:proprieta modelli effettivi} (i), we have that
$G_{\pi^{\gamma_2},1}$  is the effective model of the action of
$\Z/p\Z\In \Z/p^2\Z$ on $Y$. Now if $Y\too X$ is a $\clG$-torsor,
then it satisfies the hypothesis of \ref{prop:proprieta modelli
effettivi} (iii), then $G_{\pi^{\kappa},1}$ is the effective model
of the action of $\Z/p^2\Z/\Z/p \Z$ on $Y_1$. But by the
definition of $\gamma_1$, we have that $Y_1\too X$ is a
$G_{\pi^{\gamma_1},1}$-torsor.  Then, again by \ref{rem:modelli
effettivi e torsori}, we have $G_{\pi^{\kappa},1}\simeq
G_{\pi^{\gamma_1},1}$, which implies $\kappa=\gamma_1$.

Let us now suppose that $\kappa=\gamma_1$. We recall that
$$
Y=\Sp(A[T_1,T_2]/(\frac{(1+\pi^{\gamma_1}
T_1)^p-1}{\pi^{p\gamma_1}}-f_1,\frac{(1+\pi^{\gamma_2}T_2)^p-1}{\pi^{p\gamma_2}}-\frac{g
H(T_1)^{-p}(1+\pi^{\gamma_1}T_1)-1}{\pi^{p\gamma_2}}))
$$
Moreover, by definition of $\kappa$, we have $\alpha_\kappa
H(T)\equiv H'(T)\mod\pi^{\gamma_2}$. Then, by \ref{rem:soluzione
equazione diff}, it follows that \mbox{$H(T)\equiv
E_p(\alpha_\kappa
T)\mod \pi^{\gamma_2}$.} 
Now let us substitute ${T_2}{H(T_1)^{-1}}$ to $T_2$. Then we
obtain
$$
Y=\Sp(A[T_1,T_2]/(\frac{(1+\pi^{\gamma_1}
T_1)^p-1}{\pi^{p\gamma_1}}-f_1,\frac{(H(T_1)+\pi^{\gamma_2}T_2)^p(1+\pi^{\gamma_1}T_1)^{-1}-g}{\pi^{p\gamma_2}}))
$$
By definition of
$\clE^{(\pi^{\gamma_1},\pi^{\gamma_2};E_p(\alpha_{\kappa}S),1)}$
 there exists
$G\in\Hom_{gr}({\clG^{(\pi^{\gamma_1})}}_{|S_{\pi^{p\gamma_2}}},{\gm}_{|S_{\pi^{p\gamma_2}}})$
such that
$$E_p(\alpha_{\kappa}S)^p(1+\pi^{\gamma_1}
S)^{-1}=G(\frac{(1+\pi^{\gamma_1} S)^p-1}{\pi^{p\gamma_1}})\in
\Hom_{gr}({\clG^{(\pi^{\gamma_1})}}_{|S_{\pi^{p\gamma_2}}},{\gm}_{|S_{{\pi^{p\gamma_2}}}}).
$$
We remark that, if we think $E_p(\alpha_{\kappa} T_1),G(T_1)\in
B_1^*$, the previous equation gives
$$E_p(\alpha_{\kappa}T_1)^p(1+\pi^{\gamma_1} T_1)^{-1}\equiv
G(f_1)\mod \pi^{p\gamma_2}B_1.$$  However, we have that
$H(T_1)^p(1+\pi^{\gamma_1} T_1)^{-1}\equiv
E_p(\alpha_{\kappa}T_1)^p(1+\pi^{\gamma_1} T_1)^{-1}\equiv g\mod
\pi^{p\gamma_2}B_1$. Therefore, using \ref{lem:B/A flat},
$$
g\equiv G(f_1)\mod \pi^{p\gamma_2} A
$$
i.e. $g=G(f_1)+\pi^{p\gamma_2}f_2$ for some $f_2\in A$. Hence, by
\S\ref{sec:clE torsors}, $Y\too X$ is a
$\clE^{(\pi^{\gamma_1},\pi^{\gamma_2};E_p(\alpha_{\kappa}S),1)}$-torsor.

We now have, by definition of $\kappa$, that $\kappa=\gamma_1$ if
and only if there exists $\alpha_\kappa\in \pi R$ such that
\begin{equation}\label{eq:(*)kappa}
\alpha_\kappa H(T_1)\equiv H'(T_1)\mod\pi^{\gamma_2}
\end{equation} We remark that, since $\kappa\ge \gamma_2$,
$\kappa=\gamma_1$ only if $\gamma_1\ge \gamma_2$.  In such a case,
by \ref{rem:soluzione equazione diff}, $H(T_1)$ satisfies
\eqref{eq:(*)kappa} if and only if there exist $\alpha_\kappa\in
\pi R$ such that $\alpha_\kappa^p\equiv 0\mod\pi^{\gamma_2}$ and
$H(T_1)\equiv E_p(\alpha_\kappa T_1)\mod \pi^{\gamma_2}$.

\end{proof}
In particular we obtain the following result.
\begin{cor}\label{cor:cond.suff per no torsore}If $\gamma_1<\gamma_2$ (or equivalently $v(\clD(h_1))>v(\clD(h_2))$), then $Y\too X$ has no structure of torsor.
\end{cor}
\begin{rem}Unfortunately we have no example of coverings with
$\gamma_1<\gamma_2$. So we don't know if this case can really occur.
\end{rem}
\begin{ex}\label{ex:non strongly extendible torsor} We here give an example, for any $p\ge 3$,
where $Y\too X$ is not a $\clG$-torsor under any group scheme
$\clG$. Notation is as above. We moreover suppose that there
exists $f_1\in A_k^*\setminus {A^*_k}^{p}$ such that
$A_k\setminus A_k^p[f_1]\neq \emptyset$. 
 Take $\gamma_1,\gamma_2$ such that
$v(p)>p\gamma_1>p^2\gamma_2>0$. In particular we have
$v(p)>(p-1)\gamma_1+p\gamma_2$. Let us  take  liftings in $A$ of
$f_1 \in A_k^*\setminus {A^*_k}^{p}$
 and $f_2\in A_k\setminus A_k^p[f_1]$
which we will again denote by  $f_1$ and $f_2$.
 Moreover, let
us consider $g=f_1+\pi^{p\gamma_2}f_2\in A^*$. We consider the
$\Z/p^2\Z$-torsor $Y_K\too X_K$ with
$$ Y_K\simeq \Sp(A_K[T]/(T^{p^2}-(1+\pi^{p\gamma_1} f_1)g^{p})).$$
For instance, we can take $A=R[Z]_{(\pi,Z)}$, $\gamma_1=p+1$,
$\gamma_2=1$, $f_1=1+Z^2$, $f_2=Z$, \mbox{$g=1+Z^2+\pi^p Z$.}


Then we consider $\Sp(B_1)=\Sp(A[T_1]/(\frac{(1+\pi^{\gamma_1}
T_1)^p-1}{\pi^{p\gamma_1}}-f_1))$. Since $f_1$ is not a $p$-{th}
power $\mod\pi$ then $\Sp(B_1)_k$ is integral and $\Sp(B_1)$ is
normal (see \ref{lem:normality's criterion}). So $Y_1=\Sp(B_1)$.
We remark that by hypothesis we have that $T_1^p\equiv f_1\mod
\pi^{p\gamma_2+1}$. We now take $H(T_1)=T_1\in B_1^*$. Then we
have, by construction,
$$ \frac{g (1+\pi^{\gamma_1}T_1)-H(T_1)^p}
{\pi^{p\gamma_2}}\equiv f_2\mod\pi.$$ So we consider
$$
\Sp(B_1[T_2]/(\frac{(1+\pi^{\gamma_2}
T_2)^p-1}{\pi^{{p\gamma_2}}}-\frac{H(T_1)^{-p}g
(1+\pi^{\gamma_1}T_1)-1}{\pi^{p\gamma_2}}))
$$
Hence
$$
\Sp(B_k)=\Sp(A[T_1,T_2]/(T_1^p-f_1,T_2^p-\frac{f_2}{f_1}))
$$
Since $(B_1)_k^p=A_k^p[f_1]$ and $f_2\not\in A_k^p[f_1]$, then
 $\Sp(B_k)$ is integral, therefore $\Sp(B)$ is normal.  
Hence, \mbox{ $Y=\Sp(B)$.} The degeneration type of $Y_K\too X_K$
is $(0,\gamma_1,\gamma_2,\gamma_1+\gamma_2)$. Indeed $H'(T_1)=1$,
so
$$
a T_1\equiv \pi^{\kappa-\gamma_1} \mod\pi^{\gamma_2}
$$
if and only if $a\equiv 0\mod\pi^{\gamma_2}$ and
$\kappa-\gamma_1\ge \gamma_2$.  Since $\kappa\le
\gamma_1+\gamma_2-j$   this means $\kappa =\gamma_1+\gamma_2$ and
$j=0$. The effective model is
$$
\clG=\clE^{(\pi^{\gamma_1+\gamma_2},\pi^{\gamma_2};1,1)}.
$$
Since $\kappa=\gamma_1+\gamma_2>\gamma_1$ then $Y$ is not a
$\clG$-torsor by \ref{cor:criterio torsore}.

\end{ex}

\section{Realization of degeneration types}
\stepcounter{subsection}\setcounter{subsection}{0}
 We have shown in the above section that the degeneration type has
to satisfy some restrictions. We here want to study the problem of
determining  the elements of $\N^4$ which can be degeneration type
of some cover $Y\too X$. The notation and the hypothesis are the
same  as in the previous section.

\begin{defn}\label{def:admissible degeneration type} Any $4$-uple
$(j,\gamma_1,\gamma_2,\kappa)\in\N^4$ with the following properties:
\begin{enumerate}
\item[i)]$\max\{\gamma_1,\gamma_2\} \le \kappa\le v(\lb_{(1)})$;
\item[ii)] $\gamma_2\le p(\kappa -\gamma_1+j)\le
p\gamma_2$; 
\item[iii)] if $\kappa<p\gamma_2$ then
$\gamma_1-j=v(\lb_{(1)})-v(\lb_{(2)})=\frac{v(p)}{p}$; 
if $\kappa\ge p\gamma_2$ then  $0\le p(\gamma_2-j)\le
v(p)-p\gamma_1+\kappa$;
 \item[iv)] $pj\le \gamma_1$;
\end{enumerate}
will be called an \textit{admissible  degeneration type}. 
\end{defn}
\begin{rem}\label{rem:su degeneration type}  We remark that if $\kappa<p \gamma_2$,
then $j$ is uniquely determined from $\gamma_1$ and moreover  i)
and iii) imply iv). The first assertion follows from iii). For the
second, we note that, if $\kappa<p\gamma_2$, multiplying iii) by
$p$, we have $p\gamma_1-pj=(p-1)v(\lb_{(1)})$, since
$pv(\lb_{(2)})=v(\lb_{(1)})$. Therefore, by i), we have
$$\gamma_1-pj=(p-1)(v(\lb_{(1)})-\gamma_1)\ge 0.$$
Moreover, we remark that we always have 
$$0\le \kappa-\gamma_1+j\le \min\{\gamma_2,v(\lb_{(2)})\}.$$ By ii) we have
 to prove only  $ \kappa-\gamma_1+j\le v(\lb_{(2)})$. Moreover, since
$v(\lb_{(1)})\ge \kappa\ge p\gamma_2$ implies $\gamma_2\le
v(\lb_{(2)})$, we have  to consider only the case
$\kappa<p\gamma_2$. But by iii) and i), it follows that
$$
\kappa-\gamma_1+j=\kappa-\frac{v(p)}{p}\le
v(\lb_{(1)})-v(\lb_{(1)})+v(\lb_{(2)})=v(\lb_{(2)}).
$$
\end{rem}
\begin{lem}\label{lem:degeneration type are admissible}Any degeneration type $(j,\gamma_1,\gamma_2,\kappa)$ attached to a $\Z/p^2\Z$-torsor $Y_K\too X_K$
is admissible.
\end{lem}
\begin{proof}
i) 
comes from definitions and \ref{lem:restrictions on gamma2}, 
while iv) has been proved in \ref{lem:restriction degeneration
 type}(i). We now prove ii).
 By \ref{teoremone}, it follows that the effective model of the
action of $\Z/p^2\Z$ on $Y$ is
$\clE^{(\pi^{\kappa},\pi^{\gamma_1};E_p(\alpha_\kappa),1)}$ with
$v(\alpha_\kappa)=\kappa-\gamma_1+j$, if $\alpha_\kappa\neq 0$; and
$\kappa-\gamma_1+j=\gamma_2$ if $\alpha_\kappa=0$.
Since, by \ref{lem:restrictions on gamma2},
$\alpha_\kappa\in \Phi_{\pi^\kappa,\pi^{\gamma_2}}^1$, then  
\begin{equation}
\alpha_\kappa^p\equiv 0\mod\pi^{\gamma_2}.\end{equation}
Hence, we have
$$
\gamma_2\le p(\kappa -\gamma_1+j).
$$
From \ref{lem:restriction degeneration type II}, it follows that
$\kappa
-\gamma_1+j\le \gamma_2$. This proves ii). 

Let us now suppose $\kappa<p\gamma_2$. 
Since $\alpha_\kappa\in \Phi_{\pi^{\kappa},\pi^{\gamma_2}}^1$, by
\ref{cor:clE se lb divide mu} and \ref{rem:valutazione di elementi
di Phi1}, we have that
$$
\kappa-\gamma_1+j=\kappa-\frac{v(p)}{p},
$$
which implies $\gamma_1-j=\frac{v(p)}{p}$.
%
While, if $\kappa\ge p\gamma_2$,  by \ref{cor:clE se lb divide
mu}, we have that
$$
pv(\alpha_\kappa)=p(\kappa-\gamma_1+j)\ge p\gamma_2+(p-1)\kappa-v(p)
$$
which gives
$$
p(\gamma_2-j)\le v(p)-p\gamma_1 +\kappa.
$$
We remark that $\gamma_2-j\ge 0$ comes from \ref{lem:restriction
degeneration type}(ii). Hence iii) is proved.
\end{proof}

\begin{defn}Any admissible degeneration type, which is the degeneration type attached to a $\Z/p^2\Z$-torsor $Y_K\too X_K$,
such that the normalization $Y$ of $X$ in $Y_K$ has integral
special fiber, will be called \textit{realizable}.
\end{defn}
We now see, as a  consequence of  theorem \ref{teoremone}, what
happens in some particular cases. 
\begin{prop}\label{prop:some particulare cases}Let us suppose $Y_K\too X_K$ has $(j,\gamma_1,\gamma_2,\kappa)$ as
degeneration type.
\begin{itemize}
\item[i)] If $j<v(\lb_{(2)})$ then $pj=\gamma_1$ if and only if
$Y$ is a $G_{\pi^j,2}$-torsor. Moreover, the degeneration type is
$(j,pj,j,pj)$. In particular $Y$ is a $\mu_{p^2}$-torsor if and
only if $\gamma_1=0$, i.e. $v(\mathcal{D}(h_1))=v(p)$.
\item[ii)] $j=v(\lb_{(2)})$ if and only if  $Y$ is an
$\clE^{(\pi^{v(\lb_{(1)})},\pi^{\gamma_2};E_p(\eta_{\pi}
S),1)}$-torsor. Necessarily  $\gamma_2\ge v(\lb_{(2)})$ and  the
degeneration type is
$(v(\lb_{(2)}),v(\lb_{(1)}),\gamma_2,v(\lb_{(1)}))$. 

 \item[iii)]  $\gamma_2=j$ if and only if  $Y$ is a
$\clE^{(\pi^{\gamma_1},\pi^{\gamma_2};1,1)}$-torsor. Necessarily
$\gamma_1\ge p\gamma_2$ and the degeneration type is $(
\gamma_2,\gamma_1,\gamma_2,\gamma_1)$. In particular $Y$ is a
$\clE^{(\pi^{\gamma_1},1;1,1)}$-torsor if and only if
$\gamma_2=0$, i.e. $v(\clD(h_2))=v(p)$. \item[iv)] Let $j=0$. Then
$Y_K\too X_K$ is strongly extendible if and only if $\gamma_2=0$.
 \item[v)] $Y$ is a
$\Z/p^2\Z$-torsor if and only if $\gamma_2=v(\lb_{(1)})$, i.e.
$v(\mathcal{D}(h_2))=0$. And the degeneration type is
$(v(\lb_{(2)}),v(\lb_{(1)}),v(\lb_{(1)}),v(\lb_{(1)}))$. 
\item[vi)]If $\gamma_1=v(\lb_{(1)})$, i.e. $v(\clD(h_1))=0$, then
$j=\min\{\gamma_2,v(\lb_{(2)})\}$. So we are in the case $(ii)$ or
$(iii)$.
\end{itemize}
\end{prop}
\begin{rem}The example \ref{ex:non strongly extendible torsor} is in the case $iv)$. 
\end{rem}
\begin{proof}
\noindent By the previous lemma $(j,\gamma_1,\gamma_2,\kappa)$ is
an admissible degeneration type.
\begin{itemize}
\item[i)]
Let us suppose $\gamma_1=pj$. If $\kappa<p\gamma_2$, then by
\ref{def:admissible degeneration type}(iii) it follows that
$$
(p-1)j=\gamma_1-j=\frac{v(p)}{p}=(p-1)v(\lb_{(2)}),
$$
but this is in contradiction with $j<v(\lb_{(2)})$. Hence $\kappa\ge
p\gamma_2$. Therefore, by \ref{def:admissible degeneration
type}(ii),
$$
(p-1)\gamma_2\le \kappa-\gamma_2\le\gamma_1-j=(p-1)j.
$$
But, by \ref{def:admissible degeneration type}(iii), $\gamma_2\ge
j$. 
  Hence,
$\gamma_2=j$. So, by \ref{lem:restriction degeneration type II},
$\kappa=\gamma_1$. Then, by \ref{cor:criterio torsore}, we have
that $Y$ is a $\clE^{(\pi^j,\pi^{pj};1,1)}$-torsor. But, as we
have seen in the example \ref{ex:Gmun},
$$
\clE^{(\pi^j,\pi^{pj};1,1)}\simeq G_{\pi^j,2}
$$
Conversely, as remarked in \ref{rem:degeneration type of G torsors}
$(j,pj,j,pj)$ is the degeneration type of a $G_{\pi^j,2}$-torsor.

 We now observe that, in particular,  $Y$ is a
$\mu_{p^2}$-torsor if and only if $\gamma_1=pj=0$. But, since
$pj\le \gamma_1$ (see \ref{def:admissible degeneration type}(iv)),
then it is true if and only if $\gamma_1=0$, as stated.
\item[ii)]Let us suppose $j=v(\lb_{(2)})$. By \ref{def:admissible
degeneration type}(i),(iv) we have $\gamma_1=pj=v(\lb_{(1)})$ and
$\kappa=v(\lb_{(1)})$. Therefore, by \ref{teoremone}, we have that
$\clE^{(\pi^{v(\lb_{(1)})},\pi^{\gamma_2};E_p(\alpha_{\kappa}S),1)}$
is the effective model. In particular, there is a model map
$$
\Z/p^2\Z\simeq\clE^{(\pi^{v(\lb_{(1)})},\pi^{v(\lb_{(1)})};E_p(\eta_\pi
S),1)}\too
\clE^{(\pi^{v(\lb_{(1)})},\pi^{\gamma_2};E_p(\alpha_{\kappa}S),1)}.
$$
Hence, by \ref{lem:abbasso valutazione con blow-up}, it follows
that
$$\alpha_\kappa\equiv \eta_\pi\mod\pi^{\gamma_2}.$$  So, by \ref{cor:criterio torsore}, $Y$
is a $
\clE^{(\pi^{v(\lb_{(1)})},\pi^{\gamma_2};E_p(\eta_{\pi}S),1)}
$-torsor. Conversely, if $Y$ is a $
\clE^{(\pi^{v(\lb_{(1)})},\pi^{\gamma_2};E_p(\eta_{\pi}S),1)}$-torsor
then, by \ref{rem:degeneration type of G torsors}, the
degeneration type is
$$(v(\eta_\pi),v(\lb_{(1)}),\gamma_2,v(\lb_{(1)})).$$ So
$j=v(\eta_{\pi})=v(\lb_{(2)})$. We observe that 
$j=v(\lb_{(2)})\le \gamma_2$ by \ref{lem:restriction degeneration
type}(ii). 
 \item[iii)]
From \ref{cor:criterio torsore} and \ref{rem:degeneration type of
G torsors}, it follows that $Y$ is a
$\clE^{(\pi^{\gamma_1},\pi^{\gamma_2};1,1)}$-torsor if and only if
$\kappa=\gamma_1$ and $\gamma_2=j$. But, by \ref{lem:restriction
degeneration type II}, $\gamma_2=j$ implies $\kappa=\gamma_1$. The
thesis follows. From \ref{ex:Gmun}, it follows that $\gamma_1\ge
p\gamma_2$. And by \ref{rem:degeneration type of G torsors} we
have that the degeneration type is
$(\gamma_2,\gamma_1,\gamma_2,\gamma_1)$. Now if $\gamma_2=0$ then,
by \ref{lem:restriction degeneration type}, $j=0$ and we have the
last sentence. \item[iv)] If $j=0$, then from \ref{def:admissible
degeneration type}(ii) we have
$$
\gamma_2\le p(\kappa-\gamma_1)\le p\gamma_2.
$$
The thesis easily follows from \ref{cor:criterio torsore}.
 \item[v)] By \ref{rem:degeneration type of  G torsors}, it
 follows that a $\Z/p^2\Z$-torsor $Y\too X$ has
 $$(v(\lb_{(2)}),v(\lb_{(1)}),v(\lb_{(1)}),v(\lb_{(1)}))$$ as degeneration type.
Now let us suppose  $\gamma_2=v(\lb_{(1)})$. Since,  by
\ref{def:admissible degeneration type}$(i)$, $v(\lb_{(1)})\ge
\kappa\ge \gamma_2$ then $\kappa=v(\lb_{(1)})$. Therefore, the
effective model for $Y$ is $\Z/p^2\Z$, since it is a model of
$\Z/p^2\Z$ which is an extension of $\Z/p\Z$ by $\Z/p\Z$ (see
\ref{lem:abbasso
valutazione con blow-up}).  
Let $\sigma$ be a generator of $\Z/p^2\Z$. Since
$\gamma_2=v(\lb_{(1)})$, then,  by \ref{prop:degenerazione Z/pZ
torsori}, $Y\too Y_1$ is a $<\sigma^p>$-torsor. In particular,
$<\sigma^p>$ has no inertia at the generic point of the special
fiber. This implies that $\Z/p^2\Z=<\sigma>$ has no inertia at the
generic point of the special fiber. Let us now consider the
action of $<\sigma>/<\sigma^p>$ on $Y_1=Y/<\sigma^p>$.  
If $\sigma_{|(Y_1)_k}=\id$ then we will have the following
commutative diagram
$$
\xymatrix{Y_k\ar[rd]\ar[rr]&\ar[r]^(.01){\sigma}&\ar[ld]Y_k\\
             &(Y_k)/<\sigma^p>&}
$$
This is a contradiction, since $\sigma_{|Y_k}\neq\id$. So
$<\sigma>/<\sigma^p>$ has no inertia at the generic fiber,
therefore $Y_1\too X$ is a $\Z/p\Z$-torsor, by
\ref{prop:degenerazione Z/pZ torsori}. Hence
$\gamma_1=\kappa=\lb_{(1)}$, which implies, by \ref{cor:criterio
torsore}, that $Y\too X$ is a $\Z/p^2\Z$-torsor.

\item[vi)] Since, by \ref{def:admissible degeneration type}(i), $
v(\lb_{(1)})\ge \kappa\ge \gamma_1$, if $\gamma_1=v(\lb_{(1)})$
then $\kappa=\gamma_1$. So  by \ref{def:admissible degeneration
type}(iii) 
$$
j=v(\lb_{(2)})
$$
if $\gamma_1<p\gamma_2$ (i.e. $\gamma_2>v(\lb_{(2)})$), and
$$
j=\gamma_2
$$
if $\gamma_1\ge p\gamma_2$ (i.e. $\gamma_2\le v(\lb_{(2)})$).

\end{itemize}
\end{proof}

\begin{rem}\label{rem:ammissibile non realizzabile}
Let us suppose that $p|v(p)$. Then, for instance,
$$(j,\frac{v(p)}{p}+j,v(\lb_{(1)}),v(\lb_{(1)}))$$ is admissible
for $0\le j\le v(\lb_{(2)})$, but is not realizable, if $j\neq
v(\lb_{(2)})$, by the point $(v)$ of the proposition.
\end{rem}

So we have  seen that, in general, not all the degeneration types
are realizable. But we now see that it is true for admissible
degeneration types with $\kappa=\gamma_1<v(\lb_{(1)})$. 
They are degeneration types
 attached to $(\Z/p^2\Z)_K$-torsors which are strongly extendible.

\begin{thm}\label{teo:realizzazione tipi di degenerazione} Let us suppose that $R$ contains a primitive $p^2$-th root of
unity and that
 $p>2$. 
Let $X:=\Sp A$ be  a normal  essentially semireflexive scheme
 over $R$  with integral fibers  such that $\pi \in \mathcal{R}_A$. We
 assume \mbox{$_{p^2} Pic(X_K)=0$}, 
$H^1(X_k,\mu_p)\neq 0$ and that there exists $f\in A_k\setminus
{A_k}^p$ such that $A_k^*\not\In A_k^p[f]$. Then any admissible
degeneration type $(j,\gamma_1,\gamma_2,\kappa)$ with
$\kappa=\gamma_1<v(\lb_{(1)})$ is realizable.
\end{thm}

\begin{rem} 

 The hypothesis
$H^1(X_k,\mu_p)\neq 0$ is necessary; otherwise, for instance, it
is not possible to construct $\mu_{p^2}$-torsors  with integral
special fibers. 
Let us now consider a $\clE^{(\gamma_1,1;1,1)}$-torsor with
\mbox{$v(\lb_{(1)})>\gamma_1>0$.} From \S\ref{sec:clE torsors} it
is of the form
$$
Y=\Sp(A[T_1,T_2]/(\frac{(1+\pi^{\gamma_1}T_1)^p-1}{\pi^{p\gamma_1}}-f_1,T_2^p(1+\pi^{\gamma_1}
T_1)^{-1}-f_2))
$$
for some $f_1\in A$ and $f_2\in A^*$. We stress that for any $g\in
A$ we will denote again by ${g}$ its image in $A_k$.  We remark
that $Y_k=\Sp(A_k[T_1,T_2]/(T_1^p-f_1,T_2^p-f_2))$ is integral if
and only if $f_1\not \in {A_k}^p$ and $f_2\in {A_k^*}$ but
$f_2\not \in {A_k^p}[f_1]$. So if the last condition would not be
satisfied, there would not be $\clE^{(\gamma_1,1;1,1)}$-torsors
with integral special fiber. While the hypothesis
$\gamma_1<v(\lb_{(1)})$ is technical and it should be removed
adding
appropriate hypothesis on $X$. 
\end{rem}
\begin{proof}
We recall that in this case to be an admissible degeneration type
means
\begin{enumerate}
\item[i)]$ \gamma_1< v(\lb_{(1)})$; \item[ii)] $\gamma_2\le pj\le
p\gamma_2 $;
 \item[iii)] if $\gamma_1<p\gamma_2$ then
$\gamma_1-j=v(\lb_{(1)})-v(\lb_{(2)})=\frac{v(p)}{p}$; 
if $\gamma_1\ge p\gamma_2$ then $ p(\gamma_2-j)\le
(p-1)(v(\lb_{(1)})-\gamma_1)$;
  \item[iv)] $pj\le \gamma_1$;
\end{enumerate}

We remark that (iv) is in fact implied by the others. 
Indeed, let us suppose that $pj>\gamma_1$. Then by (ii) we have
$p\gamma_2\ge pj>\gamma_1$. But we know,  by \ref{rem:su
degeneration type} that if $p\gamma_2>\gamma_1$ then $pj\le
\gamma_1.$


Since $\kappa=\gamma_1$, it follows, by \ref{rem:degeneration type
of G torsors}, that if $(j,\gamma_1,\gamma_2,\kappa)$ is
realizable it is the degeneration type of a
$\clE^{(\gamma_1,\gamma_2;E_p(\alpha_{\gamma_1}S),1)}$-torsor,
with $v(\alpha_{\gamma_1})=j$ if $\alpha_{\gamma_1}\neq 0$ and
$j=\gamma_2$ if $\alpha_{\gamma_1}=0$.
For any $\gamma_1,\gamma_2$ as in the degeneration type, by
\ref{cor:clE se lb divide mu} 
there exists  a group scheme
$\clE^{(\pi^{\gamma_1},\pi^{\gamma_2};E_p(a S),1)}$. If $a\neq 0$
then we can choose $a$ such that  $v(\tilde{a})=j$, where
$\tilde{a}\in R$ is a lifting of $a$. In fact, if
$\gamma_1<p\gamma_2$ it is automatic, by \ref{rem:valutazione di
elementi di Phi1} and iii), that $v(\tilde{a})=j$. We call
$a=\alpha_{\gamma_1}$.
We now construct a normal $
\clE^{(\gamma_1,\gamma_2;E_p(\alpha_{\gamma_1} S),1)} $-torsor
with integral special fiber. First of all, we remark that if
$\gamma_1=0$ then $\gamma_2=0$ and
\mbox{$\clE^{(\gamma_1,\gamma_2;E_p(\alpha_{\gamma_1} S),1)}\simeq
\mu_{p^2} $.} So if we take $Y=\Sp(A[T]/(T^{p^2}-f))$ with $f\in
A^*$ not a $p$-th power $\mod\pi$, then $Y_k$ is integral. Hence
$Y$ is normal by \ref{lem:normality's criterion}.

We now suppose $\gamma_1>0$. 
 Let us take  $f_1,f_2\in A$.
Since $\pi\in \mathcal{R}_A$, then $1+\pi^{p\gamma_1} f_1\in A^*$
and \mbox{$E_p(\alpha_{\gamma_1}^p f_1)+\pi^{p\gamma_2}f_2\in
A^*$.} Then we can define  the
$\clE^{(\gamma_1,\gamma_2;E_p(\alpha_{\gamma_1} S),1)} $-torsor
$$
Y=\Sp \bigg(A[T_1,T_2]/\big(\frac{(1+\pi^{\gamma_1}
T_1)^p-1}{\pi^{p\gamma_1}}-f_1,\frac{(E_p(\alpha_{\gamma_1}
T_1)+\pi^{\gamma_2} T_2)^p(1+\pi^{\gamma_1}
T_1)^{-1}-E_p(\alpha_{\gamma_1}^p
f_1)}{\pi^{p\gamma_2}}-f_2\big)\bigg). $$ See \S\ref{sec:clE
torsors} and  \ref{rem:G=E(ap S)}.
 We have only to find
$f_1$ and $f_2$ such that $Y$ has integral special fiber.
Take $f_1$ such that $f_1$ is not a $p$-{th} power $\mod \pi$ and
$A_k^*\not\In A_k^p[f_1]$. Then we have that the special fiber of
$$
Y_1=\Sp(B_1)=\Sp \bigg(A[T_1]/(\frac{(1+\pi^{\gamma_1}
T_1)^p-1}{\pi^{p\gamma_1}}-f_1)\bigg)
$$
is integral. We now consider
$$
Y=\Sp\bigg(B_1[T_2]/(\frac{(E_p(\alpha_{\gamma_1}
T_1)+\pi^{\gamma_2} T_2)^p(1+\pi^{\gamma_1}
T_1)^{-1}-E_p(\alpha_{\gamma_1}^p
f_1)}{\pi^{p\gamma_2}}-f_2)\bigg).
$$
Since $\gamma_2\le \gamma_1$ then 
 $\gamma_2<v(\lb_{(1)})$. The special fiber is
$$
Y_k=\Sp\bigg((B_1)_k[T_2]/\big(T_2^p- \frac{E_p(\alpha_{\gamma_1}
T_1)^{-p}(1+\pi^{\gamma_1}T_1) E_p(\alpha_{\gamma_1}^p
f_1)-1}{\pi^{p\gamma_2}}-f_2 \big)\bigg)
$$
Let $P(T_1):=\frac{E_p(\alpha_{\gamma_1}
T_1)^{-p}(1+\pi^{\gamma_1}T_1) E_p(\alpha_{\gamma_1}^p
f_1)-1}{\pi^{p\gamma_2}}$. If $P(T_1)+f_2$ is not a $p$-{th} power
$\mod \pi B_1$,  then $Y_k$ is reduced. While if
 $P(T_1)+f_2\equiv P_1(T_1)^p\mod \pi $ for some $P_1(T_1)\in B_1^*$ then we substitute
$f_2+f_3$ to $f_2$
 with $f_3$ not a $p$-{th} power $\mod\pi B_1$. The fact that $f_3$ is not a $p$-th power $\mod \pi B_1$ means, as one can easily check, that
$f_3 \mod \pi$ does not belong to $A_k^p[f_1]$. Such $f_3$ there
exists by hypothesis on $f_1$. Now, if $$
 P(T_1)+f_2+f_3\equiv P_2(T_1)^p\mod
 \pi B_1
 $$ for some $P_2(T_1)\in B_1^*$, then
 $$
f_3\equiv (P_2(T_1)-P_1(T_1))^p\mod \pi B_1.
 $$
Then it follows that $f_3$ is a $p$-power $\mod \pi B_1$, against
hypothesis on $f_3$.


Finally, we verify that $Y$ has $(j,\gamma_1,\gamma_2,\gamma_1)$
as degeneration type. Since $\kappa=\gamma_1$, by \ref{teoremone},
we have that the degeneration  type is
$(v(\alpha_{\gamma_1}),\gamma_1,\gamma_2,\gamma_1)$ if
$\alpha_{\gamma_1}\neq 0$ and
$(\gamma_2,\gamma_1,\gamma_2,\gamma_1)$ if $\alpha_{\gamma_1}=0$.
But, since we have chosen $\alpha_{\gamma_1}$ such that
$\alpha_{\gamma_1}=0$ and $j=\gamma_2$ or $v(\alpha_{\gamma_1})=j$
and $\alpha_{\gamma_1}\neq 0$, then we have the thesis.

\end{proof}

\end{document}